\documentclass{article}
\usepackage[margin=2cm]{geometry}

\usepackage{import}
\usepackage{setspace}
\spacing{1.5}

\usepackage{amsmath, amssymb, bbold, mathtools}
\usepackage[ruled, vlined]{algorithm2e}
\usepackage{multirow}
\usepackage{graphicx}
\usepackage{tcolorbox}
\usepackage[]{hyperref}
\usepackage{bbding}
\usepackage[export]{adjustbox}
\usepackage{subcaption}
\usepackage{enumitem}

\usepackage{lmodern} 

\usepackage{natbib}
 \bibpunct[, ]{(}{)}{,}{a}{}{,}

\usepackage{tcolorbox}
\usepackage{tabularx}
\usepackage{tabularray}

\newcommand{\keywords}[1]{\textbf{Keywords:} #1}

\usepackage{tikz}
\usetikzlibrary{calc}
\usetikzlibrary{shapes}

\usepackage{amsthm}
\theoremstyle{plain}
\newtheorem{theorem}{Theorem}[section]
\newtheorem{lemma}[theorem]{Lemma}

\newtheorem{corollary}[theorem]{Corollary}

\newtheorem{example}[theorem]{Example}
\newtheorem{proposition}[theorem]{Proposition}

\newtheorem{definition}[theorem]{Definition}

\usepackage{acronym}
\acrodef{wlog}[WLOG]{without loss of generality}
\acrodef{lsc}[lsc]{lower semi-continuous}
\acrodef{DRO}[DRO]{distributionally robust optimization}
\acrodef{KL}[KL]{Kullback-Leibler}
\acrodef{LP}[LP]{L\'evy-Prokhorov}

\definecolor{red}{RGB}{163, 31, 52}
\definecolor{gray}{RGB}{194, 192, 191}
\definecolor{blue}{RGB}{59, 89, 152}
\definecolor{green}{RGB}{0, 179, 0}

\newcommand{\norm}[1]{\left\|#1\right\|}

\newcommand{\Prob}[1]{{\mathrm{Prob}}\left[ #1 \right]}
\newcommand{\abs}[1]{\left|#1\right|}
\newcommand{\E}[2]{{\mathbb E}_{#1} \left[ #2 \right]}

\newcommand{\one}[1]{\mathbb{1}\left\{#1\right\}}
\renewcommand{\tfrac}[2]{{#1}/{#2}}

\newcommand{\set}[2]{\left\{ #1\ : \ #2 \right\}}
\newcommand{\tset}[2]{\{ #1\ : \ #2 \}}
\newcommand{\LP}[2]{{\rm{LP}}(#1,#2)}

\newcommand{\defn}[0]{:=}

\newcommand{\mc}{\mathcal}

\newcommand{\mr}{\mathrm}

\renewcommand{\d}{{\mathrm{d}}}

\renewcommand{\Re}{\mathbb{R}}
\newcommand{\nRe}{\Re\cup\{+\infty\}}

\newcommand{\st}{\mr{s.t.}}

\renewcommand{\emph}{\textbf}

\DeclareMathOperator{\esssup}{{\rm{ess}}\sup}

\DeclareMathOperator{\kl}{KL}

\DeclareMathOperator{\lk}{LK}
\DeclareMathOperator{\lp}{LP}
\DeclareMathOperator{\dro}{DRO}
\DeclareMathOperator{\diam}{diam}

\DeclareMathOperator{\dom}{dom}

\DeclareMathOperator{\diag}{diag}

\title{Smoothed $f$-Divergence Distributionally Robust Optimization}
\author{Zhenyuan Liu \and Bart P.G.\ Van Parys \and Henry Lam}

\begin{document}

\maketitle

\begin{abstract}
 In data-driven optimization, sample average approximation (SAA) is known to suffer from the so-called optimizer's curse that causes an over-optimistic evaluation of the solution performance. 
 We argue that a special type of distributionallly robust optimization (DRO) formulation offers theoretical advantages in correcting for this optimizer's curse compared to simple ``margin'' adjustments to SAA and other DRO approaches: It attains a statistical bound on the out-of-sample performance, for a wide class of objective functions and distributions, that is nearly tightest in terms of exponential decay rate.
 This DRO uses an ambiguity set based on a Kullback Leibler (KL) divergence smoothed by the Wasserstein or L\'evy-Prokhorov (LP) distance via a suitable distance optimization. Computationally, we also show that such a DRO, and its generalized versions using smoothed $f$-divergence, are not harder than DRO problems based on $f$-divergence or Wasserstein distances, rendering our DRO formulations both statistically optimal and computationally viable.
\end{abstract}





\section{Introduction}
In this paper we discuss data-driven stochastic optimization problems in the following generic form.
Let $X\subseteq \Re^m$ be a set of feasible decisions and $\xi$ an uncertain parameter or random variable realizing in a compact event set $\Sigma\subseteq \Re^d$ with distribution $P$.
The loss incurred for a given decision $x\in X$ is $\ell(x, \xi)$. A decision with minimal expected loss, or cost $c(x, P) \defn \E{P}{\ell(x, \xi)}$, can be found as the solution to the stochastic minimization problem
\begin{equation}
  \label{eq:so}
  \textstyle\inf_{x\in X} c(x, P).
\end{equation}
Such formulations come about naturally in a wide variety of applications ranging from risk minimization in machine learning to portfolio selection in finance.
We will assume here that $\ell(x, \xi)$ is bounded below and thus the expectation $\E{P}{\ell(x, \xi)}$ is always well defined. We further assume that $\E{P}{\ell(x, \xi)}<\infty$ for some $x\in X$ and hence the considered stochastic optimization problem is feasible.

We consider the data-driven situation where the distribution $P$ is unknown but rather learned from i.i.d.\ historical observations $\{\xi_1,\ldots,\xi_N\}$. Summarizing these observations as the empirical distribution $P_N \defn \sum_{i=1}^N \delta_{\xi_i}/N$, a data-driven optimization formulation takes the form $\inf_{x\in X}\hat c(x,P_N)$, where $\hat c(\cdot,P_N)$ is a proxy or \textit{predictor} of $c(\cdot,P)$. Our main quantity of study in this paper is the probability
\begin{equation}
\Prob{\exists x\in X~\st~c(x,P)>\hat c(x, P_N)}\label{beginning prob}
\end{equation}
which represents the chance that a predictor $\hat c(x,P_N)$ cannot upper bound the correct objective value $c(x,P)$ for at least part of $X$. This probability is in general intricate depending on the form of the loss function $\ell(\cdot,\cdot)$ and the sample size $N$, and has spurred much interest in the machine learning and statistical optimization literature (as we explain next). On a high level, our goal in this paper is to show that, as the sample size grows, there is a data-driven formulation represented by $\hat c$ that controls \eqref{beginning prob} with good theoretical properties, for a wide class of loss functions $\ell$ and ground-truth distribution $P$. Moreover, this data-driven formulation is a special type of distributionally robust optimization (DRO) using what we call smoothed Kullback-Leibler (KL) divergence.

In the following, let us first motivate our study of quantity \eqref{beginning prob} in several stages.

\subsection{Beginning Motivation: Optimizer's Curse}
As a starting point, we discuss empirical optimization or sample average approximation (SAA) \cite{shapiro2021lectures,birge2011introduction,royset2022optimization}, arguably the most direct data-driven approach. This approach replaces the unknown $P$ by $P_N$ in \eqref{eq:so} to obtain the solution of
\begin{equation}
  \label{eq:saa}
  \textstyle\inf_{x\in X} c(x, P_N).
\end{equation}
Because of the luck of the draw, the obtained optimal value in \eqref{eq:saa} differs, and in fact is smaller in expectation, than the true optimal value in \eqref{eq:so}, thus giving an overly optimistic impression to the optimizer \citep{mak1999monte}. This optimistic bias phenomenon is colloquially known as the \textit{optimizer's curse}. Note that this also implies an optimistic impression at the \textit{solution level}. Denote $\hat x^{SAA}$ as the SAA solution obtained from \eqref{eq:saa}. Then, the true performance of $\hat x^{SAA}$ is typically worse than the optimal value obtained in \eqref{eq:saa}, because
\begin{equation}
\Prob{\textstyle\inf_{x\in X}c(x,P)>\inf_{x\in X}c(x, P_N)}\leq\Prob{c(\hat x^{SAA},P)>c(\hat x^{SAA}, P_N)}\label{prob saa1}
\end{equation}
where the LHS indicates the probability of being overly optimistic with respect to the optimal value, while the RHS denotes the chance of being overly optimistic with respect to the considered solution $\hat x^{SAA}$. The RHS, in particular, 
is closely related to \textit{overfitting} in machine learning. In the latter setting, $c(x,P)$ is called the risk or expected loss and problem \eqref{eq:saa} corresponds to empirical risk minimization. The optimal value $c(\hat x^{SAA}, P_N)$ obtained by SAA is denoted as the training loss, while its testing loss is taken precisely as $c(\hat x^{SAA},P)$. The RHS in \eqref{prob saa1} thus quantifies the generalization or out-of-sample performance of SAA in terms of the probability that its testing loss is larger than the training loss. 

Note that the events in \eqref{prob saa1} involve random functions evaluated at random solutions. To bound these probabilities, a simple and popular approach is to use the uniform deviation probability \eqref{beginning prob}, where the predictor $\hat c$ considered here is $c$, i.e.,
\begin{equation}
\Prob{c(\hat x^{SAA},P)>c(\hat x^{SAA}, P_N)}\leq\Prob{\exists x\in X~\st~c(x,P)>c(x, P_N)}.\label{bounding relation}
\end{equation}
In other words, \eqref{beginning prob} provides an upper bound for the probabilities of being overly optimistic both at the value level and solution level. Note also that, because of the optimizer's curse, the probability in either side of \eqref{bounding relation} is non-vanishing as the sample size $N$ grows. Because of this, typically one considers the probabilities
\begin{equation}
\Prob{c(\hat x^{SAA},P)>c(\hat x^{SAA}, P_N)+\delta}\leq\Prob{\exists x\in X~\st~c(x,P)>c(x, P_N)+\delta}\label{bounding relation1}
\end{equation}
for a small $\delta>0$. The RHS probability in \eqref{bounding relation1} in particular characterizes how likely the event occurs in which the predicted cost underestimates the unknown actual cost of the decision $\hat x^{SAA}$ by at least some positive margin $\delta$. 
However, the RHS probability in \eqref{bounding relation1}, which can be rewritten as
\begin{equation}
\Prob{\textstyle\sup_{x\in X}\{c(x,P)-c(x, P_N)\}>\delta},\label{EP}
\end{equation}
contains the supremum of an empirical process $\sup_{x\in X}\{c(x,P)-c(x, P_N)\}$. Bounding \eqref{EP} thus requires measurements of the \textit{complexity} of the loss function class $\{\ell(x,\cdot):x\in X\}$ such as the Vapnik-Chervonenkis (VC) dimension or other variants \cite[Chapters 4, 5]{wainwright2019high}. More concretely, these bounds typically take the form
\begin{equation}
\Prob{\textstyle\sup_{x\in X}\{c(x,P)-c(x, P_N)\}>\delta}\leq g(N,\mathcal C)\exp(-\delta^2 N/\mathcal V)\label{ERM bound}
\end{equation}
where $\mathcal C$ is a complexity measure of the function class $\{\ell(x,\cdot):x\in X\}$, $g(N,\mathcal C)$ is a quantity that depends subexponentially in $N$, and $\mathcal V$ is some variability measure of $\ell$ such as its variance or support bound.

\subsection{Exponential Decay Rate and Limitations of SAA}\label{sec:complexity}

Our main focus in this paper is to devise good predictors such that the RHS of \eqref{bounding relation1} decays gracefully as $N\to\infty$. In particular, we take a large deviations viewpoint in that we seek to attain a prescribed exponential decay rate in $N$ of the RHS of \eqref{bounding relation1}. In the case of SAA discussed above for instance, to attain an exponential rate $r$, we can use $\inf_{x\in X}c(x,P_N)+\delta$ as the predicted optimal value, where $\delta$ is suitably calibrated in terms of sample size $N$ and $\mathcal V$. This approach, while straightforward, has two limitations. First is that \eqref{ERM bound} can be a loose bound for any particular problem instance, so that the margin adjustment $\delta$ can be larger than needed. This is ignoring the fact that $\mathcal V$ might need to be estimated and this estimate is also loose. Second is that the bound \eqref{ERM bound} only holds for loss function classes with finite $\mathcal C$, which can be difficult to estimate or even determine for finiteness. Given this, our main goal is to provide an alternative approach based on DRO such that, \textit{as far as the exponential decay rate is concerned, could achieve nearly tightest bounds and, moreover, using minimal information on both the loss function class and underlying distribution.}

\subsection{Literature Review on DRO}\label{sec:DRO literature}
We review the DRO literature that connects to the focus of this paper. DRO is a paradigm for decision-making under uncertainty that, instead of simply ``plugging in'' the empirical distribution as in SAA, DRO advocates a worst-case viewpoint. More precisely, when facing unknown or uncertain parameters in decision-making, DRO solves for a solution that optimizes under the worst-case scenario, over a so-called ambiguity set or uncertainty set that, at least intuitively, contains the true parameter with high confidence. This typically results in a min-max formulation, where the outer minimization is over the decision and the inner maximization over the ambiguity set. As such, it can be viewed as a generalization of classical robust optimization for stochastic problems where the uncertain parameter is the underlying probability distribution. DRO originates in \cite{scarf1958min} and since then has been applied in various disciplines across control \citep{petersen2000minimax}, economics \citep{hansen2008robustness}, finance \citep{wozabal2012framework,wang2013bounds,glasserman2014robust,yuen2020distributionally,dhara2021worst}, queueing \citep{jain2010optimality,van2022distributionally} and revenue management \citep{rusmevichientong2012robust}, and recently experiences a fast-surging interest in machine learning \citep{chen2018robust,kuhn2019wasserstein,blanchet2019robust}, partly because of reasons that we will discuss next.

To be concrete, we write the DRO formulation as
\begin{equation}
\textstyle\inf_{x\in X}\sup_{Q\in\mathcal A_N}c(x,Q)\label{DRO general}
\end{equation}
where $\mathcal A_N$ is the ambiguity set. In the data-driven context, this set is calibrated from data and hence has the dependence on the sample size $N$. There are two main approaches to construct $\mathcal A_N$ in the DRO literature, one via moments, support and shape constraints \citep{delage2010distributionally,bertsimas2005optimal,wiesemann2014distributionally,goh2010distributionally,ghaoui2003worst,van2016generalized,li2019ambiguous,chen2021discrete,bai2023distributionally} that appear especially useful in applications with only partial observability and when dealing with extreme events \citep{lam2017tail,lam2021orthounimodal}. However, this approach does not guarantee statistical consistency of the obtained solution, i.e., convergence to the ground-truth solution as the sample size grows, because these constraints  intrinsically discard some information about the underlying distribution. The second approach takes $\mathcal A_N$ as a neighborhood ball surrounding some baseline distribution (often the empirical distribution $P_N$) where the size of the ball is characterized by some statistical distance, i.e., $\mathcal A_N=\{Q:D(Q,P_N)\leq\eta\}$ for some distance $D$ and ball size $\eta$. The most popular choices of distances are $f$-divergences \citep{gupta2019near,bayraksan2015data,Iyengar2005,hu2013kullback,gotoh2018robust,Lam2016robust,Lam2018sensitivity,ghosh2019robust} and the Wasserstein distance \citep{mohajerin2018data,chen2018robust,xie2019tractable,shafieezadeh2019regularization,blanchet2019quantifying,blanchet2021sample,gao2022distributionally}, with growing interest in other choices such as maximum mean discrepancy \citep{staib2019distributionally,zeng2022generalization}, integral probability metrics \citep{husain2020distributional,zhu2021kernel}, and further variants. Unlike the first approach, DROs using neighborhood balls are often statistically consistent, and for this reason constitute a reasonable competitor to SAA.
Circling back to our main question discussed at the end of Section \ref{sec:complexity}, let us now review the two main paths to obtain statistical guarantees for DRO.

\paragraph{Covering DRO}
The most basic approach is to construct an ambiguity set $\mathcal A_N$ that covers the ground-truth $P$ with high confidence, i.e., $\mathcal A_N$ is a confidence region for $P$ in the sense $\Prob{P\in\mathcal A_N}\geq1-\alpha$ for some level $1-\alpha$. This then translates into a confidence bound guarantee for the true performance of the DRO solution, thus providing an alternative bound to \eqref{bounding relation} \cite{bertsimas2018robust,delage2010distributionally}. More precisely, denoting $\hat x^{DRO}$ as an optimal solution to \eqref{DRO general}, the condition $\Prob{P\in\mathcal A_N}\geq1-\alpha$ implies immediately
\begin{equation}
\Prob{c(\hat x^{DRO},P)\leq \textstyle\sup_{Q\in\mathcal A_N}c(\hat x^{DRO},Q)}\geq1-\alpha\label{DRO guarantee general}
\end{equation}
since $P\in\mathcal A_N$ implies $c(\hat x^{DRO},P)\leq\sup \tset{c(\hat x^{DRO},Q)}{{Q\in\mathcal A_N}}$ by the worst-case definition in constructing the DRO objective. While straightforward, this reasoning could be restrictive or dimension dependent when applying to particular choices of $\mathcal A_N$. For instance, when using $f$-divergences, the confidence region only applies to ground-truth distributions that are discrete or otherwise the ball must contain all probability distributions in order to contain the ground-truth distribution \cite[Chapter 7]{polyanskiy2023information}. This implies that the statistical guarantee in \eqref{DRO guarantee general} with a non-trivial $\mathcal{A}_N$ only applies to discrete distributions. When using the $1$-Wasserstein distance $W(\cdot,\cdot)$, we have for $d>2$ and any number of samples $N\geq 1$ that
\begin{equation}
  \label{eq: Wasserstein bound}
  \Prob{W(P,P_N)> \eta} \leq C \exp\left(-c N \eta^{d}\right) \text{ if } \eta\leq 1,
\end{equation}
where $C>0$ and $c>0$ are known constants which depend only on the dimension $d$ and the true distribution $P$ \cite[Theorem 2]{fournier2015rate}. Consequently, to ensure $P\in\mathcal{A}_N=\set{Q}{W(Q,P_N)\leq \eta}$ and thus obtain guarantee \eqref{DRO guarantee general} with a prescribed fixed $\alpha$, we need the ball size to scale with the distribution dimension $d$, i.e., $\eta\geq(\log(C/\alpha)/(cN))^{1/d}$. Equivalently, if we need to attain a prescribed exponential decay rate $r$, i.e.,
\begin{equation}
  \label{eq:oos-guar}
  \textstyle\lim_{N\to\infty}\frac 1N \log \Prob{c(\hat x^{DRO},P)>\sup_{Q\in\mathcal A_N}c(\hat x^{DRO},Q)}\leq -r,
\end{equation}
we need the ball size to be $\eta\geq(r/c)^{1/d}$.
This curse of dimensionality seems unavoidable as $\E{}{W(P,P_N)}$ tends to zero not faster than $N^{-1/d}$ for any measure $P$ absolutely continuous with respect to the Lebesgue measure \citep{dudley1969speed, weed2019sharp}.
Thus, the classical rationale of ambiguity set coverage appears to run into either distribution restriction or distribution dimension dependence. We should make clear that the latter does not \textit{a priori} mean that the approach is conservative, but the conservativeness arguably appears when compared to the empirical DRO perspective that we discuss next.

\paragraph{Empirical DRO}
One can observe that in order to provide confidence bounds for $c(\hat x^{DRO},P)$, it is not necessary to enforce distributional coverage of the ambiguity set. That is, while constructing an ambiguity set to cover the ground-truth distribution does give rise to the confidence bound guarantee \eqref{DRO guarantee general}, the obtained bound is loose which means the suggested ball is chosen too large. In contrast, empirical DRO views robustness as a form of variability regularization \citep{duchi2018variance,gotoh2018robust,shafieezadeh2019regularization,gao2022wasserstein}, and bears a duality relation with the so-called empirical likelihood method in statistics \citep{lam2017empirical,lam2019recovering,duchi2021statistics}.  In particular, when the ball size $\eta$ of the ambiguity set surrounding the empirical distribution is small, roughly speaking 
\begin{equation}
\textstyle\sup_{Q\in\mathcal A_N}c(x,Q)\approx c(x,P_N)+V\sqrt{\eta}+\cdots\label{var reg}
\end{equation}
where $V$ is a variation measure. For example, $\mc V=\sqrt{\mathrm{Var}_P[\ell(x,\xi)]}$ (multiplied by a constant depending on $f$) for $f$-divergences \citep{Lam2016robust,duchi2018variance,duchi2021statistics,gotoh2018robust}, or $\mc V$ is the Lipschitz norm of $\ell(x,\cdot)$ for the 1-Wasserstein distance \citep{shafieezadeh2019regularization,blanchet2019robust,gao2022wasserstein,blanchet2022confidence}. This view allows us to deduce  that
$$\textstyle\Prob{c(\hat x^{DRO},P)>\sup_{Q\in\mathcal A_N}c(\hat x^{DRO},Q)}$$
and subsequently conclude the strengths of DRO when the loss function has low variability \citep{duchi2018variance}. This deduction arises through \eqref{var reg} and yields the parametric rate $O(1/\sqrt N)$ of $\sup_{Q\in \mathcal A_N} c(\hat x^{DRO},Q)-c(\hat x^{DRO}, P_N)$ regardless of the support of the true distribution in the case of $f$-divergence \citep{duchi2021statistics,lam2019recovering,gupta2019near} in lieu of nonparamteric rate $O(N^{-1/d})$ in the case of covering Wasserstein DRO \citep{blanchet2019robust,gao2022finite}.
However, as here the ambiguity set does not cover the unknown distribution, it typically requires conditions on the complexity of the loss function class analogous to the SAA paradigm \citep{duchi2021statistics}. 
 
\section{Main Contributions}\label{sec:contributions}
In contrast to the above two DRO paths, our work is most related to the recent works \citep{vanparys2021data} and \citep{bennouna2021learning} in which a large deviations perspective is considered. In particular, \citep{vanparys2021data} and \citep{bennouna2021learning} show that, in achieving a given exponential decay rate $r$ in \eqref{eq:oos-guar}, DRO using the KL-divergence attains the tightest possible bound in a ``meta-optimization'' sense. However, these results are argued only in the case of discrete distributions where the KL-divergence naturally applies (a straightforward extension to continuous distributions is not trivial \citep{vanparys2021data}) and, in this situation, the uniform bound \eqref{beginning prob} is not substantially different from a fixed-decision bound. Compared to these works, our paper can be viewed as a substantial generalization in understanding the uniform bound \eqref{beginning prob} for a much wider class of distributions and general loss functions. Moreover, to achieve this generalization, a main contribution of this paper is to introduce what we call a \textit{smoothed KL-divergence} that suitably combines KL with either the Wasserstein or LP distance.

To make our discussion concrete, let $\mc P(\Sigma)$ denote the space of probability measures on $\Sigma$. We call a cost predictor $\hat c:X\times \mc P(\Sigma)\to\Re$ \textit{uniformly feasible} at rate $r>0$ if
\begin{equation}
  \label{eq:feasibility-uniform}
  \textstyle\sup_{P\in \mc P(\Sigma)}~\limsup_{N\to\infty}\frac{1}{N}\log \Prob{\exists x\in X~\st~ c(x,P)>\hat c(x, P_N)} \leq -r.
\end{equation}
In other words, the cost prediction $\hat c(x, P_N)$ constitutes an upper bound on the actual unknown cost $c(x,P)$ uniformly for all $x\in X$, with probability at least $1-\exp(-r N + o(N))$ whatever distribution $P\in \mc P(\Sigma)$ generates the data. A closely related but much weaker requirement is \textit{pointwise feasibility}. Namely, we call a cost predictor $\hat c:X\times \mc P(\Sigma)\to\Re$ pointwise feasible at rate $r>0$ if
\begin{equation}
  \label{eq:feasibility}
  \textstyle\sup_{P\in \mc P(\Sigma), \, x\in X}~\limsup_{N\to\infty}\frac{1}{N}\log \Prob{c(x,P)>\hat c(x, P_N)} \leq -r.
\end{equation}
Compared to \eqref{eq:feasibility-uniform}, the requirement for all $x\in X$ is now outside the probability in \eqref{eq:feasibility} and, as such, \eqref{eq:feasibility} does not bound the probability of committing the optimizer's curse that we discussed previously. The normalized sequence $\{N\}$ in (\ref{eq:feasibility-uniform}) and (\ref{eq:feasibility}) essentially imposes a minimal speed with which the probability of the optimistic bias tends to zero. While these expressions resemble classical large deviations analyses \citep{dembo2009large}, uniform feasibility in (\ref{eq:feasibility-uniform}) asserts the strong property that the deviation holds for all decision $x$, which as we discussed is our main focus and is beyond this classical literature.

With the above rate-feasibility definitions, we can introduce our concept of efficiency. We call a cost predictor $\hat c(\cdot,\cdot)$ to be \textit{least conservative} among a class of predictors $\mathcal C$ in achieving \eqref{eq:feasibility} or \eqref{eq:feasibility-uniform} if, for any $\hat c'(\cdot,\cdot)\in\mathcal C$ that satisfies \eqref{eq:feasibility} or \eqref{eq:feasibility-uniform} at rate $r>0$, we must have $\hat c(x,\hat P)\leq\hat{c}'(x,\hat P)$ for any $x\in X,\hat P\in\mathcal{P}(\Sigma)$. Note that our least conservative notion is strong in the sense that we require $\hat c(x,\hat P)\leq\hat{c}'(x,\hat P)$ to hold for \textit{any} $x\in X,\hat P\in\mathcal{P}(\Sigma)$ instead of only \textit{some} of $x\in X,\hat P\in\mathcal{P}(\Sigma)$. As it turns out, this strong notion is appropriate for our analysis as our proposed cost predictor would indeed satisfy this strong efficiency notion.

Now we are ready to state the main contributions and implications of this paper as follows. 

\begin{itemize}[leftmargin=1em]
\item \textbf{Least Conservative Uniform Feasibility of KL-DRO in the Finite Support Case ($\abs{\Sigma}<\infty$):} We show that the KL-DRO predictor is the least conservative predictor among all upper semicontinuous cost predictors $\hat c$ satisfying uniform feasibility (\ref{eq:feasibility-uniform}). Moreover, along with this result, we also show that in terms of the DRO ambiguity set, the KL-divergence ambiguity set is the smallest among all uniformly feasible closed convex ambiguity sets. As a comparison, \cite{vanparys2021data} only considers least conservativeness for pointwise feasibility and also does not consider optimality in terms of the ambiguity set.

\item \textbf{Infinite Support Case ($\abs{\Sigma}=\infty$):} The infinite support case is substantially more challenging than the finite support case because the KL-divergence does not naturally apply. In this regard, our main contribution is a smoothing enhancement of the KL-divergence to bypass this applicability issue and attain least conservative uniform feasibility. More concretely, our contribution is four-fold: 

\begin{itemize}
    \item \textbf{Lack of Uniform Feasibility of KL-DRO:} We prove through a counterexample that the KL-DRO predictor may not be uniformly feasible even if the loss function is assumed to be jointly continuous in $x$ and $\xi$. Interestingly, this counterexample relies on the concept of ``uniformly distributed modulo 1'' \citep[Chapter 1 Section 6]{kuipers2012uniform} in analytic number theory and uses a generalization of the denseness property of $\{\sin{n},n\in\mathbb{N}\}$ on $[-1,1]$.
    \item \textbf{Uniform Feasibility of Smoothed KL-DRO:} We create two smoothed KL-divergences to resolve the feasibility issue, one via the injection of the Wasserstein distance and one via the LP distance. We show that their associated ambiguity sets enjoy finite-sample large deviations guarantees similar to those of the KL-divergence in the finite support case. In particular, regardless of the smoothing parameter, these smoothed KL-DRO predictors are uniformly feasible with minimal assumptions on the loss function beyond those needed to establish measurability of the events in (\ref{eq:feasibility-uniform}) and (\ref{eq:feasibility}).
    \item \textbf{Asymptotic Least Conservativeness of Smoothed KL-DRO:} We show that our smoothed KL-DRO predictors are the least conservative, asymptotically in the sense that the smoothing parameter goes to zero, among all upper semicontinuous cost predictors $\hat c$ satisfying uniform feasibility (\ref{eq:feasibility-uniform}). Additionally, in terms of the DRO ambiguity set, we show that our smoothed KL-divergence ambiguity sets are the smallest, again asymptotically as the smoothing parameter goes to zero, among all uniformly feasible closed convex ambiguity sets. In other words, our smoothed KL-DRO and associated ambiguity sets are the least conservative up to this smoothing parameter.

    \item \textbf{Exponential-Rate Efficiency Compared to SAA:} We contrast the above two contributions with SAA that, even to construct a predictor that attains uniform feasibility, requires information on the loss function complexity and other variability parameters. In contrast, our smoothed DRO predictors are ``automatic'' in that no such parameters are needed. Moreover, up to the smoothing parameter, our predictors are the least conservative while there is no optimality result of this sort in SAA due to the looseness of their associated concentration bounds.
\end{itemize}

In addition to the above results on least conservative uniform feasibility and their statistical implications, we also make two contributions regarding our smoothed DRO, not only for KL but more general $f$-divergences.

\item \textbf{Relative Efficiency Quantification:} We propose a quantification scheme to characterize the loss of statistical performance of any DRO formulation relative to our smoothed KL-DRO predictors. Our loss is characterized as the Pompeiu-Hausdorff distance \cite[Example 4.13]{rockafellar2009variational} measured by the Wasserstein distance between the ambiguity sets of the compared formulations. 

\item \textbf{Computation:} We show that our smoothed divergence DRO predictors are computationally not much harder than standard divergence DRO and Wasserstein DRO problems. Our smoothed $f$-divergence DRO predictor admits a tractable reformulation in terms of an inflated loss function \citep{conn2000trust} and the convex conjugate $f^*$ and can be computed to any desired level of relative accuracy as long as this inflated loss function can be efficiently evaluated. 
  
\end{itemize}

\section{Smoothed $f$-Divergence DRO}

In this section, we introduce our smoothed $f$-divergences for which the smoothed KL-divergence is a special case, and their associated DRO cost predictors. We first briefly discuss the statistical distances which appear in this paper and recap their properties.

\subsection{Common Statistical Distances}

\subsubsection{Optimal Transport Distance}
\label{ssec:wass-dist}

Denote with
\[
  W_d(P', P) = \inf \set{\int d(\xi', \xi) \, \d \gamma(\xi', \xi)}{\gamma \in \Gamma(P', P)}
\]
the optimal transport distance associated with a transport cost $d:\Sigma\times\Sigma\to\Re_+$ where $\Gamma(P', P)$ is the transport polytope of all joint distributions on $\Sigma\times\Sigma$ with given marginals $P'$ and $P$.
We assume that the transport cost $d$ is lower semicontinuous and satisfies $d(\xi, \xi)=0$ for all $\xi\in \Sigma$.
In particular, the Wasserstein distance $W(P', P)=\inf \set{\int \norm{\xi'- \xi} \, \d \gamma(\xi', \xi)}{\gamma \in \Gamma(P', P)}$, which we utilize in this paper, is associated with the transport cost $d(\xi, \xi')=\norm{\xi-\xi'}$.
This distance metrizes the weak topology on the set $\mc P(\Sigma)$ as we assume here that $\Sigma$ is compact.

\subsubsection{L\'evy-Prokhorov Metric}
\label{ssec:LP-dist}

The LP probability metric is defined as
\[
  \LP{P'}{P} \defn \inf \tset{\epsilon> 0}{P(S) \leq P'(S^\epsilon) + \epsilon, ~P'(S) \leq P(S^\epsilon) + \epsilon \quad \forall S\in \mc B(\Sigma)}
\]
where we denote with $S^\epsilon \defn \set{\xi'\in \Sigma}{\norm{\xi'-\xi}\leq \epsilon\text{ for some }\xi\in S}$ the $\epsilon$-norm inflation of a set $S$.
Like the Wasserstein distance, this metric also metrizes the weak topology when $\Sigma$ is a separable metric space \citep{prokhorov1956convergence}. Hence, the Wasserstein and LP metrics are topologically equivalent as both are compatible with the weak topology on $\mc P(\Sigma)$.
Furthermore, following \cite{strassen1965existence} the sublevel set of the LP probability metric can be written as
\begin{align*}
  & \set{P\in\mc P(\Sigma)}{\LP{P'}{P}\leq \epsilon}\\
  &\hspace{5em} =\set{P\in\mc P(\Sigma)}{\exists \gamma \in \Gamma(P, P'),~ \textstyle \int \one{\norm{\xi- \xi'}> \epsilon} \, \d \gamma(\xi, \xi') \leq \epsilon}.
\end{align*}
Hence, the LP probability metric can be associated with an optimal transport distance for a transport cost function that indicates whether some predetermined cutoff distance $\epsilon$ is exceeded or not.

\subsubsection{$f$-Divergence}\label{sec:f def}
Let $f:\Re_+\to\Re$ denote a convex and lower semicontinuous function with $f(1)=0$. First assume $\Sigma$ is a finite set. The $f$-divergence between two distributions $P$ and $P'$ in $\mathcal{P}(\Sigma)$ is defined as
\[
D_f(P', P) \defn \textstyle\sum_{\xi\in \Sigma} f\left(\tfrac{P(\xi)}{P'(\xi)}\right) P'(\xi)
\]
where we take \( 0 f(0/0) \defn 0 \) and $0 f(p/0)\defn p f^\infty$ with $f^\infty\defn \lim_{u\to\infty} f(u)/u \in \nRe$ \citep{ahmadi2012entropic}.
We remark that as $f$ is convex, $f^\infty$ bounds the subgradients of the function $f$ from above.
Next, for distributions supported on an infinite subset $\Sigma$ of $\Re^d$, the situation is slightly more complicated. By the Lebesgue decomposition theorem \citep[Section 6.9]{rudin1970real} we can decompose $P=P_c+P_s$ uniquely into a continuous part $P_c\ll P'$ and singular part $P_s \bot P'$.\footnote{Here $P_c\ll P'$ is equivalent to $P'(S)=0\implies P_c(S)=0$ for all measurable subsets $S$ in $\Sigma$ whereas $P_s \bot P'$ is equivalent to the existence of a set $S$ so that for all measurable subsets $B$ of $S$ we have $P_c(B)=0$ and for all measurable subsets $B'$ of $\Sigma\setminus S$ we have $P'(B')=0$.} In this case, the $f$-divergence is defined as
\[
  D_f(P', P) \defn \textstyle\int f\left(\tfrac{\d P_c}{\d P'}(\xi )\right) \, \d P'(\xi) + \int f^\infty \, \d P_s(\xi) 
\]
where $\tfrac{\d P_c}{\d P'}$ denotes the Radon-Nikodym derivative between $P_c$  and $P'$ \citep{nikodym_sur_1930}.
Remark that $f$-divergences are generally not symmetric in that $D_f(P', P)\neq D_f(P, P')$.
Nevertheless, it is well known that with the help of an adjoint function $f^{\circ}:t\mapsto t f(\tfrac{1}{t})$ we do have $D_f(P', P)= D_{f^{\circ}}(P, P')$.
We give several particular $f$-divergences in Table \ref{table:divergence_functions} although most attention will be reserved for the KL-divergence associated with $f(t)=-\log(t)+t-1$.

It is well known that $f$-divergences are lower semicontinuous but generally fail to be upper semicontinuous. Hence, unlike the Wasserstein distance and LP distance, they are not compatible with the weak topology on $\mc P(\Sigma)$ and hence are topologically distinct from the other two distances discussed previously.

\subsection{Smoothed $f$-Divergences and DRO Predictors}
\label{ssec:smooth-f-diverg}

The inherent topological differences between Wasserstein and LP distances and $f$-divergences result in wildly different statistical behavior.
Since Wasserstein and LP distances metrize weak convergence on $\mc P(\Sigma)$, they are ``smooth'' distances in the sense that they can capture the closeness of the empirical distribution $P_N$ to the true distribution $P$, i.e., $W(P_N,P)\rightarrow0$ and $\LP{P_N}{P}\rightarrow0$ a.s. by the Glivenko–Cantelli theorem. Moreover, as (\ref{eq: Wasserstein bound}) indicates, the Wasserstein distance enjoys a finite-sample coverage guarantee. 
Similarly, \cite[Corollary 4.2]{dudley1969speed} shows that the LP distance also admits a finite-sample guarantee
\begin{equation}
  \label{eq: LP bound}
  \Prob{\LP{P_N}{P}>\epsilon} 
  \leq
  \tfrac{\E{}{\LP{P_N}{P}}}{\epsilon}
  \leq
  C' N^{-1/(d+2)} /\epsilon
\end{equation}
where $C'$ is a known constant which depends only on the dimension $d$ and the diameter $\diam(\Sigma)=\sup_{\xi,\xi^{'}\in\Sigma} ||\xi-\xi'||$. 

On the other hand, for $f$-divergences, $D_f(P_N, Q)=f(0)+f^\infty$ whenever $Q$ is a continuous distribution as the empirical distribution $P_N$ is always finitely supported.
Hence, when the ball radius is chosen as $\eta< f(0)+f^\infty$, the associated $f$-divergence ambiguity set $\set{P'\in \mc P(\Sigma)}{D_f(P_N, P')\leq \eta}$ does not contain any continuous distribution whereas for $r\geq f(0)+f^\infty$ it contains all distributions and is hence trivial.
Only when the set $\Sigma$ is finite are finite-sample bounds for the $\kl$-divergence ambiguity sets parallel to (\ref{eq: Wasserstein bound}) available. More concretely, a finite version of Sanov's theorem guarantees that for all $N\geq 1$ we have
\begin{equation}
  \label{eq:finite:sanov}
  \Prob{D_{\kl}(P_N, P)> r}
  \leq
  (N+1)^{\abs{\Sigma}} \exp(-r N),
\end{equation}
where $|\Sigma|$ denotes the cardinality of $\Sigma$.
However, thanks to the key role that KL-divergence plays in the large deviations theory, \cite{vanparys2021data} shows that despite the lack of confidence region interpretation, the KL-DRO predictor $\hat c_{\kl, r}(x, \hat P)$, as the special case of the following $f$-divergence DRO predictor 
\begin{equation}
  \label{eq:dro:f}
   \hat c_{f, r}(x, \hat P) \defn
  \left\{
    \begin{array}{rl}
      \sup & \int \ell(x, \xi) \, \d P(\xi)\\
      \st  & D_{f}(\hat P, P) \leq r,
    \end{array}
  \right.
\end{equation}
is the least conservative predictor compared to any other pointwise feasible continuous predictors when the loss function is jointly continuous in $x$ and $\xi$. Unfortunately, although KL-DRO predictor continues to be least conservative even in the sense of uniform feasibility for finite $\Sigma$ (see Section \ref{sec:finite-alphabets}), it fails to be uniform feasible for infinite $\Sigma$ (Section \ref{sec:failure_KL_DRO}) due to the lack of guarantee in the form of \eqref{eq:finite:sanov}.

To salvage the infinite support case, we consider a  smoothing of the $f$-divergences through combining with Wasserstein or LP distances. More concretely, we define the \textit{Wasserstein-smoothed $f$-divergence} and \textit{LP-smoothed $f$-divergence} by
\begin{equation}
  \label{def:dist-type-i}
  D^\epsilon_{f, W}(P', P) = \inf\set{D_f(Q, P)}{Q\in \mc P(\Sigma), ~W(P', Q)\leq \epsilon}
\end{equation}
and
\begin{equation}
  \label{def:LPdist-type-i}
  D^\epsilon_{f, \lp}(P', P) = \inf\set{D_f(Q, P)}{Q\in \mc P(\Sigma), ~\LP{P'}{Q}\leq \epsilon}
\end{equation}
respectively, for $P\in \mc P(\Sigma)$, $P'\in \mc P(\Sigma)$ and some parameter $\epsilon\geq 0$. It is clear that both distances reduce to an $f$-divergence for the parameter $\epsilon=0$. Associated with the two smoothed $f$-divergences, we consider the two DROs:
\begin{equation}
  \label{eq:formulation-type-1}
  \hat c_{f, W, r}^{\epsilon}(x, \hat P) \defn
  \left\{
    \begin{array}{rl}
      \sup & \int \ell(x, \xi) \, \d P(\xi)\\
      \st  & P\in \mc P(\Sigma),~D^\epsilon_{f, W}(\hat P, P) \leq r
    \end{array}
  \right.
\end{equation}
and
\begin{equation}
  \label{eq:LPformulation-type-1}
  \hat c_{f, \lp, r}^{\epsilon}(x, \hat P) \defn
  \left\{
    \begin{array}{rl}
      \sup & \int \ell(x, \xi) \, \d P(\xi)\\
      \st  & P\in \mc P(\Sigma),~D^\epsilon_{f, \lp}(\hat P, P) \leq r
    \end{array}
  \right.
\end{equation}
In what follows we will call $r\geq 0$ and $\epsilon\geq 0$ the robustness and smoothness parameter, respectively.

Finally, we note that one could also consider an alternative approach to combine $f$-divergences and smooth distances that reverses the roles of the arguments, namely
\begin{equation}
  \label{def:dist-type-ii}
  D^\epsilon_{W, f}(P', P) = \inf\set{D_f(P', Q)}{Q\in \mc P(\Sigma), ~W(Q, P) \leq \epsilon}
\end{equation}
and
\begin{equation}
  \label{def:LPdist-type-ii}
  D^\epsilon_{\lp, f}(P', P) = \inf\set{D_f(P', Q)}{Q\in \mc P(\Sigma), ~\LP{Q}{P} \leq \epsilon}.
\end{equation}
Although we will mostly be concerned with the distances defined in \eqref{def:dist-type-i} and \eqref{def:LPdist-type-i}, we will discuss their alternatives \eqref{def:dist-type-ii} and \eqref{def:LPdist-type-ii} as well at the end of the paper in Section \ref{sec:computation}.

\section{Exponential-Rate Feasibility and Efficiency under Finite Support}
\label{sec:finite-alphabets}

In this section, we assume $\Sigma$ is a finite set with cardinality $2\leq|\Sigma|<\infty$ and equipped with all subsets as the $\sigma$-algebra so that any event is measurable. The probability simplex $\mathcal{P}(\Sigma)$ can be defined as $\mathcal{P}(\Sigma)=\{P\in\mathbb{R}^{|\Sigma|}:\sum_{i\in\Sigma}P(i)=1,P(i)\geq0~\forall i\in\Sigma\}$ equipped with the induced topology of $\mathbb{R}^{|\Sigma|}$ and hence the weak convergence in $\mathcal{P}(\Sigma)$ is equivalent to the usual convergence in $\mathbb{R}^{|\Sigma|}$. Under this setting, we will show that the KL-DRO enjoys exponential-rate feasibility and efficiency properties in terms of both the predictor and ambiguity set.

First, we have the uniform feasibility of the KL-DRO predictor $\hat{c}_{\kl,r}(x,\hat P)$:
\begin{theorem}[Uniform feasibility of KL-DRO under finite support]
  \label{KL_uni_feasibility}
  Suppose $|\Sigma|<\infty$. The KL-DRO predictor $\hat{c}_{\kl,r}(x,\hat P)$ satisfies the uniform feasibility \eqref{eq:feasibility-uniform} at rate $r\geq0$ for any loss function $\ell$.
\end{theorem}

Theorem \ref{KL_uni_feasibility} is derived based on a large deviations argument and is an enhancement of Theorem 3 in \citep{vanparys2021data} that states the pointwise feasibility of the same DRO. Next, we will show the least conservativeness of $\hat{c}_{\kl,r}(x,\hat P)$. We consider the class $\mathcal{C}_{u}(X\times\mathcal{P}(\Sigma))$ of upper semicontinuous predictors defined by
\begin{equation*}
\tset{\hat{c}:X\times\mathcal{P}(\Sigma)\mapsto\mathbb{R}}{\hat{c}(x,\hat P)\text{ is upper semicontinuous in }\hat P\text{ for each }x\in X}.
\end{equation*}
This upper semicontinuity restriction makes sure that the event set of interest, $\{\exists x\in X~\st~ c(x,P)>\hat c(x, P_N)\}$ in (\ref{eq:feasibility-uniform}), is measurable for any predictor. This measure-theoretic consideration will become more important when $\Sigma$ is no longer a finite set later. Here, in the finite support case, for any loss function $\ell(x,\xi)$, any $f$-divergence predictor $\hat{c}_{f,r}(x,\hat P)$ is in $\mathcal{C}_{u}(X\times\mathcal{P}(\Sigma))$ (see Lemma \ref{usc_smooth_f_DRO} with $\epsilon=0$). The following theorem shows that KL-DRO formulation is the least conservative among all upper semicontinuous predictors that are pointwise feasible or uniformly feasible.
\begin{theorem}[Least conservativeness of KL-DRO under finite support]
\label{prediction_optimality_finite}Suppose $|\Sigma|<\infty$. For any loss function $\ell(x,\xi)$, if $\hat{c}(x,\hat P)\in \mathcal{C}_{u}(X\times\mathcal{P}(\Sigma))$ satisfies \eqref{eq:feasibility} or \eqref{eq:feasibility-uniform} at rate $r>0$, then $\hat{c}_{\kl,r}(x,\hat P)\leq\hat{c}(x,\hat P)$ for any $x\in X,\hat P\in\mathcal{P}(\Sigma)$.
\end{theorem}
Theorem \ref{prediction_optimality_finite} is proven along the lines of Theorem 4 in \citep{vanparys2021data} but generalizes the latter result in two distinct directions: The considered class of predictors are extended to include upper semicontinuous predictors rather than merely continuous predictors, and the imposed requirement includes both pointwise feasibility and uniform feasibility. We remark that the out-of-sample guarantees (\ref{eq:feasibility-uniform}) and (\ref{eq:feasibility}) can always be obtained by simply associating extremely inflated cost predictions $\hat{c}(x,\hat P)$ for all $x$ and $\hat P$. In words, Theorem \ref{prediction_optimality_finite} shows that the KL predictor $\hat{c}_{\kl,r}(x,\hat P)$ is the least inflated cost prediction among those predictors in $\mathcal{C}_{u}(X\times\mathcal{P}(\Sigma))$ enjoying either (\ref{eq:feasibility}) or (\ref{eq:feasibility-uniform}).

We can also interpret the conservativeness of DRO in terms of the ambiguity set.
We define the class of DRO predictors by
\begin{align}
&\mathcal{C}_{\dro}(X\!\times\!\mathcal{P}(\Sigma)) \!=\!\{\hat{c}\in\mathcal{C}_{u}(X\!\times\!\mathcal{P}(\Sigma)):\hat{c}(x,\hat{P})\!=\!\sup\{c(x,P)\!:\!P\!\in\!\mathcal{A}(\hat{P})\},\label{DRO_predictors}\\
  & \hspace{8.4em}\,\mathcal{A}(\hat{P})\subseteq\mathcal{P}(\Sigma) \text{ is nonempty, closed and convex for any }\hat{P}\},\nonumber
\end{align}
where the ambiguity set $\mathcal{A}(\hat{P})$ only depends on $\hat P$ and is independent of the loss function. The ambiguity set is assumed to be closed and convex which is a common requirement when one is interested in DRO predictors which admit tractable algorithms. We would like to find the smallest ambiguity set in terms of the set inclusion among all pointwise feasible or uniformly feasible DRO predictors. However, a DRO predictor can happen to be feasible for some trivial loss functions. For example, if $\ell(x,\xi)$ is a constant function, even a singleton ambiguity set leads to a uniformly feasible DRO predictor. Therefore, we define the notion of \textit{reasonably feasible DRO predictor} by requiring it to be feasible \textit{at least} for all decision-free loss functions.
\begin{definition}[Decision-free loss]
A loss function $\ell(x,\xi)$ is called decision-free if $\ell(x,\xi)\equiv\ell_0(\xi)$ for some continuous function $\ell_0:\Sigma\mapsto\mathbb{R}$.
\end{definition}
\begin{definition}[Reasonable feasibility]\label{reasonable_DRO}
A DRO predictor $\hat c\in \mathcal{C}_{\dro}(X\times\mathcal{P}(\Sigma))$ is called reasonably feasible at rate $r$ if it satisfies (\ref{eq:feasibility-uniform}) at rate $r$ for any decision-free loss function $\ell(x,\xi)$.
\end{definition}
Although continuity essentially imposes no loss of generality in the finite support case, this will become important in Section \ref{sec:continuous-alphabets}. With a decision-free loss function, the pointwise feasibility (\ref{eq:feasibility}) and uniform feasibility (\ref{eq:feasibility-uniform}) are the same for DRO predictors. The next theorem establishes the least conservativeness of the KL-divergence ambiguity set, in terms of set inclusion, among the ambiguity sets of all reasonably feasible DRO predictors.

\begin{theorem}[Least conservativeness of KL ambiguity set under finite support]
\label{optimality_KL_ball}
Suppose $|\Sigma|<\infty$. If $\hat c\in \mathcal{C}_{\dro}(X\times\mathcal{P}(\Sigma))$ is a reasonably feasible DRO predictor at rate $r>0$, then we have $\{P:D_{\kl}(\hat P,P)\leq r\}\subseteq \mathcal{A}(\hat P)$ for any $\hat P\in\mathcal{P}(\Sigma)$.
\end{theorem}

Theorem \ref{optimality_KL_ball} is established via contradiction. If $\{P:D_{\kl}(\hat P,P)\leq r\}\nsubseteq \mathcal{A}(\hat P)$, then a separating hyperplane theorem asserts that there exists a decision-free loss function $\ell(x,\xi)$ such that $\hat c(x,\hat P)$ is not feasible. This theorem shows that the KL predictor $\hat{c}_{\kl,r}(x,\hat P)$ is the least conservative DRO predictor in terms of the ambiguity set inclusion. In other words, any DRO formulation which enjoys an exponential out-of-sample guarantee is necessarily at least as conservative as the KL-DRO formulation considered here.

\section{Exponential-Rate Feasibility and Efficiency under Infinite Support}
\label{sec:continuous-alphabets}
We consider the infinite support case, which is the main focus of this paper.
Before we study the properties of our predictors, we first need to address the measurability of the events considered in \eqref{eq:feasibility-uniform} and \eqref{eq:feasibility}, namely the measurability of $\{\E{P}{\ell(x, \xi)}>\hat c(x, P_N)\}$ for any $x\in X$ and $\{\exists x\in X~\st~ \E{P}{\ell(x, \xi)}>\hat c(x, P_N)\}$. We cannot guarantee that these events are measurable for all predictors $\hat c$ but it can be shown that they are indeed measurable for any $\hat c\in \mathcal{C}_{u}(X\times\mathcal{P}(\Sigma))$ (see Lemma \ref{event_measurability}, where the topology on $\mathcal{P}(\Sigma)$ is the weak topology discussed in Section \ref{sec:infinite_LDP}). Moreover, if the loss function $\ell$ is upper semicontinuous in $\xi$ for any $x\in X$, then any (smoothed) $f$-divergence DRO predictor is in $\mathcal{C}_{u}(X\times\mathcal{P}(\Sigma))$ (see Lemma \ref{usc_smooth_f_DRO}). Therefore, throughout Section \ref{sec:continuous-alphabets}, the loss function $\ell(x,\xi)$ is assumed to be upper semicontinuous in $\xi$ for any $x\in X$ without further clarification.

\subsection{Lack of Feasibility of the KL-DRO Predictor}\label{sec:failure_KL_DRO}

First, we demonstrate that even if the loss function is jointly continuous, the KL-DRO predictor can still fail to be uniformly feasible:
\begin{theorem}[Lack of feasibility of KL-DRO under infinite support]
\label{counterexample_any_r}Let $X=[0,1]$, $\Sigma=[-1,1]$ and the true probability $P$ be the uniform distribution on $\Sigma$. For any $r>0$, there exists a loss function $\ell_r(x,\xi)$ jointly continuous on $X\times \Sigma$ such that
\[
\Prob{\exists x\in X~\text{s.t.}~ c(x, P)>\hat{c}_{\kl,r}(x,P_{N})}=1.
\]
Moreover, there exists a loss function $\ell(x,\xi)$ jointly continuous on $X\times \Sigma$ such that
\begin{equation}
  \Prob{c(\hat{x}_{r}(P_{N}),P)>\hat{c}_{\kl,r}(\hat{x}_{r}(P_{N}),P_{N})}  =1\label{open_question}  
\end{equation}
and
\begin{equation}
\Prob{\textstyle\min_{x\in X}c(x,P)>\min_{x\in X}\hat{c}_{\kl,r}(x,P_{N})}  =1\label{failed_bound}
\end{equation}
for any $0<r<\log(2)$, where $\hat{x}_{r}(P_{N})\in\arg\min_{x\in X}\hat{c}_{\kl,r}(x,P_{N})$.
\end{theorem}

The previous result is established by a counterexample. Intriguingly, our constructed counterexample relies on the concept of ``uniformly distributed modulo 1'' \citep[Chapter 1 Section 6]{kuipers2012uniform} in analytic number theory and uses a generalization of the interesting fact that $\{\sin{n},n\in\mathbb{N}\}$ is dense on $[-1,1]$.

We make two key observations regarding Theorem \ref{counterexample_any_r}: 1) The KL-DRO predictor $\hat{c}_{\kl,r}$ can fail, rather miserably, to be uniformly feasible at any rate $r>0$ in \eqref{eq:feasibility-uniform} even if the loss function is jointly continuous. This should be contrasted with \cite[Theorem 10]{vanparys2021data} which establishes the pointwise feasibility of the KL-DRO predictor for a general compact set $\Sigma$ when the loss function is continuous. This also reveals that uniform feasibility is a much stronger requirement than pointwise feasibility. 2) The DRO optimal value can fail to bound the true cost evaluated at the DRO optimal decision, as shown in \eqref{open_question}, and the true optimal value, as shown in \eqref{failed_bound}. In fact, the failure in point 2) is completely due to the failure of uniform feasibility discussed in point 1), which also explains why uniform feasibility is a suitable requirement for decision making. Additionally, (\ref{open_question}) also gives a negative answer to an open question raised in the penultimate paragraph of \cite{vanparys2021data} about whether $\Prob{c(\hat{x}_{r}(P_{N}),P)>\hat{c}_{\kl,r}(\hat{x}_{r}(P_{N}),P_{N})}$ can decay to zero exponentially fast at rate $r$.

\subsection{Recovering Feasibility and Efficiency through Smoothing}
As shown in the previous subsection, the KL-DRO predictor is no longer a good choice from the perspective of uniform feasibility when $P$ is not finitely supported. In this subsection, we show that by using the smoothed $\kl$-divergence ambiguity set, we can salvage the uniform feasibility, and up to a certain extent efficiency, with mild conditions on the loss function $\ell(x,\xi)$.

\begin{theorem}[Uniform feasibility of smoothed KL-DRO under infinite support]
\label{smooth_KL_uni_feasibility}Suppose $\Sigma\subset\Re^d$ is a compact set and the loss function $\ell(x,\xi)$ is upper semicontinuous in $\xi$ for any $x\in X$. For any $\epsilon>0$ and $r\geq0$, $\hat{c}_{\kl, W, r}^{\epsilon}(x, \hat P)$ and $\hat{c}_{\kl, \lp, r}^{\epsilon}(x, \hat P)$ satisfy the uniform feasibility (\ref{eq:feasibility-uniform}) at rate $r$.
\end{theorem}

Theorem \ref{smooth_KL_uni_feasibility} is proven along the same lines as Theorem \ref{KL_uni_feasibility} but uses a modified version of Sanov's theorem stated in Theorem \ref{thm:sanov-smooth}.
Furthermore, we can strengthen the asymptotic result in Theorem \ref{smooth_KL_uni_feasibility} to finite-sample guarantees, by using the finite-sample upper bound in Theorem \ref{thm:sanov-smooth} instead of its asymptotic version. Let $C(\epsilon,\Sigma,||\cdot||)$ be the minimal covering number of $\Sigma$ under the Euclidean norm $||\cdot||$ using balls of radius $\epsilon$.
In case of LP smoothing we can show that
\[
  \Prob{\exists x\in X~\st~c(x,P)>\hat c_{\kl, \lp, r}(x, P_N)}\leq \left(\tfrac{8}{\epsilon}\right)^{C(\epsilon/2, \Sigma, \norm{\cdot})}\exp(-rN)\quad \forall N\geq 1.
\]
A structurally identical guarantee can be established for the Wasserstein smoothed formulation $\hat{c}_{\kl, W, r}^{\epsilon}(x, \hat P)$. This implies that if one is interested in merely bounding \eqref{beginning prob} below a desired probability level $\alpha\in (0,1)$, as done for instance in \citep{lam2019recovering, duchi2021statistics}, a radius $r_N = \mc O(1/N)$ suffices. However, whereas these latter results require a Donsker class restriction on function class $\set{\ell(x, \cdot)}{x\in X}$, no such restriction is imposed here. Nonetheless, there is no free lunch and instead we require here a bound on the complexity of the event set $\Sigma$ (but not the loss function $\ell$) through its covering number. We do remark that bounding covering numbers is relatively straightforward and we can for instance use the generic bound $C(\epsilon, \Sigma, \norm{\cdot})\leq (2 \diag(\Sigma) \sqrt{d}/\epsilon)^d$ \cite{shalev2014understanding}.

Next we focus on efficiency. Like in the finite support case, we consider least conservativeness both for the prediction value compared to predictors in $\mathcal{C}_{u}(X\times\mathcal{P}(\Sigma))$, and the ambiguity set compared to those in reasonably feasible DRO predictors defined in Definition \ref{reasonable_DRO}. However, unlike the finite support case, an \textit{exact} least conservative predictor does not always exist for a general compact set $\Sigma$. In fact, notice that both the ambiguity sets and prediction values of $\hat{c}_{\kl, W, r}^{\epsilon}(x, \hat P)$ and $\hat{c}_{\kl, \lp, r}^{\epsilon}(x, \hat P)$ are non-increasing when $\epsilon$ decreases, which means no predictor among $\hat{c}_{\kl, W, r}^{\epsilon}(x, \hat P)$ and $\hat{c}_{\kl, \lp, r}^{\epsilon}(x, \hat P)$ can be least conservative. This issue can be addressed by focusing on the efficiency as $\epsilon\rightarrow0$ asymptotically.
The following relation between the smoothed $f$-divergence DRO and the usual $f$-divergence DRO is essential to establish our desired asymptotic efficiency.
\begin{theorem}[Properties of smoothed ambiguity sets and DRO predictors]\label{monotonicity_DRO}
Suppose $\Sigma\subset\Re^d$ is a compact set and $r\geq0$. The ambiguity sets $\{P\in\mathcal{P}(\Sigma):D^\epsilon_{f, W}(\hat P, P)\leq r\}$ and $\{P\in\mathcal{P}(\Sigma):D^\epsilon_{f, \lp}(\hat P, P)\leq r\}$ are non-increasing when $\epsilon$ decreases and have the limit
\begin{align*}
  & \textstyle\cap_{\epsilon>0}\{P\in\mathcal{P}(\Sigma):D^\epsilon_{f, W}(\hat P, P)\leq r\}=\cap_{\epsilon>0}\{P\in\mathcal{P}(\Sigma):D^\epsilon_{f, \lp}(\hat P, P)\leq r\}\\
  & \hspace{22em}\textstyle=\{P\in\mathcal{P}(\Sigma):D_{f}(\hat P,P)\leq r\}.
\end{align*}
Similarly, if the loss function $\ell(x,\xi)$ is upper semicontinuous in $\xi$ for any $x\in X$, $\hat{c}_{f, W, r}^{\epsilon}(x, \hat P)$ and $\hat{c}_{f, \lp, r}^{\epsilon}(x, \hat P)$ are non-increasing when $\epsilon$ decreases and have the limit
\[
\lim_{\epsilon\downarrow0}\hat{c}_{f, W, r}^{\epsilon}(x, \hat P)=\lim_{\epsilon\downarrow0}\hat{c}_{f, \lp, r}^{\epsilon}(x, \hat P)=\hat{c}_{f,r}(x, \hat P)
\]
for each $x\in X, \hat P\in\mathcal{P}(\Sigma)$.
\end{theorem}

Now we consider asymptotic efficiency in terms of the prediction value:

\begin{theorem}[Asymptotic least conservativeness of smoothed KL-DRO]\label{prediction_optimality_infinite}
Suppose $\Sigma\subset\Re^d$ is a compact set and the loss function $\ell(x,\xi)$ is upper semicontinuous in $\xi$ for any $x\in X$. If $\hat{c}(x,\hat{P})\in\mathcal{C}_{u}(X\times\mathcal{P}(\Sigma))$ satisfies (\ref{eq:feasibility}) or (\ref{eq:feasibility-uniform}) at rate $r>0$, then we have%
\[
\lim_{\epsilon\downarrow0}\hat{c}_{\kl,W,r}^{\epsilon}(x,\hat{P})=\lim_{\epsilon\downarrow0}\hat{c}_{\kl,\lp,r}^{\epsilon}(x,\hat{P})=\hat{c}_{\kl,r}(x,\hat{P})\leq\hat{c}(x,\hat{P})
\]
for any $x\in X,\hat{P}\in\mathcal{P}(\Sigma)$. Consequently, the following are equivalent:

\begin{itemize}[leftmargin=1em]
\item There exists a least conservative predictor $\hat{c}(x,\hat{P})\in\mathcal{C}_{u}(X\times\mathcal{P}(\Sigma))$ satisfying (\ref{eq:feasibility}) (or (\ref{eq:feasibility-uniform})), i.e., $\hat{c}(x,\hat{P})$ satisfies (\ref{eq:feasibility}) (or (\ref{eq:feasibility-uniform})) and $\hat{c}(x,\hat{P})\leq\hat{c}^{\prime}(x,\hat{P}),\forall(x,\hat{P})\in X\times\mathcal{P}(\Sigma)$ for any $\hat{c}^{\prime}\in\mathcal{C}_{u}(X\times\mathcal{P}(\Sigma))$ that also satisfies (\ref{eq:feasibility}) (or (\ref{eq:feasibility-uniform})).

\item $\hat{c}_{\kl,r}(x,\hat{P})$ is feasible in (\ref{eq:feasibility}) (or (\ref{eq:feasibility-uniform})) and $\hat{c}_{\kl,r}(x,\hat{P})$ is the least conservative predictor.
\end{itemize}
\end{theorem}

Theorem \ref{prediction_optimality_infinite} shows that the smoothed KL-DRO predictors are asymptotically least conservative predictors among $\mathcal{C}_u(X\times \mathcal{P}(\Sigma))$ in the following sense. If there exists a least conservative predictor, then the usual KL-DRO predictor is the least conservative. On the other hand, if the usual KL-DRO predictor $\hat{c}_{\kl,r}(x,\hat{P})$ is not feasible (e.g., under the counterexample in Section \ref{sec:failure_KL_DRO}), then there will be no least conservative predictor. By letting $\epsilon\downarrow0$ in $\hat{c}_{\kl,W,r}^{\epsilon}(x,\hat{P})$ and $\hat{c}_{\kl,\lp,r}^{\epsilon}(x,\hat{P})$ we can find a sequence of predictors that are asymptotically least conservative. Putting these together, if we use the cost predictor $\hat{c}_{\kl,W,r}^{\epsilon}(x,\hat{P})$ or $\hat{c}_{\kl,\lp,r}^{\epsilon}(x,\hat{P})$ with a small $\epsilon>0$, then we simultaneously enjoy uniform feasibility and achieve nearly least conservativeness. In comparison, in all the previous approaches described in Sections \ref{sec:complexity} and \ref{sec:DRO literature}, namely adding margin $\delta$ to SAA, DRO that selects the ambiguity set as a confidence region of $P$, and DRO that utilizes variability regularization, attaining uniform feasibility would require calibrations that depend on at least some model parameters such as variability. Moreover, for these approaches, there are no corresponding conservativeness results in the sense we described. On the other hand, we caution that when the loss function class is complex, we anticipate that our predictor could blow up as $\epsilon\to0$, i.e., our cost predictor becomes highly sensitive to $\epsilon$. That is, while our DRO is ``automatic" to attain uniform feasible and nearly least conservative predictors, it is anticipated to suffer the same pathological issues when the function class is too complex to have reasonable predictors in SAA.

Finally, we investigate asymptotic efficiency in terms of the ambiguity set. As in the finite support case, we require the DRO predictor to be reasonably feasible defined in Definition \ref{reasonable_DRO}.
\begin{theorem}[Asymptotic least conservativeness of smoothed KL ambiguity set]\label{optimality_KL_ball_continuous}
Suppose $\Sigma\subset\Re^d$ is a compact set. If $\hat c(x,\hat P)\in \mathcal{C}_{\dro}(X\times\mathcal{P}(\Sigma))$ is a reasonably feasible DRO predictor at rate $r>0$ associated with the ambiguity sets $\mathcal{A}(\hat P)$, then we have 
\begin{align*}
  &\cap_{\epsilon>0}\{P\in\mathcal{P}(\Sigma):D^\epsilon_{\kl, W}(\hat P, P)\leq r\}=\cap_{\epsilon>0}\{P\in\mathcal{P}(\Sigma):D^\epsilon_{\kl, \lp}(\hat P, P)\leq r\}\\
  &\hspace{18.5em}\textstyle=\{P\in\mathcal{P}(\Sigma):D_{\kl}(\hat P,P)\leq r\}\subseteq \mathcal{A}(\hat P).
\end{align*}
for any $\hat P\in\mathcal{P}(\Sigma)$.
\end{theorem}

Similar to the finite support case, Theorem \ref{optimality_KL_ball_continuous} is proved by contradiction. In particular, when $\{P\in\mathcal{P}(\Sigma):D_{\kl}(\hat P,P)\leq r\}\nsubseteq \mathcal{A}(\hat P)$, the Hahn-Banach theorem (which can be seen as the functional version of the separating hyperplane theorem) asserts that there exists a decison-free loss function such that $\hat c(x,\hat P)$ is not feasible.

\section{Connection with Other DRO Formulations}
\label{sec:empirical-predictors}
We connect and compare our smoothed DRO with other DRO formulations, including the empirical DRO that we have discussed in the introduction, as well as other smoothing approaches. 

\subsection{Empirical DRO Predictors}
Empirical $f$-divergence DRO takes the empirical distribution $P_N$ as its ball center, which is similar to our $\hat c_f(x,P_N)$. However, it restricts the support of the distribution in the ambiguity set to be that of $P_N$. That is, the empirical $f$-divergence DRO  predictor is given by
\begin{equation}
  \label{eq:dro:ef}
  \tilde c_{f, r}(x, \hat P) \defn \left\{
    \begin{array}{rl}
      \sup & \int \ell(x, \xi) \, \d P(\xi)\\
      \st  & P\ll \hat P,~ D_{f}(\hat P, P) \leq r
    \end{array}
  \right.
\end{equation}
where here the ambiguity set contains only distributions supported on the historically observed data points and hence no cost predictions will exceed the essential supremum of the loss with respect to the empirical distribution.
Such empirical $f$-divergence DRO formulations have been studied by \cite{ahmadi2012entropic, shapiro2017distributionally} who prove a strong convex dual representation, and also by \citep{lam2017empirical,duchi2018variance,lam2019recovering,gotoh2018robust,duchi2021statistics} who derive its asymptotic properties and statistical equivalence to the empirical likelihood. While $\tilde c_{f, r}(x, \hat P)$ appears very similar to our $\hat c_{f, r}(x, \hat P)$, with the only difference being the support of $P$, this difference and the subsequent inability to exceed the corresponding essential supremum will in fact have substantial impacts when considering the exponential decay of our statistical bounds.

To understand the above more concretely, we consider a more general notion of empirical DRO predictor. We denote any predictor $\tilde c:X\times \mc P(\Sigma)\to\Re$ as \textit{empirical} when it is the case that $\tilde c(x, \hat P)\leq \esssup_{\xi\sim \hat P} \ell(x,\xi)=\inf \tset{a\in\Re}{\hat P[\ell(x, \xi)\leq a]=1}$.
In particular, an empirical cost prediction never exceeds the worst-case observed loss. The following result indicates that all empirical DRO predictors generally fail to be pointwise feasible in (\ref{eq:feasibility}), let alone uniformly feasible in \eqref{eq:feasibility-uniform}.

\begin{theorem}[Infeasibility of empirical DRO predictors]
  \label{proposition:el-infeasibility}
  Let $\Sigma\subset\Re^d$ be a set with cardinality $|\Sigma|\geq2$. Consider any measurable predictor $\tilde c:X\times \mc P(\Sigma)\to \Re$ with the property that $\tilde c(x, P_N)\leq \esssup_{\xi\sim P_N}\ell(x, \xi)$.
  If $\inf_{\xi\in\Sigma} \ell(x, \xi) < \sup_{\xi\in \Sigma}\ell(x, \xi)$, we have
  \[
    \textstyle\sup_{P\in \mc P(\Sigma)} \limsup_{N\to\infty}\frac{1}{N}\log \Prob{c(z, P)>\tilde c(x, P_N)}=0.
  \]
\end{theorem}

In contrast, we show that the smoothed counterparts of empirical DRO, given by
\[
  \tilde{c}^\epsilon_{f, W, r}(x,\hat P)
  = 
  \left\{
    \begin{array}{rl}
      \sup & \tilde c_{f, r}(x,  Q)\\
      \st  & W(\hat P, Q) \leq \epsilon
    \end{array}
  \right.
\]
and
\[
  \tilde{c}^\epsilon_{f, \lp, r}(x,\hat P)
  = 
  \left\{
    \begin{array}{rl}
      \sup & \tilde c_{f, r}(x,  Q)\\
      \st  & \lp(\hat P, Q) \leq \epsilon
    \end{array}
  \right.
\]
are always feasible for any $\epsilon>0$. In fact, we have the following equivalence.

\begin{theorem}[Equivalence between DRO and empirical DRO]
  \label{lemma:hr:primal:reduction}
  For all $x\in X$ and $\epsilon>0$ we have
  $\hat c^{\epsilon}_{f, \lp, r}(x,P_N) = \tilde{c}^\epsilon_{f, \lp, r}(x,P_N)$ for any empirical distribution $P_N$.
\end{theorem}

Theorem \ref{lemma:hr:primal:reduction} directly implies that the LP smoothed empirical DRO predictor $\tilde{c}^\epsilon_{\kl, \lp, r}(x,\hat P)$ will inherit the uniform feasibility criteria (\ref{eq:feasibility-uniform}) as shown in Theorem \ref{smooth_KL_uni_feasibility}. 
Because of a standard inequality between LP and Wasserstein metrics (c.f., inequality \eqref{ineq_LP_W}), 
the Wasserstein smoothed predictor $\tilde{c}^\epsilon_{\kl, \lp, r}(x,\hat P)$ must satisfy the uniform feasibility criteria (\ref{eq:feasibility-uniform}) as well.

\subsection{Other DRO Smoothing Approaches}

We compare our smoothed DRO formulations to other smoothing approaches in the literature.
First, in the context of hypothesis testing, the smoothed distance $D^\epsilon_{\kl, \lp}$ is implicitly proposed by \cite{zeitouni1991universal} to derive universally optimal hypothesis tests (see also \cite{huber1965robust}) and, like our setting, infinite support in that problem also requires a much more careful treatment then its finite support counterpart. This line of work was picked up by \cite{yang2018robust}, who remark that at a fundamental level most of the difficulty stems from a lack of upper semicontinuity of the KL-divergence in the infinite support case. They propose in the context of univariate distributions the smoothed distance
  \(
  D^\epsilon_{{\rm{L}}, \kl}(P', P) \defn \inf\set{D_f(P', Q)}{Q\in \mc P(\Sigma), ~{\rm{L}}(Q, P) \leq \epsilon}
  \)
  which is  similar to our smoothed distance defined in \eqref{def:LPdist-type-ii} in which however the topologically equivalent L\'evy metric
  \[
    {\rm{L}}(Q, P) \defn \inf\set{\epsilon>0}{P([-\infty, x-\epsilon]) - \epsilon \leq Q([-\infty, x])\leq P([-\infty, x+\epsilon]) + \epsilon},
  \]
  is used instead. They show that this smoothed distance $P'\mapsto D^\epsilon_{{\rm{L}}, \kl}(P', P)$ when $P$ is a continuous distribution is upper semicontinuous from which a universally optimal hypothesis tests can be derived. Nonetheless, with a different motivation, all these works do not touch upon our core study of uniform feasibility and least conservativeness and their connections with DRO.

A seemingly orthogonal way to combine an $f$-divergence and Wasserstein or LP distances would be to consider their epi-addition which, taking the Wasserstein distance as an example, is defined as
\[
  D'_{f, W}(P', P) = \inf\set{W(P', Q) + D_f(Q, P) }{Q\in \mc P(\Sigma)}.
\]
Such distances are proposed by \cite{reid2011information, dupuis2022formulation}. We can associate with this distance the DRO formulation
\begin{equation}
  \label{eq:formulation-type-prime}
  \hat c'_{f, W, r}(x, \hat P) \defn
  \left\{
    \begin{array}{rl}
      \sup & \int \ell(x, \xi) \, \d P(\xi)\\
      \st  & P\in \mc P(\Sigma),~D'_{f, W}(\hat P, P) \leq r
    \end{array}
  \right.
\end{equation}
As the next result reveals, such a DRO formulation can be reduced to our smoothed DRO in (\ref{def:dist-type-i}).

\begin{theorem}[Reduction of epi-addition-DRO to smoothed DRO]
  \label{lemma:smoothing-dro-equivalence}
We have
\[
\hat c'_{f, W, r}(x, \hat P)=\sup \tset{\hat c^\epsilon_{f,  W,  r-\epsilon}(x, \hat P)}{0\leq \epsilon \leq r}
\]
for all $x\in X$ and $\hat P\in \mc P(\Sigma)$.
\end{theorem}

Hence, the alternative family of DRO formulations in (\ref{eq:formulation-type-prime}) is in some sense subsumed in our more general proposed family.  Despite that both \citep{reid2011information} and \citep{dupuis2022formulation} discuss several statistical merits of this family of distances, as remarked in \cite{dupuis2022formulation} its computational tractability has not been investigated. In Section \ref{sec:computation} we will establish the tractability of our smoothed DRO $\hat c^\epsilon_{f,  W,  r}(x, \hat P)$ and hence also $\hat c'_{f, W, r}(x, \hat P)$.

Lastly, recall that the KL-divergence ball does not admit a confidence set interpretation as $\set{Q}{D_{\kl}(P_N, Q)\leq r}$ cannot contain the data generating distribution $P$ if it is continuous.  When using the smoothed distance $D^\epsilon_{\kl,\lp}$ or $D^\epsilon_{\kl,W}$ the sets $\tset{Q}{D^\epsilon_{\kl,\lp}(P_N, Q)\leq r}$ and $\tset{Q}{D^\epsilon_{\kl,W}(P_N, Q)\leq r}$ contain $P$ with probability tending to one exponentially fast with increasing sample size for any $\epsilon>0$ following Lemma \ref{thm:sanov-smooth}.
Intuitively, as we have that $D^\epsilon_{\kl,\lp}$ and $D^\epsilon_{\kl,W}$ are both smaller than $D_{\kl}$, this stipulates that although the KL divergence ball enjoys no confidence region interpretation itself, such interpretation can be salvaged by a judiciously chosen inflation. A similar observation was made recently in \citep{goldfeld2020convergence} by studying the distance $D_{KL}(P\otimes \mc N(0, \epsilon^2), P'\otimes \mc N(0, \epsilon^2))$ where $\otimes$ stands for convolution and $\mc N(0, \epsilon^2)$ denotes a zero mean Gaussian distribution with variance $\epsilon^2$. The data processing inequality guarantees that $D_{KL}(P\otimes \mc N(0, \epsilon^2), P'\otimes \mc N(0, \epsilon^2))\leq D_{KL}(P, P')$ and hence also the set $\set{Q}{D_{KL}(P\otimes \mc N(0, \epsilon^2), Q\otimes \mc N(0, \epsilon^2))\leq r}$ can be interpreted as an inflation of the KL divergence ball.
We have from \cite[Equation 19b]{goldfeld2020convergence} that $\E{}{D_{KL}(P_N\otimes \mc N(0, \epsilon^2), P\otimes \mc N(0, \epsilon^2))}\leq \exp(\tfrac{\diam(\Sigma)^2}{\epsilon^2})/N$ which via Markov's inequality implies that
\[
\Prob{D_{KL}(P_N\otimes \mc N(0, \epsilon^2), P\otimes \mc N(0, \epsilon^2))\leq r}\geq ~1-\exp(\tfrac{\diam(\Sigma)^2}{\epsilon^2})/(r N).
\]

\section{Relative Efficiency}

Previous sections have established that our smoothed KL divergences $D^{\epsilon}_{\kl,W}$ and $D^{\epsilon}_{\kl,\lp}$ appear efficient in the sense of being the least conservative. Here, we develop a scheme to measure the efficiency loss compared to $D^{\epsilon}_{\kl,W}$ and $D^{\epsilon}_{\kl,\lp}$ if we use other DRO distances, including other $f$-divergences and the Wasserstein distance.
First of all, we describe a general procedure to scale the ambiguity set of any alternative DRO formulation
\begin{equation}
  \hat c_{D, r}(x, \hat P) \defn
  \left\{
    \begin{array}{rl}
      \sup & \int \ell(x, \xi) \, \d P(\xi)\\
      \st  & P\in \mc P(\Sigma),~D(\hat P, P) \leq r
    \end{array}
  \right.
\end{equation}
associated with a convex lower semicontinous statistical distance $D$ to be uniformly feasible in \eqref{eq:feasibility-uniform}.
For a given statistical distance function $D$ define
\[
  R^\epsilon_{D,W{\rm{\,or\,}}\lp}(r)\defn\left\{
    \begin{array}{rl}
      \sup_{P,\hat P\in\mathcal{P}(\Sigma)}  &  D(\hat P,P)\\
      \st  &  D^\epsilon_{\kl, {W{\rm{\,or\,}}\lp}}(\hat P,P)\leq r
    \end{array}\right.
\]
as well as the limiting case
\[
  R_{D}(r)\defn\left\{
    \begin{array}{rl}
      \sup_{P,\hat P\in\mathcal{P}(\Sigma)}  &  D(\hat P,P)\\
      \st  &  D_{\kl}(\hat P,P)\leq r.
    \end{array}\right.  
\]
Clearly, we have both $R_{D}(r)\leq R^\epsilon_{D,W}(r)$ and $R_{D}(r)\leq R^\epsilon_{D,\lp}(r)$ for all $r\geq 0$ and $\epsilon\geq 0$.
The next result implies the infimal radius such that a DRO predictor based on certain distance function is uniformly feasible in \eqref{eq:feasibility-uniform} at rate $r$ is precisely $R_D(r)$.
From Theorems \ref{smooth_KL_uni_feasibility} and \ref{optimality_KL_ball_continuous}, the following corollary follows almost immediately.

\begin{corollary}[Uniform feasibility of general DRO]\label{feasibility_smooth_f_DRO}
  Suppose $\Sigma\subset\Re^d$ is a compact set and the loss function $\ell(x,\xi)$ is upper semicontinuous in $\xi$ for any $x\in X$.
  The smoothed DRO predictors $\hat{c}_{D, R_{D,W}^{\epsilon}(r)}(x,\hat{P})$ and $\hat{c}_{D, R_{D,\lp}^{\epsilon}(r)}(x,\hat{P})$ for $\epsilon>0$ are uniformly feasible in (\ref{eq:feasibility-uniform}) at rate $r\geq 0$.
  On the other hand, if $0\leq\eta<R_{D}(r)$ and $r>0$, then there exists a decision-free loss function $\ell(x,\xi)$ such that $\hat{c}_{D,\eta}(x,\hat P)$ is not pointwise feasible at rate $r$. 
\end{corollary}

Let $\hat{c}_1(x,\hat P)=\sup_{P\in \mc A_1(\hat P)}\int \ell(x, \xi)\d P(\xi)$, $\hat{c}_2(x,\hat P)=\sup_{P\in \mc A_2(\hat P)}\int \ell(x, \xi)\d P(\xi)$. We now introduce a distance between both DRO formulations that defines a loss of efficiency when comparing a DRO formulation $\hat{c}_1$ with our efficient KL-DRO formulation $\hat{c}_2=\hat{c}_{\kl,r}$. In order to do so we consider the worst-case pointwise excess $|\hat{c}_1(x,\hat P)-\hat{c}_{2}(x,\hat P)|$ for any given $x\in X,\hat P\in \mathcal{P}(\Sigma)$ when the loss function varies. Since this excess is not scale-invariant, we need to restrict the scale of the loss function to obtain a meaningful worst-case quantification. 
In particular, we consider the collection of all $1$-Lipschitz loss functions, i.e., $\ell\in \mc L_1$ where $\|\ell(x,\xi')-\ell(x, \xi)\| \leq \norm{\xi'-\xi}$ for all $\xi,\xi'\in \Sigma$ and $x\in X$. Among these functions, the worst-case pointwise excess can be characterized by the following theorem.
\begin{theorem}[Relative efficiency loss]\label{quantification2}
Let $\hat{c}_1(x,\hat P),\hat{c}_2(x,\hat P)\in\mathcal{C}_{\dro}(X\times\mathcal{P}(\Sigma))$ be two DRO predictors. Let $\mathcal{A}_1(\hat P)$ and $\mathcal{A}_2(\hat P)$ denote the associated ambiguity set of $\hat{c}_1$ and $\hat{c}_2$ respectively. Then, for any $x\in X$ and $\hat P\in\mathcal{P}(\Sigma)$
\begin{align}
  & \textstyle\sup_{\ell\in \mc L_1}|\hat{c}_{1}(x,\hat{P})-\hat{c}_{2}(x,\hat{P})|\nonumber\\
  &\hspace{0.7em}=\textstyle\max\left\{  \sup_{P\in\mathcal{A}_{1}(\hat{P})}\inf_{Q\in\mathcal{A}_{2}(\hat{P})}W(P,Q),\sup_{Q\in\mathcal{A}_{2}(\hat{P})}\inf_{P\in\mathcal{A}_{1}(\hat{P})}W(P,Q)\right\}.
\label{PH_distance}
\end{align}
In particular, if $\hat{c}_1(x,\hat P)$ is a reasonably DRO predictor at rate $r$ and $\hat{c}_2(x,\hat P)=\hat{c}_{\kl,r}(x,\hat P)$, then we have
\[
0\leq \hat{c}_{1}(x,\hat{P})-\hat{c}_{\kl,r}(x,\hat{P})\leq  \sup_{\xi'\neq\xi\in \Sigma}{ \textstyle\frac{\|\ell(x,\xi')-\ell(x, \xi)\|}{\norm{\xi'-\xi}}}\cdot \sup_{P\in\mathcal{A}_{1}(\hat{P})}\inf_{Q:D_{\kl}(\hat{P},Q)\leq r}W(P,Q).
\]
\end{theorem}
Theorem \ref{quantification2} is proved via the Sion's minimax theorem \cite[Theorem 4.2']{sion1958general} and the Kantorovich-Rubinstein dual representation \begin{equation*}
   W(P',P) \!\defn\! \sup \set{\textstyle\int \phi(\xi') \, \d P'(\xi')\!-\!\int \phi(\xi) \, \d P(\xi)\!\!}{\!\!|\phi(\xi')\!-\!\phi(\xi)|\!\leq\! \norm{\xi'\!-\!\xi}~ \forall \xi', \xi\in \Sigma}.
 \end{equation*}
In particular, Kantorovich-Rubinstein dual representation connects the $1$-Lipschitz function class $\mc L_1$ and the resulting Wasserstein distance $W(P,Q)$. The excess in (\ref{PH_distance}) is better known as the Pompeiu-Hausdorff distance between $\mathcal{A}_{1}(\hat{P})$ and $\mathcal{A}_{2}(\hat{P})$ (see e.g. \cite[Example 4.13]{rockafellar2009variational}). So Theorem \ref{quantification2} tells us that the relative efficiency between any two DRO predictors can be quantified as the Pompeiu-Hausdorff distance between their ambiguity sets. We present two important examples in using Theorem \ref{quantification2} below.

\begin{example}[$f$-divergence]
The theory of joint range of $f$-divergences asserts that the radius $R_{D_f}(r)$ does not depend on the cardinality of $\Sigma$ as long as $|\Sigma|\geq4$ (see, e.g., \citep{harremoes2011on} and \citep[Section 7.4]{polyanskiy2023information}). More concretely, given any two $f$-divergences $D_{f_1}$ and $D_{f_2}$, it is proved that the following joint ranges are identical
\begin{align*}
  & \{(D_{f_{1}}(\hat P,P),D_{f_{2}}(\hat P,P)):\\
  & \hspace{5.5em}\hat P\text{ and }P\text{ are probability distributions on some measure space }(\Omega,\mathcal{F})\}\\
=  & \{(D_{f_{1}}(\hat P,P),D_{f_{2}}(\hat P,P)):\hat P\text{ and }P\text{ are probability distributions on }\{1,\ldots,4\}\}.\\
=  & \{(D_{f_{1}}(\hat P,P),D_{f_{2}}(\hat P,P)):\hat P,P\in\mathcal{P}(\Sigma)\},|\Sigma|\geq4.
\end{align*}%
Therefore, it suffices to consider $\hat P$ and $P$ as 4-dimensional probability vectors when computing $R_{D_f}(r)$ for $|\Sigma|\geq4$ (which we will assume in this example). Here we compute $R_{D_f}(r)$ for $f$-divergences listed in Table \ref{table:divergence_functions} either by numerical optimization or the known theoretical characterization of the joint range. In particular, we have the trivial cases $R_{D_{\lk}}(r)=\infty,\forall r>0$ and $R_{D_{\chi^2}}(r)=\infty,\forall r>0$, which is obtained by choosing $\hat{P}=(\exp(-1/t^{2}) ,1-\exp(-1/t^{2}) ,0,0)$ and $P=(t,1-t,0,0)$ for $D_{\lk}$ (giving $D_{\lk}(\hat{P},P)\rightarrow\infty$ and $D_{\kl}(\hat{P},P)\rightarrow0$ as $t\rightarrow0^{+}$) and choosing $\hat{P}=(t,1-t,0,0)$ and $P=(t^{3},1-t^{3},0,0)$ for $D_{\chi^2}$ (giving $D_{\chi^2}(\hat{P},P)\rightarrow\infty$ and $D_{\kl}(\hat{P},P)\rightarrow0$ as $t\rightarrow0^{+}$). Figure \ref{figure:f_radius} displays the radius $R_{D_f}(r)$ for other $f$-divergences.

\begin{figure}[htp]
\centering
\includegraphics[width=0.5\textwidth]{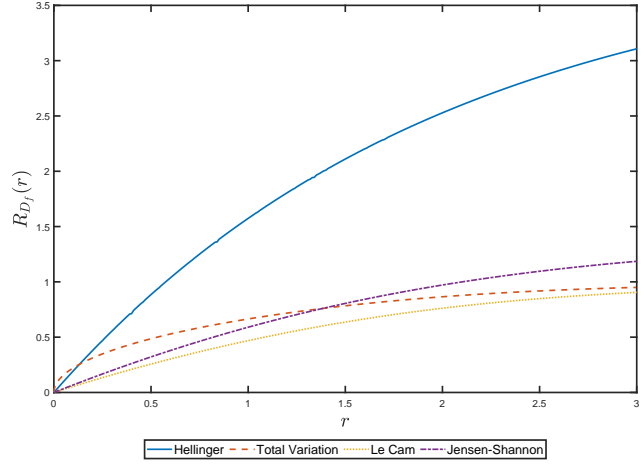}
\caption{$R_{D_f}(r),r\in[0,3]$ for different $f$-divergences listed in Table \ref{table:divergence_functions}.}
\label{figure:f_radius}
\end{figure}

The loss of efficiency quantified in Theorem \ref{quantification2} is generally hard to compute because it involves a three-layer optimization problem (the Wasserstein distance $W(P,Q)$ is also computed via an optimization problem). Here we only consider the trivial cases $\hat c_{\lk, R_{D_{\lk}}(r)}(x, \hat P)$ and $\hat c_{\chi^2, R_{D_{\chi^2}}(r)}(x, \hat P)$ whose ambiguity sets are the whole space $\mathcal{P}(\Sigma)$ since $R_{D_{\lk}}(r)=R_{D_{\chi^2}}(r)=\infty$. Let $\diam(\Sigma) = \max_{\xi\in \Sigma, \, \xi'\in\Sigma} \norm{\xi-\xi'}$ again denote the diameter of the set $\Sigma$ and suppose that the maximum is attained at some $\bar \xi\in \Sigma$ and $\bar \xi'\in \Sigma$ as $\Sigma$ is compact. Notice that any distribution $Q$ in the KL-divergence ambiguity set $\{Q\in \mc P(\Sigma),~D_{\kl}(\delta_{\bar \xi}, Q) \leq r\}$ must satisfy $Q(\{\bar \xi\})\geq e^{-r}$, which implies that
\[
\inf_{Q:D_{\kl}(\delta_{\bar \xi}, Q) \leq r}W(\delta_{\bar \xi'},Q)\geq \inf_{Q:D_{\kl}(\delta_{\bar \xi}, Q) \leq r}\|\bar \xi'-\bar \xi\|Q(\{\bar \xi\})\geq \diam(\Sigma)e^{-r}.
\]
Hence, by further taking supremum over $\hat P$ in the conclusion of Theorem \ref{quantification2}, we have the following lower bound on the loss of efficiency for $\hat c_{\lk, R_{D_{\lk}}(r)}(x, \hat P)$ (same for $\hat c_{\chi^2, R_{D_{\chi^2}}(r)}(x, \hat P)$)
\begin{align*}
&\sup_{\hat P\in \mc P(\Sigma)}\sup_{\ell\in \mc L_1}\left(\hat{c}_{\lk, R_{D_{\lk}}(r)}(x,\hat{P})-\hat{c}_{\kl,r}(x,\hat{P})\right)\\
=&\sup_{\hat P\in \mc P(\Sigma)}\sup_{P\in\mathcal{P}(\Sigma)}\inf_{Q:D_{\kl}(\hat{P},Q)\leq r}W(P,Q)\\
\geq&\inf_{Q:D_{\kl}(\delta_{\bar \xi},Q)\leq r}W(\delta_{\bar \xi'},Q) ~(\text{by taking }\hat P=\delta_{\bar \xi},P=\delta_{\bar \xi'})
\geq \diam(\Sigma)e^{-r}.
\end{align*}

\end{example}

\begin{example}[Wasserstein]
  We argue here that $R_W(r) \geq \diam(\Sigma) (1-e^{-r})$ for all $r\geq 0$. Suppose $\diam(\Sigma)$ is attained at points $\bar \xi\in \Sigma$ and $\bar \xi'\in \Sigma$. Indeed, $R_W(r) \geq W(\delta_{\bar \xi}, e^{-r}\delta_{\bar \xi} + (1-e^{-r})\delta_{\bar \xi'})=\diam(\Sigma) (1-e^{-r})$ as $D_{\kl}(\delta_{\bar \xi}, e^{-r}\delta_{\bar \xi} + (1-e^{-r}\delta_{\xi'}))=r$. 

  Suppose now that $\Sigma$ is a convex set. We can construct $\bar \xi''=e^{-r}\bar \xi + (1-e^{-r})\bar \xi'\in \Sigma$ so that $\|\bar \xi''-\bar \xi\|=\diam(\Sigma) (1-e^{-r})$ and hence $W(\delta_{\bar \xi}, \delta_{\bar \xi''})=\diam(\Sigma) (1-e^{-r})$. Thus, we obtain $\delta_{\bar \xi''}\in\{P\in \mc P(\Sigma),~W(\delta_{\bar \xi}, P) \leq R_W(r)\}$ since $R_W(r) \geq \diam(\Sigma) (1-e^{-r})$. Moreover, recall that any distribution $Q$ in the KL-divergence ambiguity set $\{Q\in \mc P(\Sigma),~D_{\kl}(\delta_{\bar \xi}, Q) \leq r\}$ must satisfy $Q(\{\bar \xi\})\geq e^{-r}$, which implies that
  \[
  \inf_{Q:D_{\kl}(\delta_{\bar \xi}, Q) \leq r}W(\delta_{\bar \xi''},Q)\geq \inf_{Q:D_{\kl}(\delta_{\bar \xi}, Q) \leq r}\|\bar \xi''-\bar \xi\|Q(\{\bar \xi\})\geq \diam(\Sigma) (1-e^{-r})e^{-r}.
  \]
  Hence, by further taking supremum over $\hat P$ in the conclusion of Theorem \ref{quantification2}, we have the following lower bound on the loss of efficiency for $\hat{c}_{W,R_W(r)}(x,\hat{P})$
  \begin{align*}
      &\sup_{\hat P\in \mc P(\Sigma)}\sup_{\ell\in \mc L_1}\left(\hat{c}_{W,R_W(r)}(x,\hat{P})-\hat{c}_{\kl,r}(x,\hat{P})\right)\\
      =&\sup_{\hat P\in \mc P(\Sigma)}\sup_{P:W(\hat P, P) \leq R_W(r)}\inf_{Q:D_{\kl}(\hat{P},Q)\leq r}W(P,Q)\\
      \geq&\inf_{Q:D_{\kl}(\delta_{\bar \xi},Q)\leq r}W(\delta_{\bar \xi''},Q) ~(\text{by taking }\hat P=\delta_{\bar \xi},P=\delta_{\bar \xi''})
      \geq \diam(\Sigma) (1-e^{-r})e^{-r}.
  \end{align*}
\end{example}

\section{Computation}
\label{sec:computation}

We close our study by presenting how the predictors advanced in (\ref{eq:formulation-type-1}) and (\ref{eq:LPformulation-type-1}) admit tractable reformulations under relatively mild assumptions on the loss function $\ell$.
For the sake of exposition we consider here the distances
\begin{equation}
  \label{def:dist-type-i-gen}
  D^\epsilon_{f, d}(P', P) = \inf\set{D_f(Q, P)}{Q\in \mc P(\Sigma), ~W_d(P', Q)\leq \epsilon}
\end{equation}
and
\begin{equation}
  \label{def:dist-type-ii-gen}
  D^\epsilon_{d, f}(P', P) = \inf\set{D_f(P', Q)}{Q\in \mc P(\Sigma), ~W_d(Q, P) \leq \epsilon}
\end{equation}
where \eqref{def:dist-type-i-gen} generalizes the distances $D^\epsilon_{f, W}(P', P)$ and $D^\epsilon_{f, \lp}(P', P)$ for judiciously chosen lower semicontinous transport cost functions $d$ (see Sections~\ref{ssec:wass-dist} and \ref{ssec:LP-dist}), and \eqref{def:dist-type-ii-gen} considers the parallel notion that exchanges the roles of the arguments in $d$ and $f$ discussed in Section \ref{ssec:smooth-f-diverg}. In particular, Wasserstein smoothing can be reduced to the particular choice $d(\xi, \xi')=\norm{\xi-\xi'}$ whereas LP smoothing is identified with the transport cost $d(\xi, \xi')=\one{\norm{\xi-\xi'}> \epsilon}$.
Likewise, we associate with both distances the generalized predictors
\begin{equation}
  \label{eq:formulation-type-1-gen}
  \hat c_{f, d, r}^{\epsilon}(x, \hat P) \defn
  \left\{
    \begin{array}{rl}
      \sup & \int \ell(x, \xi) \, \d P(\xi)\\
      \st  & P\in \mc P(\Sigma),~D^\epsilon_{f, d}(\hat P, P) \leq r
    \end{array}
  \right.
\end{equation}
and
\begin{equation}
  \label{eq:formulation-type-2-gen}
  \hat c_{d, f, r}^{\epsilon}(x, \hat P) \defn
  \left\{
    \begin{array}{rl}
      \sup & \int \ell(x, \xi) \, \d P(\xi)\\
      \st  & P\in \mc P(\Sigma),~D^\epsilon_{d, f}(\hat P, P) \leq r.
    \end{array}
  \right.
\end{equation}
As our composite predictors can be interpreted to interpolate between a divergence DRO predictor and an optimal transport DRO predictor, we first recall under what circumstances both of these fundamental predictors admit a tractable reformulation.

\subsection{Tractability of General $f$-Divergence DRO}
We consider a general $f$-divergence DRO predictor
\begin{equation}
  \label{eq:f-dro}
  \hat c_{f, r}(x, \hat P) \defn
  \left\{
    \begin{array}{rl}
      \sup & \int \ell(x,\xi) \, \d P(\xi)\\
      \st  & P\in \mc P(\Sigma),\\
           & D_f(\hat P, P)\leq r.
    \end{array}
  \right.
\end{equation}
We remark that as the divergence $D_f$ is convex jointly in both its arguments, this DRO predictor is a concave function of $\hat P\in \mc P(\Sigma)$.
Furthermore, it has a tractable dual representation involving the convex conjugate of $f$ when the set $\Sigma$ is finite (see, e.g., \citep{bayraksan2015data}).
However, to the best of our knowledge this predictor has not been studied in the general case presented in Section \ref{sec:f def}.
We first show that the distribution $P$ in (\ref{eq:f-dro}) over which the supremum is taken can be decomposed into a component $P_c$ absolutely continuous with respect to $\hat P$ and a singular probability $p_s$ with worst-case associated cost $\ell^\infty(x)\defn \sup_{\xi\in \Sigma}\ell(x,\xi)$, which gives the following:

\begin{lemma}
  \label{lemma:finite-representation-f-divergence}
  Let $0 \leq \ell(x, \xi) \leq \ell^\infty(x)\defn\sup_{\xi\in \Sigma}\ell(x,\xi)$.
  Then, we have for all $\hat P\in \mc P(\Sigma)$, $x\in X$ and $r\geq 0$ that
  \begin{equation}
    \label{eq:f-dro-continuous-singular}
    \hat c_{f, r}(x, \hat P) =
    \left\{
      \begin{array}{rl}
        \sup & \int \ell(x,\xi) \, \d P_c(\xi) + p_s \ell^\infty(x)\\
        \st  & P_c\in \mc P_+(\Sigma), ~P_c \ll \hat P, ~p_s \in \Re,\\
             & \int  \, \d P_c(\xi) + p_s=1,\\
             & \int f \left(\tfrac{\d P_c}{\d \hat P}(\xi)\right) \, \d \hat P(\xi) + p_s f^\infty \leq r.
      \end{array}
    \right.
  \end{equation}
  where $\mc P_+(\Sigma)$ denotes the set of positive measures on $\Sigma$.  
\end{lemma}

We are now ready to state a dual representation for the predictor $\hat c_{f, r}$ involving the convex conjugate $f^\star: s \mapsto \sup_{t\geq 0} \, s t - f(t) \in \nRe$ of the function $f$. Lemma \ref{lemma:conjugate-function} indicates that this function is convex lower semicontinuous and nondecreasing. 
We define here the perspective function $\lambda f^\star(\eta/\lambda)$ as
\[
  \lambda f^{\ast}\left(  \frac{\eta}{\lambda}\right)
  \defn
  \begin{cases}
    \sup_{t\geq 0}\eta t-\lambda f(t) & {\rm{if}}~\lambda\geq 0, \\
    +\infty & {\rm{otherwise}}
  \end{cases}
\]
which is lower semicontinuous and jointly convex extended value function. We have the following:

\begin{theorem}[Tractable reformulation of general $f$-divergence DRO]
  \label{thm:love-dual}
  Let $0 \leq \ell(x, \xi) \leq \ell^\infty(x)\defn\sup_{\xi\in \Sigma}\ell(x,\xi)$ for all $x\in X$ and $\xi \in \Sigma$.
  Then, we have
  \begin{align}
    \label{eq:love-dual}
    \hat c_{f, r}(x, \hat P) = & \min \{\textstyle\int \lambda f^\star\left(\tfrac{(\ell(x,\xi)-\eta)}{\lambda}\right) \d \hat P(\xi) + r\lambda +\eta:\\
    & \hspace{12em}\eta\in\Re, ~\lambda\in\Re_+, ~\ell^\infty(x) - \eta \leq f^\infty \lambda \} \nonumber
  \end{align}
  for all $\hat P\in \mc P(\Sigma)$, $x\in X$ and $r> 0$.
\end{theorem}

Theorem \ref{thm:love-dual} states that the divergence DRO in \eqref{eq:f-dro} admits a tractable reformulation in terms of the convex conjugate of the divergence function $f$. In other words, if we have access to an oracle with which we can evaluate $f^\star$ efficiently then the associated DRO can be solved efficiently. Table \ref{table:divergence_functions} lists the associated explicit conjugate functions $f^\star$ for several popular $f$-divergences. Finally, comparing the result in Theorem \ref{thm:love-dual} with the counterpart for the empirical $f$-divergence predictor in \eqref{eq:dro:ef},
\begin{align}
  \tilde c_{f, r}(x, \hat P) & = \textstyle\inf_{\lambda\geq 0, \mu\in \Re}~\mu+r\lambda+\lambda \int f^\star\left(\tfrac{(\ell(x,\xi)-\mu)}{\lambda}\right) \d \hat P(\xi), \label{eq:ahmadi-javid} 
\end{align}
we see that the only difference is the final constraint $\ell^\infty(x) - \eta \leq f^\infty \lambda$.

\subsection{Tractability of Optimal Transport DRO}
\label{ssec:wasserstein-dro-computation}

Consider now a  transport-based DRO problem
\begin{equation}
  \label{eq:W-dro}
  \hat c_{d, \epsilon}(x, \hat P) \defn
  \left\{
    \begin{array}{rl}
      \sup & \int \ell(x,\xi) \, \d P(\xi)\\
      \st  & P\in \mc P(\Sigma),\\
           & W_d(\hat P, P)\leq \epsilon.
    \end{array}
  \right.
\end{equation}
We remark that, like the previous divergence case, as the metric $W_d$ is convex jointly in both its arguments, this DRO predictor is a concave function of $\hat P\in \mc P(\Sigma)$.
The tractability of such Wasserstein DRO problems has been studied extensively since the pioneering work \cite{mohajerin2018data}. A general strong dual formulations can be stated \citep{blanchet2019quantifying} in terms of the Moreau-Yosida regularization
\[
  \ell_\gamma(x, \xi) \defn \sup_{\xi'\in \Sigma} \ell(x, \xi')-\gamma d(\xi, \xi')
\]
of the loss function $\ell$. 
However, for our usage we will find it more opportune to consider a slightly different regularization instead.
Define the $\delta$-inflated loss function as $\ell^\delta(x, \xi) \defn \sup_{\xi'\in \Sigma, d(\xi, \xi')\leq \delta} \ell(x, \xi')$ for which we have that $\ell^\infty(x) = \lim_{\delta\to\infty} \ell^\delta(x, \xi)$ for all $x\in X$ and $\xi\in \Sigma$.
As we assume that the event set $\Sigma$ is compact, we have that $\ell^\infty(x) = \ell^{\diam(\Sigma)}(x, \xi)$ for all $x\in X$ and $\xi\in \Sigma$. The inflated and regularized loss functions are intimately related as indicated by the following result.

\begin{lemma}
  \label{lemma:moreau-yosida}
  We have for any $x\in X$, $\xi\in \Sigma$ and $\gamma\geq 0$ that
  \[
    \ell_\gamma(x, \xi) = \textstyle\sup_{\delta\in \Delta} \ell^\delta(x,\xi) -\gamma \delta
  \]
  with $\Delta \defn [0, \diam(\Sigma)]$.
\end{lemma}

With the above, we can rephrase \cite[Theorem 1.a]{blanchet2019quantifying} slightly in terms of $\ell^\delta(x,\xi)$ as follows:

\begin{theorem}[{\cite[Theorem 1.a]{blanchet2019quantifying}}]
  \label{thm:kuhn-dual}
  Let the distance function $d$ be lower semicontinuous so that $d(\xi, \xi)=0$ for all $\xi\in \Sigma$ and $\ell(x, \xi)$ an upper semicontinuous function in $\xi\in \Sigma$ for all $x\in X$.
  Then, we have 
  \begin{equation}
    \label{eq:kuhn-dual}
    \hat c_{d, \epsilon}(x, \hat P) = \min \set{\textstyle\int \left[\sup_{\delta\in \Delta} \ell^\delta(x,\xi) -\gamma \delta\right] \,\d \hat P(\xi)+\epsilon\gamma}{\gamma\geq 0}
  \end{equation}
  for all $x\in X$, $\hat P\in \mc P(\Sigma)$ and $\epsilon > 0$.
\end{theorem}

However, we remark that due to the univariate quantifier $\delta\in \Delta$ in (\ref{eq:kuhn-dual}) it cannot be solved directly using an off-the-shelf optimization solver just yet.
Instead, we may solve
\begin{equation}
  \label{eq:kuhn-dual-discretized}
  \hat c^K_{d, \epsilon}(x, \hat P) \defn \min \set{\textstyle\int \left[\max_{\delta\in \Delta_K} \ell^\delta(x,\xi) -\gamma \delta\right] \,\d \hat P(\xi)+\epsilon\gamma}{\gamma\geq 0}
\end{equation}
where here we consider a discretized subset $\Delta_K$ of $\Delta$. If we have access to an oracle with which we can evaluate $\ell^\delta$ efficiently then $\hat c_{d, \epsilon}(x, \hat P)$  can be evaluated efficiently to any desired level of relative accuracy as we show in the next proposition. Let
\(
  \Delta_K \defn \{0=\epsilon\cdot \tfrac{0}{K}, \dots, \epsilon = \epsilon\cdot \tfrac{K}{K} = \epsilon \left(\tfrac{\diam(\Sigma)}{\epsilon}\right)^{0/K}, \dots, \diam(\Sigma)= \epsilon \left(\tfrac{\diam(\Sigma)}{\epsilon}\right)^{K/K}  \}\subseteq \Delta
\)
for some fixed $\epsilon\in (0, \diam(\Sigma)]$. We have the following:

\begin{proposition}
  \label{lemma:wasserstein:approximation}
  Let the distance function $d$ be lower semicontinuous so that $d(\xi, \xi)=0$ for all $\xi\in \Sigma$ and $\ell(x, \xi)$ an upper semicontinuous function in $\xi\in \Sigma$ for all $x\in X$.
  Then, we have $$\hat c^K_{d, \epsilon}(x, \hat P)\leq \hat c_{d, \epsilon}(x, \hat P)\leq \left(\left(\tfrac{\diam(\Sigma)}{\epsilon}\right)^{1/K}+\tfrac{1}{K}\right) \hat c^K_{d, \epsilon}(x, \hat P)$$ for all $x\in X$, $\hat P\in \mc P(\Sigma)$ and $\epsilon \in (0, \diam(\Sigma)]$.
\end{proposition}

We remark that $\Delta_K$ is a discrete set counting only $\abs{\Delta_K} = 2K+1$ points. From the previous proposition it follows that with $K\geq 2 \left(1+\log\left(\tfrac{\diam(\Sigma)}{\epsilon}\right)\right)/\rho$ we get that \( \hat c^K_{d, \epsilon}(x, \hat P)\leq \hat c_{d, \epsilon}(x, \hat P)\leq (1+\rho) \hat c^K_{d, \epsilon}(x, \hat P) \). That is, we can approximate $\hat c_{d, \epsilon}(x, \hat P)$ as $\hat c^K_{d, \epsilon}(x, \hat P)$ within any relative accuracy $\rho>0$. %

\begin{example}[Ball Oracles]
  From the previous discussion it follows that for problem (\ref{eq:W-dro}) to be tractable we need to have access to an oracle with which the inflated loss function $\ell^\delta(x, \xi)$ can be evaluated efficiently. Such oracles are particularly well studied in the context of trust region optimization methods \citep{conn2000trust} where they are also known as ball oracles \citep{carmon2020acceleration}.
  Let $\ell(x, \xi)=\max_{t\in [1,\dots, T]}\ell_t(x, \xi)$ be a finite maximum of concave functions $\ell_t(x, \cdot)$, $d(\xi, \cdot)$ a quasiconvex function and $\Xi$ a convex set. When the concave optimization problems $\ell^\delta_t(x, \xi) = \sup \tset{\ell_t(x, \xi')}{\xi'\in \Sigma, \, d(\xi, \xi')\leq \delta}$ admit a ball oracle then clearly so does $\ell^\delta(x, \xi)=\max_{t\in [1,\dots, T]}\ell^\delta_t(x, \xi)$. This structural assumption was imposed in \cite{mohajerin2018data} to guarantee the tractability of the Wasserstein DRO formulation.
\end{example}

\subsection{Tractability of Smoothed DRO}

We now indicate that the composite DRO problems in (\ref{eq:formulation-type-1-gen}) and (\ref{eq:formulation-type-2-gen}) are not computationally harder than the standard divergence and Wasserstein DRO problems discussed in (\ref{eq:f-dro}) and (\ref{eq:W-dro}). First we consider (\ref{eq:formulation-type-1-gen}), which possesses the following reformulation.

\begin{theorem}[Tractable reformulation of smoothed DRO]
  \label{thm:main-type-I}
  Let $0 \leq \ell(x, \xi) \leq \ell^\infty(x)\defn\sup_{\xi\in \Sigma}\ell(x,\xi)$, $d$  a lower semicontinuous function so that $d(\xi, \xi)=0$ for all $\xi\in \Sigma$ and $\ell(x, \xi)$ an upper semicontinuous function in $\xi\in \Sigma$ for all $x\in X$.
  We have
  \[
    \begin{array}{r@{~\,}l}
      \hat c^\epsilon_{f, d, r}(x, \hat P) = \inf & \int \left[\sup_{\delta\in \Delta} \lambda f^\star \left(\tfrac{(\ell^{\delta}(x,\xi)-\eta)}{\lambda}\right)-\delta \gamma\right]  \, \d \hat P(\xi) + r \lambda +\eta +\epsilon\gamma\\
      \st & \eta\in \Re, ~ \lambda\in \Re_+, ~\gamma\in \Re_+, \\
                                                   & \ell^\infty(x)-\eta \leq f^\infty \lambda
    \end{array}
  \]
  for all $x\in X$, $\hat P\in \mc P(\Sigma)$, $r> 0$ and $\epsilon \geq 0$.
\end{theorem}

Theorem \ref{thm:main-type-I} states that the smoothed DRO in (\ref{eq:formulation-type-1-gen}) admits a tractable reformulation in terms of the inflated loss function $\ell^\delta$ and the convex conjugate $f^\star$.
In other words, if we have access to a ball oracle with which we can evaluate the inflated loss function efficiently, then the optimization problem in (\ref{eq:formulation-type-1}) can be solved efficiently.
Indeed, we may approximate $\hat c^\epsilon_{f,d,r}(x, \hat P)$ with
\[
  \begin{array}{rl}
    \hat c^{\epsilon, K}_{f, d, r}(x, \hat P) = \inf & \int \left[\sup_{\delta\in \Delta_K} \lambda f^\star \left(\tfrac{(\ell^{\delta}(x,\xi)-\eta)}{\lambda}\right)-\delta \gamma\right]  \, \d \hat P(\xi) + r \lambda +\eta +\epsilon\gamma \\
    \st & \eta\in \Re, ~ \lambda\in \Re_+, ~\gamma\in \Re_+,\\
                                                  & \ell^\infty(x)+\epsilon\gamma-\eta \leq f^\infty \lambda
  \end{array}
\]
to any desired level of relative accuracy where here we consider the same finite subset $\Delta_K\subseteq \Delta$ defined as in Section \ref{ssec:wasserstein-dro-computation}, as presented below.

\begin{proposition}
  \label{lemma:wasserstein-divergence:approximation}
    Let $0 \leq \ell(x, \xi) \leq \ell^\infty(x)\defn\sup_{\xi\in \Sigma}\ell(x,\xi)$, $f^\star$ an upper-semicontinous function, $d$  a lower semicontinuous function so that $d(\xi, \xi)=0$ for all $\xi\in \Sigma$ and $\ell(x, \xi)$ an upper semicontinuous function in $\xi\in \Sigma$ for all $x\in X$. We have $$\hat c^{\epsilon, K}_{f, d, r}(x, \hat P)\leq  \hat c^{\epsilon}_{f, d, r}(x, \hat P)\leq  \left(\left(\tfrac{\diam(\Sigma)}{\epsilon}\right)^{1/K}+\tfrac{1}{K}\right)  \hat c^{\epsilon, K}_{f, d, r}(x, \hat P)$$ for all $x\in X$, $\hat P\in \mc P(\Sigma)$ and $\epsilon \in (0, \diam(\Sigma)]$.
\end{proposition}

Next we consider the parallel formulation (\ref{eq:formulation-type-2-gen}) where the argument roles of $f$ and $d$ are reversed. Similar to Theorem \ref{thm:main-type-I}, we have the following:

\begin{theorem}[Tractable reformulation of reverse-smoothed DRO]
  \label{thm:main-type-II}
  Let $0 \leq \ell(x, \xi) \leq \ell^\infty(x)\defn\sup_{\xi\in \Sigma}\ell(x,\xi)<\infty$, $d$  a lower semicontinuous function so that $d(\xi, \xi)=0$ for all $\xi\in \Sigma$ and $\ell(x, \xi)$ an upper semicontinuous function in $\xi\in \Sigma$ for all $x\in X$. We have
  \[
    \begin{array}{rl}
      \hat c^{\epsilon}_{d, f, r}(x, \hat P) = \inf & \int \lambda f^\star\left(\tfrac{\left(\left[\sup_{\delta\in \Delta}\ell^{\delta}(x,\xi)-\gamma\delta\right] +\epsilon\gamma-\eta\right)}{\lambda}\right) \, \d \hat P(\xi) + r \lambda +\eta \\
      \st & \eta\in \Re, ~ \lambda\in \Re_+, ~\gamma\in \Re_+,\\
                                                    & \ell^\infty(x)+\epsilon\gamma-\eta \leq f^\infty \lambda
    \end{array}
  \]
  for all $x\in X$, $\hat P\in \mc P(\Sigma)$, $r\geq 0$ and $\epsilon > 0$.
\end{theorem}

The implication of Theorem \ref{thm:main-type-II} is the same as that of Theorem  \ref{thm:main-type-I}. Similar to before, we can approximate $\hat c^\epsilon_{d,f,r}(x, \hat P)$ with
\[
  \begin{array}{rl}
    \hat c^{\epsilon, K}_{d, f, r}(x, \hat P) = \inf & \int \lambda f^\star\left(\tfrac{\left(\left[\sup_{\delta\in \Delta_K}\ell^{\delta}(x,\xi)-\gamma\delta\right] +\epsilon\gamma-\eta\right)}{\lambda}\right) \, \d \hat P(\xi) + r \lambda +\eta \\
    \st & \eta\in \Re, ~ \lambda\in \Re_+, ~\gamma\in \Re_+,\\
                                                  & \ell^\infty(x)+\epsilon\gamma-\eta \leq f^\infty \lambda
  \end{array}
\]
to any desired level of accuracy, as stated below.

\begin{proposition}
  \label{lemma:divergence-wasserstein:approximation}
  Let $0 \leq \ell(x, \xi) \leq \ell^\infty(x)\defn\sup_{\xi\in \Sigma}\ell(x,\xi)$, $d$  a lower semicontinuous function so that $d(\xi, \xi)=0$ for all $\xi\in \Sigma$ and $\ell(x, \xi)$ an upper semicontinuous function in $\xi\in \Sigma$ for all $x\in X$. We have $$\hat c^{\epsilon, K}_{d, f, r}(x, \hat P)\leq  \hat c^{\epsilon}_{d, f, r}(x, \hat P)\leq  \left(\left(\tfrac{\diam(\Sigma)}{\epsilon}\right)^{1/K}+\tfrac{1}{K}\right)  \hat c^{\epsilon, K}_{d, f, r}(x, \hat P)$$ for all $x\in X$, $\hat P\in \mc P(\Sigma)$ and $\epsilon \in (0, \diam(\Sigma)]$.
\end{proposition}

Therefore, comparing Theorems \ref{thm:main-type-I} and \ref{thm:main-type-II} with Theorems \ref{thm:love-dual} and \ref{thm:kuhn-dual} (and similarly Propositions \ref{lemma:wasserstein-divergence:approximation} and \ref{lemma:divergence-wasserstein:approximation} with Proposition \ref{lemma:wasserstein:approximation}), we see that, regardless of the smoothed DRO formulations focused on in this paper or their reverse versions, their computation is not much harder than the counterparts in standard $f$-divergence and Wasserstein DRO problems.

\section*{Acknowledgements}
We gratefully acknowledge support from the National Science Foundation under grants CAREER CMMI-1834710 and IIS-1849280.

\bibliography{references}

\appendix

\section{Supporting Lemmas}
\label{sec:supporting_lemmas}

\begin{lemma}[Continuity and Monotonicity]
  \label{lemma:conjugate-function}
  The function $f^\star$ is convex and satisfies $s_1\leq s_2 \implies f^\star(s_1)\leq f^\star(s_2)$ and $f^\infty < s$ implies $f^\star(s)=\infty$. Finally, $f^\star$ is continuous on $[-\infty, f^\infty]$.
\end{lemma}

\begin{lemma}[Recession Inequality]
  \label{lemma:recession-inequality}
  Let $f$ be a convex function. Then,
  \(
  f((q+\Delta)/p') p' \leq f(q/p') p' + \Delta f^\infty
  \)
  for all  $q\in \Re_+$, $p'\in \Re_+$ and $\Delta\geq 0$.
\end{lemma}

The following Lemmas \ref{property_smooth_f_div}-\ref{usc_smooth_f_DRO} deal with the properties of the smoothed $f$-divergences and the corresponding DRO predictors.
In these lemmas, we assume the space $\Sigma$ is a general compact set in $\mathbb{R}^d$ which implies that $\mathcal{P}(\Sigma)$ is a compact Polish space according to \cite[Theorem D.8]{dembo2009large}. Moreover, we assume the loss function $\ell(x,\xi)$ is upper semicontinuous in $\xi$ for any $x$ (this essentially imposes no restrictions on $\ell$ if $\Sigma$ is finite) In particular, the results in Lemmas \ref{attainability_smooth_f_DRO}-\ref{usc_smooth_f_DRO} apply to the usual $f$-divergence DRO predictor $\hat{c}_{f,r}(x,\hat P)$ when $\epsilon=0$.

\begin{lemma}
  \label{property_smooth_f_div}
  Assume $\epsilon\geq0$. The infimum in the definitions of $D_{f,W}^{\epsilon}(P^{\prime},P)$ and $D_{f,\lp}^{\epsilon}(P^{\prime},P)$ is attainable. Moreover, $D_{f,W}^{\epsilon}(P^{\prime},P)$ and $D_{f,\lp}^{\epsilon}(P^{\prime},P)$ are convex and lower semicontinuous in $(P^{\prime},P)$ on $\mathcal{P}(\Sigma)\times\mathcal{P}(\Sigma)$ and satisfy the limiting property
\[
\lim_{\epsilon\downarrow0}D_{f,W}^{\epsilon}(P^{\prime},P)=\lim_{\epsilon\downarrow0}D_{f,\lp}^{\epsilon}(P^{\prime},P)=D_{f}(P^{\prime},P),
\]
for any $(P^{\prime},P)\in\mathcal{P}(\Sigma)\times\mathcal{P}(\Sigma)$.
\end{lemma}

\begin{lemma}
  \label{attainability_smooth_f_DRO}
  Suppose the loss function $\ell(x,\xi)$ is upper semicontinous in $\xi$ for any $x\in X$. Assume $\epsilon\geq0$ and $r\geq0$. Then for any $x\in X,\hat{P}\in\mathcal{P}(\Sigma)$, the supremum in the definitions of $\hat{c}_{f,W,r}^{\epsilon}(x,\hat{P})$ and $\hat{c}_{f,\lp,r}^{\epsilon}(x,\hat{P})$ is attainable.
\end{lemma}

\begin{lemma}
  \label{usc_smooth_f_DRO}
  Suppose the loss function $\ell(x,\xi)$ is upper semicontinous in $\xi$ for any $x\in X$. Assume $\epsilon\geq0$ and $r\geq0$. Then the predictors $\hat{c}_{f,W,r}^{\epsilon}(x,\hat{P})$ and $\hat{c}_{f,\lp,r}^{\epsilon}(x,\hat{P})$ are upper semicontinuous in $\hat{P}$ for any $x\in X$.
\end{lemma}

The next lemma is a generalized version of the Riemann-Lebesgue lemma for periodic functions. The proof can be found in \cite[Corollary 2.2]{andrica2004extension}.

\begin{lemma}[Generalized Riemann-Lebesgue Lemma]\label{RL_lemma}
Suppose $f$ is a continuous function on $[a,b]\subseteq[0,\infty)$ and $g$ is a continuous periodic function of period $T$ on $[0,\infty)$. We have%
\[
\lim_{n\rightarrow\infty}\int_{a}^{b}f(x)g(nx)dx=\frac{1}{T}\int_{0}^{T}g(x)dx\int_{a}^{b}f(x)dx.
\]
\end{lemma}

\begin{lemma}\label{event_measurability}
Let $\Sigma$ be a compact set and $\mathcal{C}_{u}(X\times\mathcal{P}(\Sigma))$ be the class of upper semicontinuous predictors. Given a loss function $\ell$, if a predictor $\hat{c}\in\mathcal{C}_{u}(X\times\mathcal{P}(\Sigma))$, then the events $\{\mathbb{E}_{P}[\ell(x,\xi)]>\hat{c}(x,P_{N})\}$ for any $x\in X$ and $\{\exists x\in X~\st~\mathbb{E}_{P}[\ell(x,\xi)]>\hat{c}(x,P_{N})\}$ are measurable for any $P\in\mathcal{P}(\Sigma)$.
\end{lemma}

\section{Large Deviations Theory}\label{sec:LDP_theory}
In this section, we present necessary results in large deviations theory \cite{dembo2009large}. We focus on the Sanov's theorem which provides bounds on the exponential decaying rate of $(1/N)\log \Prob{P_{N}\in\Gamma}$ for $\Gamma\in\mathcal{P}(\Sigma)$.

\subsection{Finite Support}
First, we consider the finite support setting. Suppose the true distribution $P$ is supported on a finite set $\Sigma$ with cardinality $2\leq|\Sigma|<\infty$, where the probability simplex $\mathcal{P}(\Sigma)$ can be defined as $\mathcal{P}(\Sigma)=\{P\in\mathbb{R}^{|\Sigma|}:\sum_{i\in\Sigma}P(i)=1,P(i)\geq0~\forall i\in\Sigma\}$ equipped with the induced topology of $\mathbb{R}^{|\Sigma|}$. We have the following large deviations principle for the empirical distribution $P_N$.
\begin{theorem}[Sanov's Theorem I]\label{finite_LDP}
Let $P_N$ be the empirical distribution of $N$ i.i.d. samples from the true distribution $P$ supported on a finite set $\Sigma$. For every set $\Gamma\subseteq\mathcal{P}(\Sigma)$, $P_N$ satisfies the large deviations principle
\begin{align*}
  -\inf_{\nu\in\Gamma^{\circ}}D_{\kl}(\nu,P)\leq& \liminf_{N\rightarrow\infty}\frac{1}{N}\log \Prob{P_{N}\in\Gamma}\\
  & \limsup_{N\rightarrow\infty}\frac{1}{N}\log \Prob{P_{N}\in\Gamma}\leq-\inf_{\nu\in\Gamma}D_{\kl}(\nu,P),
\end{align*}
where $\Gamma^{\circ}$ is the interior of $\Gamma$ with respect to the induced topology on $\mathcal{P}(\Sigma)$.
\end{theorem}

\subsection{Infinite Support}\label{sec:infinite_LDP}
In this section, we introduce the large deviations theory for the empirical distribution $P_N$ when the true distribution $P$ is supported on a general compact set $\Sigma\subset\mathbb{R}^d$. 

As we have seen in the finite support case, the large deviations rate depends on the topology on $\mathcal{P}(\Sigma)$, which we will introduce now. Let $\mathcal{M}(\Sigma)$ be the space of finite signed measure on $\Sigma$. Let $C_{b}(\Sigma)\equiv C(\Sigma)$ be the space of (bounded) continuous function equipped with the sup norm. We equip $\mathcal{M}(\Sigma)$ with the weak topology generated by the open sets $\{U_{\phi,y,\delta}:\phi\in C_{b}(\Sigma),y\in\mathbb{R},\delta>0\}$ where
\[
U_{\phi,y,\delta}=\left\{  \nu\in\mathcal{M}(\Sigma):\textstyle\left\vert \int_{\Sigma}\phi d\nu-y\right\vert <\delta\right\}  .
\]
According to \cite[Section 6.2]{dembo2009large}, this topology makes $\mathcal{M}(\Sigma)$ into a locally convex Hausdorff topological vector space, whose topological dual is the class of functionals $\{\nu\mapsto\int_{\Sigma}\phi d\nu:\phi\in C_{b}(\Sigma)\}$ and thus can be identified with $C_{b}(\Sigma)$. Under the weak topology, $\mathcal{P}(\Sigma)$ is a closed subset of $\mathcal{M}(\Sigma)$ and further is a compact subset due to the compactness of $\Sigma$. We equip $\mathcal{P}(\Sigma)$ with the induced topology which is proved to be compatible with the L\'{e}vy-Prokhorov metric. With the L\'{e}vy-Prokhorov metric, $\mathcal{P}(\Sigma)$ is a Polish space and the convergence in $\mathcal{P}(\Sigma)$ is indeed the weak convergence of probability distributions. Moreover, according to \cite[Theorem 2]{gibbs2002choosing}, the L\'{e}vy-Prokhorov metric and Wasserstein metric bound each other by
\begin{equation}
   \LP{P}{Q}^2\leq W(P,Q)\leq (\mathrm{diam}(\Sigma)+1)\LP{P}{Q},\label{ineq_LP_W} 
\end{equation}
where $\mathrm{diam}(\Sigma)=\sup\{||x-y||:x,y\in\Sigma\}$ is the diameter of $\Sigma$. Inequality (\ref{ineq_LP_W}) implies that L\'{e}vy-Prokhorov metric and Wasserstein metric are equivalent metrics on $\mathcal{P}(\Sigma)$ and thus they determine the same topology on $\mathcal{P}(\Sigma)$. Therefore, all the results above remain unchanged if the L\'{e}vy-Prokhorov metric is replaced by the Wasserstein metric.

Below is the large deviations principle when $\Sigma$ is a compact set.
\begin{theorem}[Sanov's Theorem II]\label{infinite_LDP}
Let $P_N$ be the empirical distribution of $N$ i.i.d. samples from the true distribution $P$ supported on a compact set $\Sigma\subset\mathbb{R}^d$. For every Borel subset $\Gamma\subseteq\mathcal{P}(\Sigma)$, $P_N$ satisfies the large deviations principle
\begin{align*}
  -\inf_{\nu\in\Gamma^{\circ}}D_{\kl}(\nu,P)\leq &\liminf_{N\rightarrow\infty}\frac{1}{N}\log \Prob{P_{N}\in\Gamma}\\
                                                 & \limsup_{N\rightarrow\infty}\frac{1}{N}\log \Prob{P_{N}\in\Gamma}\leq-\inf_{\nu\in\bar\Gamma}D_{\kl}(\nu,P),
\end{align*}
where $\Gamma^{\circ}$ is the interior of $\Gamma$ and $\bar\Gamma$ is the closure of $\Gamma$.
\end{theorem}
The upper bound in Theorem \ref{infinite_LDP} can be further strengthened to a finite sample bound by using our smoothed KL-divergence. Let $C(\epsilon,\Sigma,||\cdot||)$ be the minimal covering number of $\Sigma$ under the Euclidean norm $||\cdot||$ using balls of radius $\epsilon$. We have
\begin{theorem}[Finite-sample upper bound]\label{thm:sanov-smooth}
Let $P_{N}$ be the empirical distribution of $N$ i.i.d. samples from the true distribution $P$ supported on a compact set $\Sigma\subset\mathbb{R}^d$. For any Borel set $\Gamma\subseteq\mathcal{P}(\Sigma)$, we have%
\[
\Prob{P_{N}\in\Gamma}\leq\inf_{\epsilon>0}\left\{  \left(  \frac{8}{\epsilon}\right)  ^{C(\epsilon/2,\Sigma,||\cdot||)}\exp\left(  -N\inf_{\nu\in\Gamma}D_{\kl,\lp}^{\epsilon}(\nu,P)\right)  \right\}  ,
\]
and%
\[
\Prob{P_{N}\!\in\!\Gamma}\leq\inf_{\epsilon>0}\left\{  \!\left( \! \frac{8(\mathrm{diam}(\Sigma)+1)}{\epsilon}\right)  ^{C(\frac{\epsilon}{(2(\mathrm{diam}(\Sigma)+1)},\Sigma,||\cdot||)} \!\!\!\!\!\!\exp\left(  -N\inf_{\nu\in\Gamma}D_{\kl,W}^{\epsilon}(\nu,P)\right) \! \right\}\!.
\]
Consequently, we have%
\[
\limsup_{N\rightarrow\infty}\frac{1}{N}\log \Prob{P_{N}\in\Gamma}\leq\min\left\{-\inf_{\nu\in\Gamma}D_{\kl,\lp}^{\epsilon}(\nu,P),-\inf_{\nu\in\Gamma}D_{\kl,W}^{\epsilon}(\nu,P)\right\}
\]
for any $\epsilon>0$.
\end{theorem}

By comparing Theorems \ref{infinite_LDP} and \ref{thm:sanov-smooth}, we can see that the advantage of using the smoothed KL-divergence is that the closure $\bar\Gamma$ in large deviations upper bound is replaced by the original set $\Gamma$. Since $\bar\Gamma$ can be significantly larger than $\Gamma$, this improvement is crucial to show the uniform feasibility of the smoothed KL-DRO predictor.

\section{Proofs}

\subsection{Proofs for Section \ref{sec:supporting_lemmas}}

\subsubsection{Proof of Lemma \ref{lemma:conjugate-function}}
\begin{proof}

  The monotonicity property follows immediately from
  \[
    s_1 t - f(t) \leq s_2 t -f(t) \quad \forall t\in \Re_+
  \]
  and taking the infimum over $t\geq 0$. 
  Convexity and lower semicontinuity of $f^\star$ is immediate from the fact that $f^\star(s)=\sup_{t\geq 0} s t -f (t)$ is the maximum of linear and continuous functions in $s$.
  
  Recall that as $f$ is a convex function $t\mapsto f(t)/t$ is an increasing function. Hence, $\sup_{t\geq 0}f(t)/t = \lim_{t\to\infty}f(t)/t=f^\infty$ and consequently we have $f^\infty t-f(t)\geq 0$ for all $t\in\Re_+$. When $f^\infty$ is finite and $f^\infty < s$ we have
  \[
    \textstyle\sup_{t\geq 0} \, s t - f(t) = \sup_{t\geq 0} \, (f^\infty t - f(t))+(s-f^\infty) t \geq \sup_{t\geq 0}\, (s-f^\infty) t = +\infty.
  \]
  Furthermore, for any $s<f^\infty$ it follows that
  \[
    f(s) = -\inf_{t\geq 0} f(t) - st = s t^\star -f(t^\star) \in \Re
  \]
  for some $t^\star\geq 0$ \cite[Proposition 3.2.1]{bertsekas2009convex}. Indeed, the function $t\mapsto f(t)-st$ is lower semicontinous coercive as $\lim_{t\to\infty} f(t)-st = \lim_{t\to\infty} t\left(\tfrac{f(t)}{t}-s\right) = \infty$ as $\lim_{t\to\infty} \tfrac{f(t)}{t}-s > 0$. We have shown that $\dom(f^\star)=\set{f(s)\in \Re}{s\in \Re}$ is of the form $[-\infty, f^\infty]$ or $[-\infty, f^\infty)$. We now show continuity of $f^\star$ on $[0, f^\infty]$ in all cases,
  \begin{itemize}
  \item Case $\dom(f^\star)=\Re$ as $f^\infty=\infty$: $f^\star$ is continuous on $\Re$ as every convex function with domain $\Re$ is continuous on $\Re$ \cite[Corollary 10.1.1]{rockafellar1997convex},
  \item Case $\dom(f^\star)=[-\infty, f^\infty)$ with $f^\infty\in \Re$: Every convex function with open domain $[-\infty, f^\infty)$ for $f^\infty\in \Re$ is continuous on its domain \cite[Theorem 10.1]{rockafellar1997convex}. As we have for any lower semicontinuous convex function that $\lim_{s\to f^\infty}f^\star(s)=\infty = f^\star(f^\infty)$ \cite[Theorem 4.14]{nesterov2003introductory}, the function $f^\star$ is in fact continuous on $[0, f^\infty]$.
  \item Case $\dom(f^\star)=[-\infty, f^\infty]$ with $f^\infty\in \Re$: $f^\star$ is continuous as every lower semicontinuous convex function with polyhedral domain $[-\infty, f^\infty]$ is continuous relative to its domain $[-\infty, f^\infty]$ \cite[Theorem 10.2]{rockafellar1997convex},
  \end{itemize}
  establishing the second claim.
\end{proof}

\subsubsection{Proof of Lemma \ref{lemma:recession-inequality}}

\begin{proof}
  Indeed, consider any fixed $q\in \Re_+$, $p'\in \Re_+$ and denote with $g(\Delta) = f((q+\Delta)/p') p'$.
  The subgradient inequality for convex functions guarantees that
  \(
  g(0) \geq g(\Delta) - \Delta g'(\Delta) \iff f(q/p') p' \geq f((q+\Delta)/p') p' - \Delta g'(\Delta)
  \)
  where $g'(\Delta)$ denotes any subgradient of the function $g$ at $\Delta$.
  Remark however that $g'(\Delta) = f'((q+\Delta)/p') \leq f^\infty$ where here $f'((q+\Delta)/p')$ denotes an arbitrary subgradient of $f$ at $(q+\Delta)/p'$.
\end{proof}

\subsubsection{Proof of Lemma \ref{property_smooth_f_div}}
\begin{proof}
We prove the lemma for $D_{f,\lp}^{\epsilon}(P^{\prime},P)$ and the same argument can be applied to the Wasserstein case.

We first prove the infimum in the definition of $D_{f,\lp}^{\epsilon}(P^{\prime},P)$ is attainable. Note that $D_{f}(Q,P)$ is lower semicontinuous (see e.g. \cite[Theorem 2.34]{ambrosio2000functions}) and the closed subset $\{Q\in\mathcal{P}(\Sigma):\mathrm{LP}(P^{\prime},Q)\leq\epsilon\}\subseteq\mathcal{P}(\Sigma)$ is compact because $\mathcal{P}(\Sigma)$ is compact. Since the infimum of a lower semicontinuous function on a compact set is attainable, we know there exists $Q\in\mathcal{P}(\Sigma)$ such that $\mathrm{LP}(P^{\prime},Q)\leq\epsilon$ and $D_{f,\lp}^{\epsilon}(P^{\prime},P)=D_{f}(Q,P)$.

Next we prove the convexity. Suppose $P_{i}^{\prime},P_{i}\in\mathcal{P}(\Sigma),i=1,2$ and $\lambda\in(0,1)$. By the attainability, there exist $Q_{i}\in\mathcal{P}(\Sigma),i=1,2$ such that%
\[
D_{f,\lp}^{\epsilon}(P_{i}^{\prime},P_{i})=D_{f}(Q_{i},P_{i}),\quad \mathrm{LP}(P_{i}^{\prime},Q_{i})\leq\epsilon
\]
for $i=1,2$. Since L\'{e}vy-Prokhorov distance is quasi-convex, we have
\[
\mathrm{LP}(\lambda P_{1}^{\prime}+(1-\lambda)P_{2}^{\prime},\lambda Q_{1}+(1-\lambda)Q_{2})\leq\max\{\mathrm{LP}(P_{1}^{\prime},Q_{1}),\mathrm{LP}(P_{2}^{\prime},Q_{2})\}\leq\epsilon.
\]
Therefore, by the definition of $D_{f,\lp}^{\epsilon}$ and convexity of $D_{f}$, we have
\begin{align*}
&  D_{f,\lp}^{\epsilon}(\lambda P_{1}^{\prime}+(1-\lambda)P_{2}^{\prime},\lambda P_{1}+(1-\lambda)P_{2})\\
&  \leq D_{f}(\lambda Q_{1}+(1-\lambda)Q_{2},\lambda P_{1}+(1-\lambda)P_{2})\\
&  \leq\lambda D_{f}(Q_{1},P_{1})+(1-\lambda)D_{f}(Q_{2},P_{2})\\
&  =\lambda D_{f,\lp}^{\epsilon}(P_{1}^{\prime},P_{1})+(1-\lambda)D_{f,\lp}^{\epsilon}(P_{2}^{\prime},P_{2}).
\end{align*}

Then we show the lower semicontinuity. Suppose $P_{n}^{\prime},P_{n}\in\mathcal{P}(\Sigma),n\geq1$ which satisfy $P_{n}^{\prime}\overset{w}{\rightarrow}P^{\prime}\in\mathcal{P}(\Sigma)$ and $P_{n}\overset{w}{\rightarrow}P\in\mathcal{P}(\Sigma)$. Let $Q_{n}\in\mathcal{P}(\Sigma)$ be the probability measure attaining the infimum in $D_{f,\lp}^{\epsilon}(P_{n}^{\prime},P_{n})$, i.e.,
\begin{equation}
D_{f,\lp}^{\epsilon}(P_{n}^{\prime},P_{n})=D_{f}(Q_{n},P_{n}),\quad \mathrm{LP}(P_{n}^{\prime},Q_{n})\leq\epsilon\label{solutions _Qn}%
\end{equation}
for any $n\geq1$. To show $D_{f,\lp}^{\epsilon}(P^{\prime},P)\leq\liminf_{n\rightarrow\infty}D_{f,\lp}^{\epsilon}(P_{n}^{\prime},P_{n})$, it suffices to show that $D_{f,\lp}^{\epsilon}(P^{\prime},P)\leq\lim_{k\rightarrow\infty}D_{f,\lp}^{\epsilon}(P_{n_{k}}^{\prime},P_{n_{k}})$ for any converging subsequence $\{D_{f,\lp}^{\epsilon}(P_{n_{k}}^{\prime},P_{n_{k}}),k\geq1\}$. We take such a converging subsequence $D_{f,\lp}^{\epsilon}(P_{n_{k}}^{\prime},P_{n_{k}})=D_{f}(Q_{n_{k}},P_{n_{k}})$. Since $\{Q_{n_{k}},k\geq 1\}\subset\mathcal{P}(\Sigma)$ and $\mathcal{P}(\Sigma)$\ is compact, we can extract a converging subsequence of $\{Q_{n_{k}},k\geq1\}$. For simiplicity of the notation, we still write the converging subsequence as $\{Q_{n_{k}},k\geq1\}$, that is, $Q_{n_{k}}\overset{w}{\rightarrow}Q\in\mathcal{P}(\Sigma)$. Since $\mathrm{LP}(P_{n_{k}}^{\prime},Q_{n_{k}})\leq\epsilon$, $P_{n_{k}}^{\prime}\overset{w}{\rightarrow}P^{\prime}$ and $Q_{n_{k}}\overset{w}{\rightarrow}Q$, we know $\mathrm{LP}(P^{\prime},Q)\leq\epsilon$. Thus, we have%
\[
D_{f,\lp}^{\epsilon}(P^{\prime},P)\leq D_{f}(Q,P)\leq\lim_{k\rightarrow\infty}D_{f}(Q_{n_{k}},P_{n_{k}})=\lim_{k\rightarrow\infty}D_{f,\lp}^{\epsilon}(P_{n_{k}}^{\prime},P_{n_{k}}),
\]
where the first inequality is due to the definition of $D_{f,\lp}^{\epsilon}$, the second inequality is due to the lower semicontinuity of $D_{f}$ and the last equality follows from (\ref{solutions _Qn}). Since the above inequality holds for any converging subsequence $D_{f,\lp}^{\epsilon}(P_{n_{k}}^{\prime},P_{n_{k}})$, we obtain the lower semicontinuity of $D_{f,\lp}^{\epsilon}$.

Finally, we prove the limiting property%
\[
\lim_{\epsilon\downarrow0}D_{f,\lp}^{\epsilon}(P^{\prime},P)=D_{f}(P^{\prime},P).
\]
By the definition of $D_{f,\lp}^{\epsilon}(P^{\prime},P)$, we know that $D_{f,\lp}^{\epsilon}(P^{\prime},P)\leq D_{f}(P^{\prime},P)$ and $D_{f,\lp}^{\epsilon}(P^{\prime},P)$ is non-decreasing as $\epsilon$ decreases. So the limit exists and satisfies%
\[
\lim_{\epsilon\downarrow0}D_{f,\lp}^{\epsilon}(P^{\prime},P)\leq D_{f}(P^{\prime},P).
\]
Now we take a sequence $\epsilon_{n}\downarrow0$. By the attainability, there exists $Q_{n}\in\mathcal{P}(\Sigma)$ such that $D_{f,\lp}^{\epsilon_{n}}(P^{\prime},P)=D_{f}(Q_{n},P)$ and $\mathrm{LP}(P^{\prime},Q_{n})\leq \epsilon_{n}\rightarrow0$. Therefore, we have $Q_{n}\overset{w}{\rightarrow}P^{\prime}$. By the lower semicontinuity of $D_{f}$, we have%
\[
\lim_{\epsilon\downarrow0}D_{f,\lp}^{\epsilon}(P^{\prime},P)=\lim_{n\rightarrow\infty}D_{f,\lp}^{\epsilon_{n}}(P^{\prime},P)=\lim_{n\rightarrow\infty}D_{f}(Q_{n},P)\geq D_{f}(P^{\prime},P).
\]
Combining both directions, we obtain%
\[
\lim_{\epsilon\downarrow0}D_{f,\lp}^{\epsilon}(P^{\prime},P)=D_{f}(P^{\prime},P).
\]
\end{proof}

\subsubsection{Proof of Lemma \ref{attainability_smooth_f_DRO}}
\begin{proof}
We prove this lemma for $\hat{c}_{f,\lp,r}^{\epsilon}(x,\hat{P})$ and the same argument can be applied to the Wasserstein case. By the lower semicontinuity of $D_{f,\lp}^{\epsilon}$, $\{P\in\mathcal{P}(\Sigma):D_{f,\lp}^{\epsilon}(\hat{P},P)\leq r\}$ is closed and thus compact. Moreover, since $\ell(x,\xi)$ is upper semicontinuous in $\xi$ on a compact set $\Sigma$ for any $x$, it is bounded above for any $x$. Therefore, Portmanteau theorem ensures that $c(x,P)=\mathbb{E}_{P}[\ell(x,\xi)]$ is upper semicontinuous in $P$ for any $x$. Then the conclusion follows from the fact that the maximization problem for an upper semicontinuous function on a compact set is attainable.
\end{proof}

\subsubsection{Proof of Lemma \ref{usc_smooth_f_DRO}}
\begin{proof}
We prove this lemma for $\hat{c}_{f,\lp,r}^{\epsilon}(x,\hat{P})$ and the same argument can be applied to the Wasserstein case. We need to show for any $\hat{P}_{n}\overset{w}{\rightarrow}\hat{P}$, we have%
\[
\limsup_{n\rightarrow\infty}\hat{c}_{f,\lp,r}^{\epsilon}(x,\hat{P}_{n})\leq\hat{c}_{f,\lp,r}^{\epsilon}(x,\hat{P}).
\]
By Lemma \ref{attainability_smooth_f_DRO}, there exists $P_{n}$ s.t. $D_{f,\lp}^{\epsilon}(\hat{P}_{n},P_{n})\leq r$ and $\hat{c}_{f,\lp,r}^{\epsilon}(x,\hat{P}_{n})=c(x,P_{n})$. Therefore, it suffices to show%
\begin{equation}
\limsup_{n\rightarrow\infty}c(x,P_{n})\leq\hat{c}_{f,\lp,r}^{\epsilon}(x,\hat{P}).\label{limsup_inequality}
\end{equation}
We take any converging subsequence of $\{c(x,P_{n}),n\geq1\}$, say $\{c(x,P_{n_{k}}),k\geq1\}$. Since $\mathcal{P}(\Sigma)$ is compact, we can extract a converging subsequence of $\{P_{n_{k}}\}$. Without loss of generality, we still denote the converging subsequence as $\{P_{n_{k}}\}$, i.e., $P_{n_{k}}\overset{w}{\rightarrow}P$ for some $P\in\mathcal{P}(\Sigma)$. Since $\hat{P}_{n}\overset{w}{\rightarrow}\hat{P}$, we still have $\hat{P}_{n_{k}}\overset{w}{\rightarrow}\hat{P}$. Therefore, by the lower semicontinuity of $D_{f,\lp}^{\epsilon}$,
\[
D_{f,\lp}^{\epsilon}(\hat{P},P)\leq\liminf_{k\rightarrow\infty}D_{f,\lp}^{\epsilon}(\hat{P}_{n_{k}},P_{n_{k}})\leq r,
\]
and consequently by the upper semicontinuity of $c$ (see the proof of Lemma \ref{attainability_smooth_f_DRO}) and the definition of $\hat{c}_{f,\lp,r}^{\epsilon}(x,\hat{P})$,
\begin{equation}
\lim_{k\rightarrow\infty}c(x,P_{n_{k}})\leq c(x,P)\leq\hat{c}_{f,\lp,r}^{\epsilon}(x,\hat{P}).\label{subsequence_inequality}
\end{equation}
Since (\ref{subsequence_inequality}) holds for any converging subsequence of $\{c(x,P_{n_{k}})\}$, we obtain (\ref{limsup_inequality}). This concludes our proof.
\end{proof}

\subsubsection{Proof of Lemma \ref{event_measurability}}
\begin{proof}
Recall that the topology in $\mathcal{P}(\Sigma)$ is the weak topology (when $\Sigma$ is finite, weak topology is equivalent to the induced Euclidean topology in the simplex $\mathcal{P}(\Sigma)\subseteq\mathbb{R}^{|\Sigma|}$). For any given $x\in X$, we have%
\[
\{\mathbb{E}_{P}[\ell(x,\xi)]>\hat{c}(x,P_{N})\}=\{P_{N}\in\Gamma(x)\},
\]
where%
\[
\Gamma(x)=\{\nu\in\mathcal{P}(\Sigma):\mathbb{E}_{P}[\ell(x,\xi)]>\hat{c}(x,\nu)\}.
\]
Since $\hat{c}(x,\nu)$ is upper semicontinuous in $\nu$ for any $x$, $\Gamma(x)$ is an open subset of $\mathcal{P}(\Sigma)$, which makes $\{\mathbb{E}_{P}[\ell(x,\xi)]>\hat{c}(x,P_{N})\}=\{P_{N}\in\Gamma(x)\}$ measurable. Moreover, we can see
\[
\{\exists x\in X\text{ s.t. }\mathbb{E}_{P}[\ell(x,\xi)]>\hat{c}(x,P_{N})\}=\left\{  \hat{P}_{N}\in\cup_{x\in X}\Gamma(x)\right\}  ,
\]
where $\cup_{x\in X}\Gamma(x)$ is still an open subset of $\mathcal{P}(\Sigma)$. Therefore, $\{\exists x\in X$ s.t. $\mathbb{E}_{P}[\ell(x,\xi)]>\hat{c}(x,P_{N})\}$ is also measurable.
\end{proof}

\subsection{Proofs in Section \ref{sec:LDP_theory}}

\subsubsection{Proof of Theorem \ref{thm:sanov-smooth}}
\begin{proof}
We first consider the L\'{e}vy-Prokhorov smoothing distance $D_{\kl,\lp}^{\epsilon}$. According to \cite[Exercise 6.2.19]{dembo2009large}, we have%
\begin{equation}
\Prob{P_{N}\in\Gamma}\leq\inf_{\epsilon>0}\left\{  \left(  \frac{4}{\epsilon}\right)  ^{C(\epsilon,\Sigma,||\cdot||)}\exp\left(  -N\inf_{\nu\in\Gamma^{\epsilon}}D_{\kl}(\nu,P)\right)  \right\}  ,\label{dembo_finite_bound}%
\end{equation}
where $\Gamma^{\epsilon}$ is the closed $\epsilon$-blowup of $\Gamma$ with respect to the L\'{e}vy-Prokhorov distance. Next, we show that
\begin{equation}
\inf_{\nu\in\Gamma^{\epsilon^{\prime}}}D_{\kl}(\nu,P)\geq\inf_{\nu\in\Gamma}D_{\kl,\lp}^{\epsilon}(\nu,P)\label{bound_inf's}%
\end{equation}
for any $0<\epsilon^{\prime}<\epsilon$. In fact, for any $\nu\in\Gamma^{\epsilon^{\prime}}$, it satisfies $\mathrm{LP}(\nu,\Gamma)\leq \epsilon^{\prime}<\epsilon$. Therefore, there exists $\nu^{\prime}\in\Gamma$ such that $\mathrm{LP}(\nu,\nu^{\prime})<\epsilon$. By the definition of $D_{\kl,\lp}^{\epsilon}(\nu^{\prime},P)$, we have%
\[
D_{\kl}(\nu,P)\geq D_{\kl,\lp}^{\epsilon}(\nu^{\prime},P)\geq\inf_{\nu^{\prime}\in\Gamma}D_{\kl,\lp}^{\epsilon}(\nu^{\prime},P).
\]
Taking the infimum over $\nu\in\Gamma^{\epsilon^{\prime}}$, we get (\ref{bound_inf's}). Now we fix $\epsilon>0$. By (\ref{dembo_finite_bound}) and (\ref{bound_inf's}), we have%
\begin{align*}
\Prob{P_{N}\in\Gamma} & \leq\left(  \frac{4}{\epsilon/2}\right) ^{C(\epsilon/2,\Sigma,||\cdot||)}\exp\left(  -N\inf_{\nu\in\Gamma^{\epsilon/2}}D_{\kl}(\nu,P)\right)  \\
& \leq\left(  \frac{8}{\epsilon}\right)  ^{C(\epsilon/2,\Sigma,||\cdot||)}\exp\left(  -N\inf_{\nu\in\Gamma}D_{\kl,\lp}^{\epsilon}(\nu,P)\right)  .
\end{align*}
Taking the infimum over $\epsilon>0$, we obtain%
\begin{equation}
\Prob{P_{N}\in\Gamma}\leq\inf_{\epsilon>0}\left\{  \left(  \frac{8}{\epsilon}\right)  ^{C(\epsilon/2,\Sigma,||\cdot||)}\exp\left(  -N\inf_{\nu\in\Gamma}D_{\kl,\lp}^{\epsilon}(\nu,P)\right)  \right\}  .\label{finite_LP_bound}%
\end{equation}
Then we consider the Wasserstein smoothing distance $D_{\kl,W}^{\epsilon}$. According to \cite[Theorem 2]{gibbs2002choosing}, we have%
\[
W(P,Q)\leq(\mathrm{diam}(\Sigma)+1)\mathrm{LP}(P,Q).
\]
Therefore, we must have $D_{\kl,\lp}^{\epsilon}(\nu,P)>D_{\kl,W}^{(\mathrm{diam}(\Sigma)+1)\epsilon}(\nu,P)$. Combining it with (\ref{finite_LP_bound}), we have%
\begin{align}
  & \Prob{P_{N}\in\Gamma}\nonumber\\
\leq & \inf_{\epsilon>0}\left\{  \left(  \frac{8}{\epsilon}\right)  ^{C(\epsilon/2,\Sigma,||\cdot||)}\exp\left(  -N\inf_{\nu\in\Gamma}D_{\kl,W}^{(\mathrm{diam}(\Sigma)+1)\epsilon}(\nu,P)\right)  \right\}\nonumber\\
=& \inf_{\epsilon>0}\left\{  \left(  \frac{8(\mathrm{diam}(\Sigma)+1)}{\epsilon}\right)  ^{C(\epsilon/(2(\mathrm{diam}(\Sigma)+1)),\Sigma,||\cdot||)}\exp\left(  -N\inf_{\nu\in\Gamma}D_{\kl,W}^{\epsilon}(\nu,P)\right)  \right\}  .\label{finite_W_bound}%
\end{align}
The last asymptotic upper bound can be obtained from (\ref{finite_LP_bound}) and (\ref{finite_W_bound}).
\end{proof}

\subsection{Proofs in Section \ref{sec:finite-alphabets}}

\subsubsection{Proof of Theorem \ref{KL_uni_feasibility}}

\begin{proof}
Recall that $\hat{c}_{\kl,r}(x,\hat{P})=\sup\{c(x,P):D_{\kl}(\hat{P},P)\leq r,P\in\mathcal{P}(\Sigma)\}$. Then we have the following inclusion relation:%
\[
\{\exists x\in X\text{ such that }c(x,P)>\hat{c}_{\kl,r}(x,P_{N})\}\subset\{D_{\kl}(P_{N},P)>r\}.
\]
Therefore, for any true distribution $P$, by the large deviations principle in Theorem \ref{finite_LDP}, we have%
\begin{align*}
&  \limsup_{N\rightarrow\infty}\frac{1}{N}\log \Prob{\exists x\in X\text{ such that }c(x,P)>\hat{c}_{\kl,r}(x,P_{N})}\\
&  \leq\limsup_{N\rightarrow\infty}\frac{1}{N}\log \Prob{D_{\kl}(P_{N},P)>r}\\
&  \leq-\inf_{\nu\in\{\nu:D_{\kl}(\nu,P)>r\}}D_{\kl}(\nu,P)\\
&  \leq-r,
\end{align*}
which means $\hat{c}_{\kl,r}(x,\hat{P})$ is uniformly feasible.
\end{proof}

\subsubsection{Proof of Theorem \ref{prediction_optimality_finite}}
\begin{proof}
First, suppose $\hat c(x,\hat P)$ is pointwise feasible. We can generalize the proof of \cite[Theorem 4]{vanparys2021data} in this case. In fact, notice that the set of disappointing estimator realizations $\mathcal{D}(x,P_2)=\{P':c(x,P_2)>\hat c(x,P')\}$ in their proof is also open when $\hat c$ is upper semicontinuous for any $x$. Then their proof still holds in this case. Next, suppose $\hat c(x,\hat P)$ is uniformly feasible, which implies that it is also pointwise feasible and thus reduces it to the first case.
\end{proof}

\subsubsection{Proof of Theorem \ref{optimality_KL_ball}}
\begin{proof}
If the conclusion does not hold, there exists $\hat{P}_{0},P_{0}\in\mathcal{P}(\Sigma)$ such that $D_{\kl}(\hat{P}_{0},P_{0})\leq r$ but $P_{0}\notin\mathcal{A}(\hat{P}_{0})$. We will show this actually contradicts to that $\hat{c}$ is reasonably pointwise feasible. Notice that $\{P_{0}\}$ and $\mathcal{A}(\hat{P}_{0})$ are two disjoint compact convex set on $\mathbb{R}^{|\Sigma|}$. So by the separating hyperplane theorem, there exists $a_{0}\in\mathbb{R}^{|\Sigma|}$ such that%
\[
\langle a_{0},P_{0}\rangle>\sup_{P\in\mathcal{A}(\hat{P}_{0})}\langle a_{0},P\rangle.
\]
By the continuity of inner product, we can choose $P_{1}=(1-\lambda)P_{0}+\lambda\hat{P}_{0}$ with a small $\lambda>0$ so that%
\begin{equation}
\langle a_{0},P_{1}\rangle>\sup_{P\in\mathcal{A}(\hat{P}_{0})}\langle a_{0},P\rangle.\label{separation_finite}%
\end{equation}
Suppose the support set $\Sigma$ is enumerated as $\Sigma=\{\bar\xi_1,\ldots,\bar\xi_{|\Sigma|}\}$. We define the decision-free loss function $\ell(x,\bar\xi_i)=a_{0i}$ for any $x\in X$ and any $\bar\xi_i\in\Sigma$. Under this loss function and the true probability $P_{1}$, by the large deviations principle, we have for any $x$ that%
\begin{align*}
&  \liminf_{N\rightarrow\infty}\frac{1}{N}\log \Prob{c(x,P_{1})>\hat{c}(x,P_{N})}\\
&  \geq-\inf_{\nu\in\{\nu:c(x,P_{1})>\hat{c}(x,\nu)\}^{\circ}}D_{\kl}(\nu,P_{1})\\
&  =-\inf_{\nu\in\{\nu:c(x,P_{1})>\hat{c}(x,\nu)\}}D_{\kl}(\nu,P_{1}),
\end{align*}
where the equality follows from that $\{\nu:c(x,P_{1})>\hat{c}(x,\nu)\}$ is open due to the upper semicontinuity of $\hat{c}(x,\nu)$ in $\nu$.
Notice that (\ref{separation_finite}) can be equivalently written as%
\begin{align*}
  & c(x,P_{1})=E_{P_{1}}[\ell(x,\xi)]=\langle a_{0},P_{1}\rangle>\\
  & \hspace{10em}\sup_{P\in\mathcal{A}(\hat{P}_{0})}\langle a_{0},P\rangle=\sup\{c(x,P):P\in\mathcal{A}(\hat{P}_{0})\}=\hat{c}(x,\hat{P}_{0}),
\end{align*}
which implies that $\hat{P}_{0}\in\{\nu:c(x,P_{1})>\hat{c}(x,\nu)\}$. Therefore,
\begin{align*}
  & \inf_{\nu:c(x,P_{1})>\hat{c}(x,\nu)}D_{\kl}(\nu,P_{1})\leq D_{\kl}(\hat{P}_{0},P_{1})\leq\\
  & \hspace{10em}(1-\lambda)D_{\kl}(\hat{P}_{0},P_{0})+\lambda D_{\kl}(\hat{P}_{0},\hat{P}_{0})\leq(1-\lambda)r<r,
\end{align*}
which implies the large deviations lower bound satisfies%
\[
\liminf_{N\rightarrow\infty}\frac{1}{N}\log \Prob{c(x,P_{1})>\hat{c}(x,P_{N})}>-r.
\]
So $\hat{c}$ is not reasonably pointwise feasible at rate $r$. This contradicts the choice of $\hat{c}$. Therefore, we must have $\{P:D_{\kl}(\hat{P},P)\leq r\}\subseteq\mathcal{A}(\hat{P})$ for any $\hat{P}$.
\end{proof}

\subsection{Proofs in Section \ref{sec:continuous-alphabets}}

\subsubsection{Proof of Theorem \ref{counterexample_any_r}}
\begin{proof}
Fix $r>0$. Let $x_{r}$ be a fixed constant such that $x_{r}>e^{r}-1$. We
construct the following continuous periodic function of period $T_{r}=\pi
+\pi/x_{r}$:
\[
g_{r}(x)=\left\{
\begin{array}
[c]{l}%
\sin(x),\\
x_{r}\sin\left(  x_{r}\left(  x+\frac{\pi}{x_{r}}-\pi\right)  \right)  \\
g_{r}(x-T_{r})\\
g_{r}(x+T_{r})
\end{array}
\left.
\begin{array}
[c]{l}%
\text{if }0\leq x\leq\pi\\
\text{if }\pi<x\leq\pi+\frac{\pi}{x_{r}}\\
\text{if }x>T_{r}\\
\text{if }x<0
\end{array}
\right.  \right.  .
\]
It's easy to verify that $\int_{0}^{T_{r}}g_{r}(x)dx=0$, $\max_{x\in
\lbrack0,T_{r}]}g_{r}(x)=1$ and $\min_{x\in\lbrack0,T_{r}]}g_{r}(x)=-x_{r}$.
We define the loss function $\ell_{r}(x,\xi)$ as
\[
\ell_{r}(x,\xi)=\left\{
\begin{array}
[c]{c}%
x^{2}g_{r}\left(  \frac{\xi}{x}\right)  \\
0
\end{array}
\left.
\begin{array}
[c]{c}%
\text{if }x\neq0\\
\text{if }x=0
\end{array}
\right.  \right.  .
\]
We can verify $\ell_{r}$ is continuous on $X\times\Sigma$: for $(x_{0},\xi
_{0})\in X\times\Sigma$, if $x_{0}\neq0$, then clearly $\ell_{r}$ is
continuous at $(x_{0},\xi_{0})$; if $x_{0}=0$ and $(x,\xi)\rightarrow
(x_{0},\xi_{0})$, then $|\ell_{r}(x,\xi)|\leq x^{2}\max\{1,x_{r}%
\}\rightarrow0=\ell_{r}(x_{0},\xi_{0})$.

Before we prove the theorem, we introduce the concept of ``uniformly distributed modulo 1'' \citep[Chapter 1 Section
6]{kuipers2012uniform}. For a real number $y$, we write $[y]$ as the integral part of $y$, i.e.,
the largest integer no more than $y$. Then $\mathrm{frac}(y):=y-[y]\in
\lbrack0,1)$ is the fractional part of $y$. For a vector $y\in\mathbb{R}^{s}$,
we define $\mathrm{frac}(y)\in\mathbb{R}^{s}$ as the vector of fractional part
of each component. For a sequence of vectors $\{y_{k},k\geq1\}\subseteq
\mathbb{R}^{s}$, we say it is uniformly distributed modulo 1 in $\mathbb{R}%
^{s}$ if
\[
\lim_{n\rightarrow\infty}\frac{\#\{1\leq k\leq n:a\leq\mathrm{frac}%
(y_{k})<b\}}{n}=\prod_{j=1}^{s}(b_{j}-a_{j})
\]
for all hyperrectangles $[a,b)\subset\lbrack0,1)^{s}$, where $a,b$ are vectors
and both the inequality $a\leq\mathrm{frac}(y_{k})<b$ and hyperrectangle
$[a,b)$ are in the componentwise sense. One consequence of $\{y_{k},k\geq1\}$
being uniformly distributed modulo 1 is that the fractional part
$\{\mathrm{frac}(y_{k}),k\geq1\}$ is dense on $[0,1)^{s}$. To see this, take
any point $z\in\lbrack0,1)^{s}$ and any $\epsilon\in(0,1/2)$. Consider the
hyperrectangle $[\max(z-\epsilon,0),\min(z+\epsilon,1))$, where all the
operators are interpreted in the componentwise sense. We have%
\begin{align*}
&  \lim_{n\rightarrow\infty}\frac{\#\{1\leq k\leq n:\max(z-\epsilon
,0)\leq\mathrm{frac}(y_{k})<\min(z+\epsilon,1)\}}{n}\\
&  =\prod_{j=1}^{s}(\min(z_{j}+\epsilon,1)-\max(z_{j}-\epsilon,0))\geq
\epsilon^{s}.
\end{align*}
Therefore, for any $\epsilon\in(0,1/2)$, there are infinitely many
$\mathrm{frac}(y_{k})$'s that lie in $[\max(z-\epsilon,0),\min(z+\epsilon,1))$
which proves the density. By \cite[Example 6.1]{kuipers2012uniform}, we know
that if $1,\theta_{1},\ldots,\theta_{s}$ are linearly independent over the
rationals, then the sequence $$\{(k\theta_{1},\ldots,k\theta_{s}),k\geq1\}$$ is
uniformly distributed modulo 1 and thus $\{(\mathrm{frac}(k\theta_{1}%
),\ldots,\mathrm{frac}(k\theta_{s})),k\geq1\}$ is dense on $[0,1)^{s}$.

Next, we start to prove this theorem. Recall the true distribution $P$ is
uniformly distributed on $[-1,1]$. Consider samples $\xi_{1},\ldots,\xi_{N}$.
We have%
\begin{align*}
&  \Prob{  1,\xi_{1},\ldots,\xi_{N}\text{ are linearly
dependent over the rationals}}  \\
&  =\Prob{  \exists q_{0},\ldots,q_{N}\in\mathbb{Q}\text{ such
that }q_{0}+q_{1}\xi_{1}+\cdots+q_{N}\xi_{N}=0} \\
&  \leq\cup_{q_{0},\ldots,q_{N}\in\mathbb{Q}}\Prob{
q_{0}+q_{1}\xi_{1}+\cdots+q_{N}\xi_{N}=0}  \\
&  =0,
\end{align*}
which means $1,\xi_{1},\ldots,\xi_{N}$ are actually linearly independent over
the rationals almost surely. We fix a sample such that $1,\xi_{1},\ldots
,\xi_{N}$ are linearly independent over the rationals. Then we know that
$\{(\mathrm{frac}(k\xi_{1}),\ldots,\mathrm{frac}(k\xi_{N})),k\geq1\}$ is dense
on $[0,1)^{N}$. This further implies that
\begin{align}
(g_{r}(T_{r}k\xi_{1}),\ldots,g_{r}(T_{r}k\xi_{N})) &  =(g_{r}(T_{r}k\xi
_{1}-T_{r}[k\xi_{1}]),\ldots,g_{r}(T_{r}k\xi_{N}-T_{r}[k\xi_{N}]))\nonumber\\
&  =(g_{r}(T_{r}\mathrm{frac}(k\xi_{1})),\ldots,g_{r}(T_{r}\mathrm{frac}%
(k\xi_{N})))\label{density_sequence}%
\end{align}
is dense on $[-x_{r},1]^{N}$. Moreover, by the Generalized Riemann-Lebesgue
Lemma (Lemma \ref{RL_lemma}), we have (note that both $g_{r}(T_{r}\xi)$ and
$g_{r}(-T_{r}\xi)$\ has period $1$)%
\begin{align}
& \lim_{k\rightarrow\infty}\int_{-1}^{1}\frac{1}{2}g_{r}(kT_{r}\xi
)d\xi\nonumber\\
& =\lim_{k\rightarrow\infty}\left(  \frac{1}{2}\int_{-1}^{0}g_{r}(kT_{r}%
\xi)d\xi+\frac{1}{2}\int_{0}^{1}g_{r}(kT_{r}\xi)d\xi\right)  \nonumber\\
& =\lim_{k\rightarrow\infty}\left(  \frac{1}{2}\int_{0}^{1}g_{r}(-kT_{r}%
\xi)d\xi+\frac{1}{2}\int_{0}^{1}g_{r}(kT_{r}\xi)d\xi\right)  \nonumber\\
& =\frac{1}{2}\lim_{k\rightarrow\infty}\int_{0}^{1}g_{r}(-kT_{r}\xi)d\xi
+\frac{1}{2}\lim_{k\rightarrow\infty}\int_{0}^{1}g_{r}(kT_{r}\xi
)d\xi\nonumber\\
& =\frac{1}{2}\int_{0}^{1}g_{r}(-T_{r}\xi)d\xi\int_{0}^{1}1d\xi+\frac{1}%
{2}\int_{0}^{1}g_{r}(T_{r}\xi)d\xi\int_{0}^{1}1d\xi\nonumber\\
& =0.\label{RL_application}%
\end{align}
Therefore, $\exists$ a large enough $k_{0}\in\mathbb{N}$ (depending on
$\xi_{1},\ldots,\xi_{N}$ and $r$) such that%
\[
-x_{r}\leq g_{r}(T_{r}k_{0}\xi_{i})<-\frac{e^{r}-1+x_{r}}{2},\forall
i=1,\ldots,N
\]
according to the density of (\ref{density_sequence}) and%
\begin{equation}
\frac{1}{2}\int_{-1}^{1}g_{r}(k_{0}T_{r}\xi)d\xi>1-\frac{e^{r}+1+x_{r}}%
{2}e^{-r}\label{lower_bound_integral}%
\end{equation}
according to (\ref{RL_application}). Notice that the KL-DRO predictor $\hat
{c}_{\kl,r}(x,P_{N})$ can be equivalently written as%
\begin{align}
\sup_{p_{i}}\text{ } &  \sum_{i=1}^{N}x^{2}g_{r}\left(  \frac{\xi_{i}}%
{x}\right)  p_{i}+p_{N+1}\sup_{\xi}\left[  x^{2}g_{r}\left(  \frac{\xi}%
{x}\right)  \right]  \nonumber\\
\st &  \sum_{i=1}^{N+1}p_{i}=1\label{equi_form_KL-DRO}\\
&  p_{i}\geq0,i=1,\ldots,N+1\nonumber\\
&  \sum_{i=1}^{N}\frac{1}{N}\log\left(  \frac{1}{Np_{i}}\right)  \leq
r\nonumber
\end{align}
We can see any feasible solution in (\ref{equi_form_KL-DRO}) must satisfy
$p_{i}>0$ for $i=1,\ldots,N$. Additionally, the last constraint can be written
as%
\[
r\geq-\log N-\log\sqrt[N]{p_{1}\cdots p_{N}}\Leftrightarrow\sqrt[N]%
{p_{1}\cdots p_{N}}\geq\frac{1}{N}e^{-r},
\]
which implies%
\[
\sum_{i=1}^{N}p_{i}\geq N\sqrt[N]{p_{1}\cdots p_{N}}\geq e^{-r}.
\]
Taking $x=1/(T_{r}k_{0})$ and for any feasible solution in
(\ref{equi_form_KL-DRO}), we have%
\begin{align*}
&  \sum_{i=1}^{N}\frac{1}{(T_{r}k_{0})^{2}}g_{r}(T_{r}k_{0}\xi_{i}%
)p_{i}+p_{N+1}\sup_{\xi}\left[  \frac{1}{(T_{r}k_{0})^{2}}g_{r}(T_{r}k_{0}%
\xi)\right]  \\
&  =\frac{1}{(T_{r}k_{0})^{2}}\left(  \sum_{i=1}^{N}g_{r}(T_{r}k_{0}\xi
_{i})p_{i}+1-\sum_{i=1}^{N}p_{i}\right)  \\
&  <\frac{1}{(T_{r}k_{0})^{2}}\left(  -\frac{e^{r}-1+x_{r}}{2}\sum_{i=1}%
^{N}p_{i}+1-\sum_{i=1}^{N}p_{i}\right)  \\
&  =\frac{1}{(T_{r}k_{0})^{2}}\left(  1-\frac{e^{r}+1+x_{r}}{2}\sum_{i=1}%
^{N}p_{i}\right)  \\
&  \leq\frac{1}{(T_{r}k_{0})^{2}}\left(  1-\frac{e^{r}+1+x_{r}}{2}%
e^{-r}\right)  ,
\end{align*}
This implies that for this fixed sample $\xi_{1},\ldots,\xi_{N}$ and
$x=1/(T_{r}k_{0})$,
\begin{align*}
& \hat{c}_{\kl,r}(x,P_{N})\leq\frac{1}{(T_{r}k_{0})^{2}}\left(  1-\frac
  {e^{r}+1+x_{r}}{2}e^{-r}\right)  <\\
  & \hspace{8em}\frac{1}{(T_{r}k_{0})^{2}}\int_{-1}^{1}%
\frac{1}{2}g_{r}(T_{r}k_{0}\xi)d\xi=\mathbb{E}_{P}\left[  x^{2}g_{r}\left(
\frac{\xi}{x}\right)  \right]  \equiv\mathbb{E}_{P}[\ell_{r}(x,\xi)],
\end{align*}
where the strict inequality is from (\ref{lower_bound_integral}). Hence, we
obtain%
\begin{align*}
&  \Prob{  \exists x\in X\text{ such that }\mathbb{E}_{P}%
[\ell_{r}(x,\xi)]>\hat{c}_{\kl,r}(x,P_{N})}  \\
\geq &  \Prob{  1,\xi_{1},\ldots,\xi_{N}\text{ are linearly
independent over the rationals}}  \\
= &  1.
\end{align*}
This proves the first part of the theorem.

Now we prove the second part of the theorem. Suppose $r<\log2$. Then we can
choose $x_{r}\equiv1>e^{r}-1$. In this case, $g_{r}(x)\equiv\sin(x)$, which is
a periodic odd function. Therefore, we can actually see that
\[
c(x,P)=\mathbb{E}_{P}\left[  x^{2}g_{r}\left(  \frac{\xi}{x}\right)  \right]
=0,\forall x\in X.
\]
The first part of the theorem ensures that $\{\exists x\in X$ such that
$\hat{c}_{\kl,r}(x,P_{N})<\mathbb{E}_{P}[\ell_{r}(x,\xi)]=0\}$ happens with
probability 1. On this event, we must have%
\[
\hat{c}_{\kl,r}(\hat{x}_{r}(P_{N}),P_{N})=\min_{x\in X}\hat{c}_{\kl,r}%
(x,P_{N})<0=c(\hat{x}_{r}(P_{N}),P)=\min_{x\in X}c(x,P),
\]
which implies that
\[
\Prob{  c(\hat{x}_{r}(P_{N}),P)>\hat{c}_{\kl,r}(\hat{x}%
_{r}(P_{N}),P_{N})}  =1
\]
and%
\[
\Prob{  \min_{x\in X}c(x,P)>\min_{x\in X}\hat{c}_{\kl,r}%
(x,P_{N})} =1.
\]
\end{proof}

\subsubsection{Proof of Theorem \ref{smooth_KL_uni_feasibility}}
\begin{proof}
We prove this theorem for $\hat{c}_{\kl,\lp,r}^{\epsilon}(x,\hat{P})$ and the same argument can be applied to $\hat{c}_{\kl,W,r}^{\epsilon}(x,\hat{P})$. Recall again that $\hat{c}_{\kl,\lp,r}^{\epsilon}(x,\hat{P})=\sup\{c(x,P):P\in\mathcal{P}(\Sigma),D_{\kl,\lp}^{\epsilon}(\hat{P},P)\leq r\}$. Therefore, we have%
\[
\{\exists x\in X\text{ s.t. }c(x,P)>\hat{c}_{\kl,\lp,r}^{\epsilon}(x,P_{N})\}\subseteq\{D_{\kl,\lp}^{\epsilon}(P_{N},P)>r\}.
\]
In particular, the measurability of $\{\exists x\in X\text{ s.t. }c(x,P)>\hat{c}_{\kl,\lp,r}^{\epsilon}(x,P_{N})\}$ is guaranteed by Lemma \ref{usc_smooth_f_DRO} and Lemma \ref{event_measurability}. By Lemma \ref{property_smooth_f_div}, $D_{\kl,\lp}^{\epsilon}$ is lower semicontinuous, which means that the set $\Gamma=\{\nu\in\mathcal{P}(\Sigma):D_{\kl,\lp}^{\epsilon}(\nu,P)>r\}$ is open. So we have%
\begin{align*}
& \limsup_{N\rightarrow\infty}\frac{1}{N}\log \Prob{\exists x\in X\text{ s.t. }c(x,P)>\hat{c}_{\kl,\lp,r}^{\epsilon}(x,P_{N})}\\
& \leq\limsup_{N\rightarrow\infty}\frac{1}{N}\log \Prob{P_{N}\in\Gamma}\\
& \leq-\inf_{\nu\in\Gamma}D_{\kl,\lp}^{\epsilon}(\nu,P)\\
& \leq-r,
\end{align*}
where the second inequality follows from Theorem \ref{thm:sanov-smooth} and the last inequality is due to the definition of $\Gamma$. This proves the uniform feasibility.
\end{proof}

\subsubsection{Proof of Theorem \ref{monotonicity_DRO}}
\begin{proof}
We prove the result for LP smoothing and similar argument can be applied to Wasserstein smoothing. Recall that%
\[
D_{f,\lp}^{\epsilon}(P^{\prime},P)=\inf\{D_{f}(Q,P):Q\in\mathcal{P}(\Sigma),\mathrm{LP}(P^{\prime},Q)\leq\epsilon\}.
\]
Therefore, for $0<\epsilon_{1}<\epsilon_{2}$, we have%
\begin{equation}
D_{f,\lp}^{\epsilon_{2}}(P^{\prime},P)\leq D_{f,\lp}^{\epsilon_{1}}(P^{\prime},P)\leq D_{f}(P^{\prime},P),\label{robust_monotone}%
\end{equation}
which implies the monotonicity of $\{P\in\mathcal{P}(\Sigma):D_{f,\lp}^{\epsilon}(P^{\prime},P)\leq r\}$:
\begin{align}
  &\{P\in\mathcal{P}(\Sigma):D_{f}(P^{\prime},P)\leq r\}\subseteq\{P\in\mathcal{P}(\Sigma):D_{f,\lp}^{\epsilon_{1}}(P^{\prime},P)\leq r\}\nonumber\\
  &\hspace{17em}\subseteq\{P\in\mathcal{P}(\Sigma):D_{f,\lp}^{\epsilon_{2}}(P^{\prime},P)\leq r\}.
\label{robust_inclusion}%
\end{align}
We also obtain%
\[
\{P\in\mathcal{P}(\Sigma):D_{f}(P^{\prime},P)\leq r\}\subset\bigcap_{\epsilon>0}\{P\in\mathcal{P}(\Sigma):D_{f,\lp}^{\epsilon}(P^{\prime},P)\leq r\}.
\]
Now we prove the other direction. Take $P\in\bigcap_{\epsilon>0}\{P\in\mathcal{P}(\Sigma):D_{f,\lp}^{\epsilon}(P^{\prime},P)\leq r\}$, i.e., $D_{f,\lp}^{\epsilon}(P^{\prime},P)\leq r$ for any $\epsilon>0$. Letting $\epsilon\downarrow0$ and using Lemma \ref{property_smooth_f_div}, we have $D_{f}(P^{\prime},P)\leq r$. This proves the other direction and thus we obtain%
\[
\{P\in\mathcal{P}(\Sigma):D_{f}(P^{\prime},P)\leq r\}=\bigcap_{\epsilon>0}\{P\in\mathcal{P}(\Sigma):D_{f,\lp}^{\epsilon}(P^{\prime},P)\leq r\}.
\]

Next we consider the DRO predictors. From Equation (\ref{robust_inclusion}), we know that $\hat{c}_{f,\lp,r}^{\epsilon}(x,\hat{P})\geq\hat{c}_{f,r}(x,\hat{P})$ and $\hat{c}_{f,\lp,r}^{\epsilon}(x,\hat{P})$ is non-increasing when $\epsilon$ decreases. So the limit of $\hat{c}_{f,\lp,r}^{\epsilon}(x,\hat{P})$ for $\epsilon\downarrow 0$ exists and satisfies%
\[
\lim_{\epsilon\downarrow0}\hat{c}_{f,\lp,r}^{\epsilon}(x,\hat{P})\geq\hat{c}_{f,r}(x,\hat{P}).
\]
To prove the other direction, we take $\epsilon_{n}\downarrow0$. By the attainability of $\hat{c}_{f,\lp,r}^{\epsilon}(x,\hat{P})$ in Lemma \ref{attainability_smooth_f_DRO}, there exists $Q_{n}$ such that $\hat{c}_{f,\lp,r}^{\epsilon_{n}}(x,\hat{P})=c(x,Q_{n})$ and $D_{f,\lp}^{\epsilon_{n}}(\hat{P},Q_{n})\leq r$. By the compactness of $\mathcal{P}(\Sigma)$, we extract a converging subsequence $Q_{n_{k}}\overset{w}{\rightarrow}Q$. Fix $k_{0}\in\mathbb{N}$. When $k\geq k_{0}$, by the monotonicity in (\ref{robust_monotone}), we have%
\[
r\geq\liminf_{k\rightarrow\infty}D_{f,\lp}^{\epsilon_{n_{k}}}(\hat{P},Q_{n_{k}})\geq\liminf_{k\rightarrow\infty}D_{f,\lp}^{\epsilon_{n_{k_{0}}}}(\hat{P},Q_{n_{k}})\geq D_{f,\lp}^{\epsilon_{n_{k_{0}}}}(\hat{P},Q),
\]
where the last inequality uses the lower semicontinuity of $D_{f,\lp}^{\epsilon_{n_{k_{0}}}}$. Then letting $k_{0}\rightarrow\infty$ and using Lemma \ref{property_smooth_f_div}, we have
\(
D_{f}(\hat{P},Q)\leq r.
\)
Therefore, we have%
\[
\hat{c}_{f,r}(x,\hat{P})\geq c(x,Q)\geq\limsup_{k\rightarrow\infty}c(x,Q_{n_{k}})=\limsup_{k\rightarrow\infty}\hat{c}_{f,\lp,r}^{\epsilon_{n_{k}}}(x,\hat{P})=\lim_{\epsilon\downarrow0}\hat{c}_{f,\lp,r}^{\epsilon}(x,\hat{P}),
\]
where the first inequality follows from the definition of $\hat{c}_{f,r}(x,\hat{P})$, the second inequality follows from the upper semicontinuity of $c(x,Q)$ in $Q$ (due to the upper semicontinuity of $\ell$, see the proof of Lemma \ref{attainability_smooth_f_DRO}), the first equality follows from $\hat{c}_{f,\lp,r}^{\epsilon_{n}}(x,\hat{P})=c(x,Q_{n})$ and the last equality follows from the monotonicity of $\hat{c}_{f,\lp,r}^{\epsilon}(x,\hat{P})$ in $\epsilon$. Combining both directions, we obtain%
\(
\lim_{\epsilon\downarrow0}\hat{c}_{f,\lp,r}^{\epsilon}(x,\hat{P})=\hat{c}_{f,r}(x,\hat{P}).
\)
\end{proof}

\subsubsection{Proof of Theorem \ref{prediction_optimality_infinite}}
\begin{proof}
According to Theorem \ref{monotonicity_DRO}, to show the asymptotic efficiency, it suffices to show $\hat{c}_{\kl,r}(x,\hat{P})\leq\hat{c}(x,\hat{P})$. Moreover, since uniform feasibility implies pointwise feasibility, it suffices to prove $\hat{c}_{\kl,r}(x,\hat{P})\leq\hat{c}(x,\hat{P})$ when $\hat{c}(x,\hat{P})$ satisfies the pointwise feasibility (\ref{eq:feasibility}), which is already proved in \cite[Theorem 10]{vanparys2021data} (notice that the disappointing estimator realization $D(x,P_{2})=\{P^{\prime}:c(x,P_{2})>\hat{c}(x,P^{\prime})\}$ is still open because of the upper semicontinuity of $\hat{c}$ so their proof still holds here). 

Then we show the equivalence of the two claims. Suppose (1) holds. Then we have $\hat{c}_{\kl,r}(x,\hat{P})\leq\hat{c}(x,\hat{P})$ from the above result. On the other hand, $\hat{c}_{\kl,W,r}^{\epsilon}(x,\hat{P})\in\mathcal{C}_{u}(X\times\mathcal{P}(\Sigma))$ is feasible for both (\ref{eq:feasibility}) and (\ref{eq:feasibility-uniform}). In either case, the efficiency of $\hat{c}(x,\hat{P})$ implies $\hat{c}(x,\hat{P})\leq\hat{c}_{\kl,W,r}^{\epsilon}(x,\hat{P})$. Letting $\epsilon\downarrow0$, we have $\hat{c}(x,\hat{P})\leq\hat{c}_{\kl,r}(x,\hat{P})$. Combining both directions, we obtain $\hat{c}(x,\hat{P})\equiv\hat{c}_{\kl,r}(x,\hat{P})$, which proves (2). Now suppose (2) holds. Then (1) is obvious. In particular, the fact $\hat{c}_{\kl,r}(x,\hat{P})\in\mathcal{C}_u(X\times\mathcal{P}(\Sigma))$ is due to Lemma \ref{usc_smooth_f_DRO} by setting $\epsilon=0$.
\end{proof}

\subsubsection{Proof of Theorem \ref{optimality_KL_ball_continuous}}
\begin{proof}
According to Theorem \ref{monotonicity_DRO}, it suffices to show $\{P\in\mathcal{P}(\Sigma):D_{\kl}(\hat{P},P)\leq r\}\subseteq\mathcal{A}(\hat{P})$. Suppose the conclusion is not true. Then there exist $\hat{P}_{0},P_{0}\in\mathcal{P}(\Sigma)$ such that $D_{\kl}(\hat{P}_{0},P_{0})\leq r$ but $P_{0}\notin\mathcal{A}(\hat{P}_{0})$. Consider the two disjoint convex sets $\{P_{0}\}$ and $\mathcal{A}(\hat{P}_{0})$. They are closed on $\mathcal{P}(\Sigma)$. Let $\mathcal{M}(\Sigma)$ denote all finite signed measures on $\Sigma$ and it is a locally convex Hausdorff topological vector space when equipped with open sets $U_{\phi,y\delta}$ (see Section \ref{sec:infinite_LDP}). Since $\mathcal{P}(\Sigma)$ itself is a closed and compact subset of $\mathcal{M}(\Sigma)$, both $\{P_{0}\}$ and $S(P_{0}^{\prime})$ are also closed and compact on $\mathcal{M}(\Sigma)$. By the Hahn-Banach Theorem (see, e.g., \citep{dembo2009large} Theorem B.6), there exists a functional on $\mathcal{M}(\Sigma)$ which can be identified as a continuous function $\ell_{0}(\xi):\Sigma\mapsto\mathbb{R}$ such that%
\begin{equation}
\int_{\Sigma}\ell_{0}(\xi)dP_{0}(\xi)>\sup_{P\in\mathcal{A}(\hat{P}_{0})}\int_{\Sigma}\ell_{0}(\xi)dP(\xi).\label{separation_continuous}%
\end{equation}
We define the decision-free loss function $\ell(x,\xi):=\ell_{0}(\xi)$ with expected loss $c(x,P)=\mathbb{E}_{P}[\ell(x,\xi)]=\mathbb{E}_{P}[\ell_{0}(\xi)]$. According to (\ref{separation_continuous}) and the linearity of $c(x,P)$ in $P$, we have%
\[
\lim_{\lambda\downarrow0}c(x,(1-\lambda)P_{0}+\lambda\hat{P}_{0})=c(x,P_{0})>\sup_{P\in\mathcal{A}(\hat{P}_{0})}c(x,P).
\]
Therefore, $\exists\lambda_{0}\in(0,1)$ small enough such that%
\[
c(x,(1-\lambda_{0})P_{0}+\lambda_{0}\hat{P}_{0})>\sup_{P\in\mathcal{A}(\hat{P}_{0})}c(x,P).
\]
Moreover, by the convexity of $D_{\kl}$, we have%
\[
D_{\kl}(\hat{P}_{0},(1-\lambda_{0})P_{0}+\lambda_{0}\hat{P}_{0})\leq(1-\lambda_{0})D_{\kl}(\hat{P}_{0},P_{0})+\lambda_{0}D_{\kl}(\hat{P}_{0},\hat{P}_{0})\leq(1-\lambda_{0})r.
\]
For simplicity, we define $P_{\lambda_{0}}=(1-\lambda_{0})P_{0}+\lambda_{0}\hat{P}_{0}$. Now we will show $\hat{c}(x,\hat{P})$ is not feasible under the loss function $\ell(x,\xi)$ and true distribution $P_{\lambda_{0}}$. Since the loss is decision-free, pointwise feasibility is equivalent to uniform feasibility. So we fix some $x\in X$ and consider the pointwise feasibility in the following. Let $\Gamma(x)=\{\nu\in\mathcal{P}(\Sigma):c(x,P_{\lambda_{0}})>\hat{c}(x,\nu)\}$. Since $\hat{c}(x,\nu)\in C_{u}(X\times\mathcal{P}(\Sigma))$ is upper semicontinuous in $v$, $\Gamma(x)$ is an open subset of $\mathcal{P}(\Sigma)$. By the large deviations lower bound in Theorem \ref{infinite_LDP}, we have%
\begin{align*}
  \liminf_{N\rightarrow\infty}\frac{1}{N}\log P_{\lambda_{0}}^{\infty}(c(x,P_{\lambda_{0}})>\hat{c}(x,P_{N}))=& \liminf_{N\rightarrow\infty}\frac{1}{N}\log P_{\lambda_{0}}^{\infty}(P_{N}\in\Gamma(x))\\
  \geq&-\inf_{\nu\in\Gamma(x)}D_{\kl}(\nu,P_{\lambda_{0}}).
\end{align*}
Note that $\hat{P}_{0}\in\Gamma(x)$ because%
\[
c(x,P_{\lambda_{0}})>\sup_{P\in\mathcal{A}(\hat{P}_{0})}c(x,P)\equiv\sup\{c(x,P):P\in\mathcal{A}(\hat{P}_{0})\}=\hat{c}(x,\hat{P}_{0}).
\]
Therefore,%
\begin{align*}
  \liminf_{N\rightarrow\infty}\frac{1}{N}\log P_{\lambda_{0}}^{\infty}(c(x,P_{\lambda_{0}})>\hat{c}(x,P_{N}))\geq&-\inf_{\nu\in\Gamma(x)}D_{\kl}(\nu,P_{\lambda_{0}})\geq-D_{\kl}(\hat{P}_{0},P_{\lambda_{0}})\\
  \geq& -(1-\lambda_{0})r,
\end{align*}
which means $\hat{c}$ is not reasonably feasible and contradicts our assumption. Therefore, the conclusion must be true.
\end{proof}

\subsubsection{Proof of Corollary \ref{feasibility_smooth_f_DRO}}
\begin{proof}
  By the definition of $R_{D,W}^{\epsilon}(r)$ and $R_{D,\lp}^{\epsilon}(r)$ we have that $\hat{c}_{\kl,W,r}^\epsilon(x,\hat{P}) \leq \hat{c}_{D, R_{D,W}^{\epsilon}(r)}(x,\hat{P})$ as well as $\hat{c}_{\kl,\lp,r}^\epsilon(x,\hat{P}) \leq \hat{c}_{D, R_{D,\lp}^{\epsilon}(r)}(x,\hat{P})$. Moreover, due to the convexity and lower semicontinuity of the distance $D$, both $\hat{c}_{D, R_{D,W}^{\epsilon}(r)}(x,\hat{P})$ and $\hat{c}_{D, R_{D,\lp}^{\epsilon}(r)}(x,\hat{P})$ belong to the class $\mathcal{C}_u(X\times\mathcal{P}(\Sigma))$ by similar arguments in Lemma \ref{attainability_smooth_f_DRO} and Lemma \ref{usc_smooth_f_DRO}. So the event $\{\exists x\in X~\st~\mathbb{E}_{P}[\ell(x,\xi)]>\hat{c}(x,P_{N})\}$ are measurable for these two predictors according to Lemma \ref{event_measurability}.
  The feasibility result follows now from Theorem \ref{smooth_KL_uni_feasibility}.
  
  Now suppose $0\leq\eta<R_{D}(r)$ and $r>0$. Then there exists $\hat{P}_{0}$ and $P_{0}$ such that $D_{\kl}(\hat{P}_{0},P_{0})\leq r$ and $D(\hat{P}_{0},P_{0})>\eta$, which implies that $\{P:D_{\kl}(\hat{P}_{0},P)\leq r\}\varsubsetneq\{P:D(\hat{P}_{0},P)\leq\eta\}$. Then according to Theorem \ref{optimality_KL_ball_continuous}, there exists a decision-free loss function such that $\hat{c}_{D,\eta}(x,\hat{P})$ is not pointwise feasible nor uniformly feasible at rate $r$.
\end{proof}

\subsubsection{Proof of Theorem \ref{quantification2}}
\begin{proof}
The quantification can be written as%
\begin{equation*}
\sup_{\ell\in \mc L_1}|\hat{c}_{1}(x,\hat{P})-\hat{c}_{2}(x,\hat{P})|=\max\left\{  \sup_{\ell\in \mc L_1}\hat{c}_{1}(x,\hat{P})-\hat{c}_{2}(x,\hat{P}),\sup_{\ell\in \mc L_1}\hat{c}_{2}(x,\hat{P})-\hat{c}_{1}(x,\hat{P})\right\}  .
\end{equation*}
Now we simplify the inner superma. Take the first one as an example. Since the quantification is measured for each fixed $x$, we can omit the dependence on $x$, i.e., the first supremum can be equivalently written as%
\begin{align*}
& \sup_{\ell\in \mc L_1}\hat{c}_{1}(x,\hat{P})-\hat{c}_{2}(x,\hat{P})\\
& =\sup_{\ell\in \mc L_1}\left(  \sup_{P\in\mathcal{A}_{1}(\hat{P})}\int_{\Sigma}\ell(x,\xi)dP-\sup_{Q\in\mathcal{A}_{2}(\hat{P})}\int_{\Sigma}\ell(x,\xi)dQ\right)  \\
& =\sup_{||\ell_{0}||_{L}\leq1}\left(  \sup_{P\in\mathcal{A}_{1}(\hat{P})}\int_{\Sigma}\ell_{0}(\xi)dP-\sup_{Q\in\mathcal{A}_{2}(\hat{P})}\int_{\Sigma}\ell_{0}(\xi)dQ\right)  ,
\end{align*}
where $\ell_{0}:\Sigma\mapsto\mathbb{R}$ and $||\ell_{0}||_{L}=\sup_{\xi_{1},\xi_{2}\in\Sigma}|\ell_{0}(\xi_{1})-\ell_{0}(\xi_{2})|/||\xi_{1}-\xi_{2}||$. We can further simplify it as%
\begin{align*}
& \sup_{\ell\in \mc L_1}\hat{c}_{1}(x,\hat{P})-\hat{c}_{2}(x,\hat{P})\\
& =\sup_{||\ell_{0}||_{L}\leq1}\sup_{P\in\mathcal{A}_{1}(\hat{P})}\left(\int_{\Sigma}\ell_{0}(\xi)dP-\sup_{Q\in\mathcal{A}_{2}(\hat{P})}\int_{\Sigma}\ell_{0}(\xi)dQ\right)  \\
& =\sup_{P\in\mathcal{A}_{1}(\hat{P})}\sup_{||\ell_{0}||_{L}\leq1}\left(\int_{\Sigma}\ell_{0}(\xi)dP-\sup_{Q\in\mathcal{A}_{2}(\hat{P})}\int_{\Sigma}\ell_{0}(\xi)dQ\right)  \\
& =\sup_{P\in\mathcal{A}_{1}(\hat{P})}\sup_{||\ell_{0}||_{L}\leq1}\inf_{Q\in\mathcal{A}_{2}(\hat{P})}\left(  \int_{\Sigma}\ell_{0}(\xi)dP-\int_{\Sigma}\ell_{0}(\xi)dQ\right)  .
\end{align*}
For any fixed $P$, we define $f(\ell_{0},Q)=\int_{\Sigma}\ell_{0}(\xi)dP-\int_{\Sigma}\ell_{0}(\xi)dQ$ as a functional in $(\ell_{0},Q)$. Notice that $f$ is bilinear in $(\ell_{0},Q)$ and thus is concave in $\ell_{0}$ and convex in $Q$. Besides, $\mathcal{A}_{2}(\hat{P})$ is a compact set, and $f(\ell_{0},Q)$ is continuous in $Q$ for any $\ell_{0}$ with $||\ell_{0}||_{L}\leq1$. Then by Sion's minimax theorem \cite[Theorem 4.2']{sion1958general}, we have%
\begin{align*}
  &\sup_{||\ell_{0}||_{L}\leq1}\inf_{Q\in\mathcal{A}_{2}(\hat{P})}\left(\int_{\Sigma}\ell_{0}(\xi)dP-\int_{\Sigma}\ell_{0}(\xi)dQ\right) \\ =&\inf_{Q\in\mathcal{A}_{2}(\hat{P})}\sup_{||\ell_{0}||_{L}\leq1}\left(\int_{\Sigma}\ell_{0}(\xi)dP-\int_{\Sigma}\ell_{0}(\xi)dQ\right)  \\
  =&\inf_{Q\in\mathcal{A}_{2}(\hat{P})}W(P,Q),
\end{align*}
where the second inequality follows from the Kantorovich-Rubinstein dual representation. Therefore, we have%
\[
\sup_{\ell\in \mc L_1}\hat{c}_{1}(x,\hat{P})-\hat{c}_{2}(x,\hat{P})=\sup_{P\in\mathcal{A}_{1}(\hat{P})}\inf_{Q\in\mathcal{A}_{2}(\hat{P})}W(P,Q).
\]
Similarly, we can also obtain%
\[
\sup_{\ell\in \mc L_1}\hat{c}_{2}(x,\hat{P})-\hat{c}_{1}(x,\hat
{P})=\sup_{Q\in\mathcal{A}_{2}(\hat{P})}\inf_{P\in\mathcal{A}_{1}(\hat{P}%
)}W(P,Q),
\]
which gives the desired quantification%
\[
\sup_{\ell\in \mc L_1}|\hat{c}_{1}(x,\hat{P})-\hat{c}_{2}(x,\hat{P})|\!=\!\max\left\{ \! \sup_{P\in\mathcal{A}_{1}(\hat{P})}\inf_{Q\in\mathcal{A}_{2}(\hat{P})}\!W(P,Q),\!\!\sup_{Q\in\mathcal{A}_{2}(\hat{P})}\inf_{P\in\mathcal{A}_{1}(\hat{P})}\!W(P,Q)\right\}.
\]

When $\hat{c}_{1}$ is a reasonably DRO predictor at rate $r>0$ and $\hat{c}_{2}=\hat{c}_{\kl,r}$, we have%
\[
\sup_{Q\in\mathcal{A}_{2}(\hat{P})}\inf_{P\in\mathcal{A}_{1}(\hat{P})}W(P,Q)=0
\]
because $\mathcal{A}_{2}(\hat{P})=\{P:D_{\kl}(\hat{P},P)\leq r\}\subseteq\mathcal{A}_{1}(\hat{P})$. So the quantification reduces to%
\begin{align*}
  \sup_{\ell\in \mc L_1}|\hat{c}_{1}(x,\hat{P})-\hat{c}_{2}(x,\hat{P})|=& \sup_{\ell\in \mc L_1}\hat{c}_{1}(x,\hat{P})-\hat{c}_{\kl,r}(x,\hat{P})\\
  =&\sup_{P\in\mathcal{A}_{1}(\hat{P})}\inf_{Q\in\{Q:D_{\kl}(\hat{P},Q)\leq r\}}W(P,Q).
\end{align*}
\end{proof}%

\subsection{Proofs in Section \ref{sec:empirical-predictors}}

\subsubsection{Proof of Theorem \ref{proposition:el-infeasibility}}
\begin{proof}
  Fix $x \in X$. If $\inf_{\xi\in\Sigma} \ell(x, \xi) < \sup_{\xi\in \Sigma}\ell(x, \xi)$ there exists $\xi_a\in \Sigma$ and $\xi_b\in \Sigma$ for which we have $\ell(x, \xi_a)<\ell(x, \xi_b)$. Let us denote with $\bar \Sigma = \{\xi_a, \xi_b\}$. Consider $P_c\in \mc P(\bar \Sigma)$ so that $P_c(\xi_a)=1-c$ and $P_c(\xi_b)=c$ with $\E{P_c}{\ell(x, \xi)}=(1-c)\ell(x, \xi_a)+c\ell(x, \xi_b)>\ell(x,\xi_a)$ for all $c\in (0, 1]$. Consider the event in which $\xi_1=\dots=\xi_N=\xi_a$ so that $P_N=\delta_{\xi_a}$. Clearly, if this event occurs then $\hat c(x, P_N) \leq \ell(x, \xi_a)<\E{P_c}{\ell(x, \xi)}$.
  Hence,
  \begin{align*}
    & \sup_{P\in \mc P(\Sigma)} \limsup_{N\to\infty}\frac{1}{N}\log \Prob{\E{P}{\ell(x, \xi)}>\hat c(x, P_N)}\\
    \geq  &\sup_{P\in \mc P(\bar \Sigma)}\limsup_{N\to\infty}\frac{1}{N}\log \Prob{\E{P}{\ell(x, \xi)}>\hat c(x, P_N)}\\
    \geq & \sup_{P\in \mc P(\bar \Sigma)}\limsup_{N\to\infty}\frac{1}{N}\log \Prob{P_N=\delta_{\xi_a}}\\
     \geq &  \sup_{c\in (0, 1]} \frac{1}{N}\log\left( (1-c)^N\right) =  \sup_{c\in (0, 1]} \log (1-c) = 0.
  \end{align*}
  The claim follows by observing that as $\Prob{\E{P}{\ell(x, \xi)}>\hat c(x, P_N)}\leq 1$ for all $P\in \mc P(\Sigma)$ and $N\geq 1$ we also have \(\sup_{P\in \mc P(\Sigma)} \limsup_{N\to\infty}\frac{1}{N}\log \Prob{\E{P}{\ell(x, \xi)}>\hat c(x, P_N)}\leq 0\).
\end{proof}

\subsubsection{Proof of Theorem \ref{lemma:hr:primal:reduction}}

\begin{proof}

  First, in \eqref{eq:minimax-equality} we will show that
  \begin{align*}
    & \hat c^{\epsilon}_{f, d, r}(x, P_N)\\
    = & \inf \Big\{\sup \set{\textstyle\int \lambda f^\star\left(\frac{\ell(x, \xi)-\eta}{\lambda}\right) + r\lambda +\eta\, \d Q(\xi) }{ Q\in \mc P(\Sigma),~\lp(P_N, Q)\leq \epsilon } : \\
                                             & \hspace{20em} \eta\in\Re, ~\lambda\in\Re_+, ~\ell^\infty(x) - \eta \leq f^\infty \lambda \Big\}.
  \end{align*}
  we remark that for any $\lambda\in \Re_+$ and $\eta\in \Re$ with $\ell^\infty(x) - \eta \leq f^\infty \lambda$ we have that
  \(
  a \mapsto \lambda f^\star\left(\tfrac{(a-\eta)}{\lambda}\right)
  \)
  is a nondecreasing function following Lemma \ref{lemma:conjugate-function}.
  In particular, denote with $\{\xi_1,\ldots,\xi_K\}$ the support of $P_N \in \mc P(\Sigma)$ ordered so that $\ell^\epsilon(x, \xi_1)\leq \dots \leq \ell^\epsilon(x, \xi_K)$ then clearly we also have $\lambda f^\star\left(\tfrac{(\ell^\epsilon(x, \xi_1)-\eta)}{\lambda}\right)\leq \dots\leq \lambda f^\star\left(\tfrac{(\ell^\epsilon(x, \xi_K)-\eta)}{\lambda}\right)$. Following now Theorem 2.2 in \citep{bennouna2022holistic} that for all $\lambda\in \Re_+$ and $\eta\in \Re$ with $\ell^\infty(x) - \eta \leq f^\infty \lambda$  we have
    \begin{align*}
      & \sup \set{\textstyle\int \lambda f^\star\left(\frac{\ell(x, \xi)-\eta}{\lambda}\right) + r\lambda +\eta\, \d Q(\xi) }{ Q\in \mc P(\Sigma),~\lp(P_N, Q)\leq \epsilon }\\=& \textstyle\int \lambda f^\star\left(\frac{\ell(x, \xi)-\eta}{\lambda}\right) + r\lambda +\eta\, \d P_N^\star(\xi)
    \end{align*}
    for some distribution $P_N^\star\in \mc P(\Sigma)$, $\lp(P_N, P_N^\star)\leq \epsilon$ and $P^\star_N(\Sigma^\star(x))\geq \epsilon$ where $\Sigma^\star(x)\in \arg\max_{\xi\in \Sigma}\ell(x, \xi)$. Hence,
    \begin{align*}
      &\hat c^{\epsilon}_{f, d, r}(x, P_N)\\
      = & \inf \set{\int \lambda f^\star\left(\frac{\ell(x, \xi)-\eta}{\lambda}\right) + r\lambda +\eta\, \d P_N^\star(\xi)}{\eta\in\Re, ~\lambda\in\Re_+, ~\ell^\infty(x) - \eta \leq f^\infty \lambda}\\
      = & \hat c^{\epsilon}_{f, r}(x,P^\star_N)
    \end{align*}
    where the second equality follows from Theorem \ref{thm:love-dual}.
    We will now prove finally that $\hat c^{\epsilon}_{f, r}(x,P^\star_N) = \tilde c^{\epsilon}_{f, r}(x,P^\star_N)$. Recall that in (\ref{eq:dro:ef}) we defined
    \begin{equation}
      \label{eq:char:i}
      \tilde c^{\epsilon}_{f, r}(x,P^\star_N) = \left\{
        \begin{array}{rl}
          \sup & \int \ell(x, \xi) \, \d P_c(\xi)\\
          \st  & P_c\in \mc P_+(\Sigma), ~P_c\ll P_N^\star,\\
               & \int \d P_c(\xi) = 1,\\
               & \int f\left(\frac{\d P_c}{\d P_N^\star}(\xi )\right) \, \d P_N^\star(\xi)\leq r
        \end{array}
      \right.
    \end{equation}
    whereas following Theorem \ref{lemma:finite-representation-f-divergence} we have
    \begin{equation}
      \label{eq:char:ii}
      \hat c_{f, r}(x, P_N^\star) =
      \left\{
        \begin{array}{rl}
          \sup & \int \ell(x,\xi) \, \d P_c(\xi) + p_s \ell^\infty(x)\\
          \st  & P_c\in \mc P_+(\Sigma), ~P_c \ll P_N^\star, ~p_s \in \Re,\\
               & \int  \, \d P_c(\xi) + p_s=1,\\
               & \int f \left(\tfrac{\d P_c}{\d P_N^\star}(\xi)\right) \, \d P_N^\star(\xi) + p_s f^\infty \leq r.
        \end{array}
      \right.
    \end{equation}
    Hence, clearly we have that $\tilde c^{\epsilon}_{f, r}(x,P^\star_N)\leq \hat c_{f, r}(x, P_N^\star)$ for all $x\in X$. To prove that also  $\tilde c^{\epsilon}_{f, r}(x,P^\star_N)\geq \hat c_{f, r}(x, P_N^\star)$ for all $x\in X$ consider a feasible solution $(P_c, p_s)$ in \eqref{eq:char:ii}. Define the distribution $P^\star_{N, \Sigma^\star(x)}$ through $P^\star_{N, \Sigma^\star(x)}(B) = P^\star_{N}(B\cap \Sigma^\star(x))/P^\star_{N}(\Sigma^\star(x))$ for all measurable subsets $B$ of $\Sigma$.
    This distribution is well-defined as $P^\star_N(\Sigma^\star(x))\geq \epsilon>0$. Consider now $P_c'  = P_c + p_s P^\star_N(\Sigma^\star(x))$. Clearly, we have that by construction $P_c'\ll P^\star_N$ and $\int \d P_c'(\xi) = \int \d P_c(\xi) + p_s \int \d P^\star_{N, \Sigma^\star(x)}(\xi) = \int \d P_c(\xi) + p_s   = 1$. Lastly, we also have using Lemma \ref{lemma:recession-inequality} that
    \begin{align*}
      \int f\left(\frac{\d P'_c}{\d P_N^\star}(\xi )\right) \, \d P_N^\star(\xi) = & \int f\left(\frac{\d (P_c + p_s P^\star_N(\Sigma^\star(x)) ) }{\d P_N^\star}(\xi )\right) \, \d P_N^\star(\xi) \\
      \leq & \int f \left(\tfrac{\d P_c}{\d P_N^\star}(\xi)\right) \, \d P_N^\star(\xi) + \int p_s f^\infty\d P^\star_{N, \Sigma^\star(x)}(\xi) \\
      = & \int f \left(\tfrac{\d P_c}{\d P_N^\star}(\xi)\right) \, \d P_N^\star(\xi) + p_s f^\infty \leq r.
    \end{align*}
    Hence, as we have that $P_c'$ is feasible in \eqref{eq:char:i} and $\int \ell(x, \xi) \d P_c'(\xi) = \int \ell(x, \xi) \d P_c(\xi) + \ell^\infty(x)\int p_s \d P^\star_{N, \Sigma^\star(x)}(\xi) = \int \ell(x, \xi) \d P_c(\xi) +\ell^\infty(x)p_s$ it follows that $\tilde c^{\epsilon}_{f, r}(x,P^\star_N)\geq \hat c_{f, r}(x, P_N^\star)$ for all $x\in X$. It remains to show that $\hat c^{\epsilon}_{f, \lp, r}(x,P_N)=\tilde c^{\epsilon}_{f, r}(x,P^\star_N) = \tilde c^{\epsilon}_{f, \lp, r}(x,P_N)$. We first remark that we have trivially $\hat c^{\epsilon}_{f, \lp, r}(x,P_N)=\tilde c^{\epsilon}_{f, r}(x,P^\star_N) \geq \tilde c^{\epsilon}_{f, \lp, r}(x,P_N)$. On the other hand, as we have that $P^\star_N\in \mc P(\Sigma)$ and $\lp(P_N, P_N^\star)\leq \epsilon$ we must also have $\tilde c^{\epsilon}_{f, \lp, r}(x,P_N)\geq \tilde c^{\epsilon}_{f, r}(x,P^\star_N)=\hat c_{f, r}(x, P_N^\star)=\hat c^{\epsilon}_{f, d, r}(x, P_N)$. Hence, $\hat c^{\epsilon}_{f, \lp, r}(x,P_N)=\tilde c^{\epsilon}_{f, \lp, r}(x,P_N)$.
\end{proof}

\subsubsection{Proof of Theorem \ref{lemma:smoothing-dro-equivalence}}
\begin{proof}
  Observe indeed that we have
  \begingroup
  \allowdisplaybreaks
  \begin{align*}
    \hat c'_{f, W, r}(x, \hat P) & \defn \sup \set{\textstyle\int\ell(x, \xi)\, \d P(\xi)}{D'_{f, W}(\hat P, P)\leq r}\\
                            & = \left\{
                              \begin{array}{rl}
                                \sup & \int\ell(x, \xi)\, \d P(\xi)\\
                                \st  & P\in\mc P(\Sigma),~ Q\in\mc P(\Sigma),\\
                                     & D_{f}(Q, P)+D_W(\hat P, Q)\leq r.
                              \end{array}
                              \right.\\
                            & = \left\{
                              \begin{array}{rl}
                                \sup & \int\ell(x, \xi)\, \d P(\xi)\\
                                \st  & P\in\mc P(\Sigma), ~Q\in\mc P(\Sigma), ~\epsilon \in [0,r]\\
                                     & D_{f}(Q, P)+\epsilon\leq r\\
                                     & D_W(\hat P, Q)\leq \epsilon
                              \end{array}
                              \right.\\
                            & = \sup \set{\hat c^{\epsilon}_{f, W, r-\epsilon}(x, \hat P)}{\epsilon \in [0, r]}
  \end{align*}
  which establishes the claim.
  \endgroup
\end{proof}

\subsection{Proofs in Section \ref{sec:computation}}

\subsubsection{Proof of Lemma \ref{lemma:finite-representation-f-divergence}}
\label{proof:lemma:finite-representation-f-divergence}

The Lebesgue decomposition theorem states that for any $P\in \mc P(\Sigma)$ we have
\(
P = P_c + P_s
\)
for some unique $P_c\ll \hat P$ and $P_s\bot \hat P$. Hence, following (\ref{eq:f-dro}) we must have
\begin{equation}
  \label{eq:f-dro-continuous-singular-2}
  \begin{array}{rl}
    \hat c_{f,r}(x, \hat P) = \sup &  \int \ell(x, \xi)\, \d P_c(\xi)+ \int \ell(x, \xi)\, \d P_s(\xi) \\
    \st &  P_c\ll \hat P \in \mc P_+(\Sigma),~P_s\bot \hat P\in \mc P_+(\Sigma),\\
                                   & \int \d P_c(\xi) + \int \d P_s(\xi) = 1,\\
                                   & \int f\left(\tfrac{\d P_c}{\d \hat P}(\xi )\right) \, \d \hat P(\xi) + \int f^\infty \, \d P_s(\xi) \leq r.
  \end{array}
\end{equation}
Define now
\begin{equation}
  \label{eq:f-dro-continuous-singular-3}
  \begin{array}{rl}
    \hat c'_{f,r}(x, \hat P) = \sup & \int \ell(x, \xi)\, \d P_c(\xi)+ \int \ell(x, \xi)\, \d P_s(\xi) \\
    \st &  P_c\ll \hat P \in \mc M_+(\Sigma),~P_s\in \mc M_+(\Sigma),\\
                                    & \int \d P_c(\xi) + \int \d P_s(\xi) = 1,\\
                                    & \int f\left(\tfrac{\d P_c}{\d P'}(\xi )\right) \, \d P'(\xi) + \int f^\infty \, \d P_s(\xi) \leq r.
  \end{array}
\end{equation}
Clearly, we have that $\hat c'_{f,r}(x, \hat P)\geq \hat c_{f,r}(x, \hat P)$ for all $x\in X$, $\hat P\in \mc P(\Sigma)$ and $r\geq 0$ as $\hat c_{f,\epsilon,r}(x, \hat P)$ is characterized using an optimization problem which includes the additional constrains $P_c\ll \hat P$ and $P_s\bot \hat P$ compared to the optimization problem characterizing $\hat c'_{f,r}(x, \hat P)$ in (\ref{eq:f-dro-continuous-singular-2}). We will now prove that also $\hat c'_{f,r}(x, \hat P)\leq \hat c_{f,r}(x, \hat P)$. To do so we consider any feasible point $(P_c, P_s)$ in the optimization problem characterizing $\hat c'_{f,r}(x, \hat P)$. Again using the Lebesgue decomposition theorem we have
\( P_s = P_{cc} + P_{ss}\) for some $P_{cc}\ll \hat P$ and $P_{ss}\bot \hat P$. We shall now establish that the point $(P_c+P_{cc}, P_{ss})$ is feasible and attains the same objective value in the optimization problem characterizing $\hat c_{f,r}(x, \hat P)$ in (\ref{eq:f-dro-continuous-singular-2}). Indeed, we have following Lemma \ref{lemma:recession-inequality} that
\(
\int f\left(\tfrac{\d (P_c+P_{cc})}{\d \hat P}(\xi )\right) \, \d \hat P(\xi) + \int f^\infty \, \d P_{ss}(\xi)
\leq \int f\left(\tfrac{\d (P_c)}{\d \hat P}(\xi )\right) \, \d \hat P(\xi) + \int f^\infty \, \d (P_{cc}+P_{ss})(\xi) \leq r.
\)
The ultimate inequality follows from the feasibility of $(P_c, P_{ss}+P_{cc})=(P_c, P_s)$ in the optimization problem characterizing $\hat c'_{f,r}(\hat P)$.
Finally, 
\[
\textstyle\sup \set{ \int \ell(x, \xi)\, \d P_s(\xi)}{P_s\in \mc P_+(\Sigma), ~p_s\in \Re_+, ~p_s=\int_\Sigma \d P_s(\xi)} =p_s \sup_{\xi\in \Sigma}\ell(x, \xi) = p_s \ell^\infty(x)
\]
from which the claimed result follows immediately.

\subsubsection{Proof of Theorem \ref{thm:love-dual}}
\label{proof:thm:love-dual}

According to Lemma \ref{lemma:finite-representation-f-divergence} and the change of variables $(1-p_s) P'_c = P_c$ and $P'_c\in \mathcal P(\Sigma)$ we have that
\begin{align}
  \hat c_{f, r}(x, \hat P) = & \left\{
                               \begin{array}{rl}
                                 \sup & \int \ell(x,\xi) \, \d P_c(\xi) + p_s \ell^\infty(x)\\
                                 \st  & P_c\in \mc P_+(\Sigma), ~P_c \ll \hat P, ~p_s \in [0, 1],\\
                                      & \int  \, \d P_c(\xi) + p_s=1,\\
                                      & \int f \left(\tfrac{\d P_c}{\d \hat P}(\xi)\right) \, \d \hat P(\xi) + p_s f^\infty \leq r.
                               \end{array}
                               \right.\nonumber\\
  = & \left\{
      \begin{array}{rl}
        \sup & \int (1-p_s)\ell(x,\xi) \, \d P'_c(\xi) + p_s \ell^\infty(x)\\
        \st  & P'_c\in \mc P(\Sigma), ~P'_c \ll \hat P, ~p_s \in [0, 1],\\
             & \int f \left((1-p_s)\tfrac{\d P'_c}{\d \hat P}(\xi)\right) \, \d \hat P(\xi) \leq r-p_s f^\infty.
      \end{array}
      \right. \label{eq:continuous-f-divergence-reduction}
\end{align}
We remark that for a fixed $p_s$ the optimization problem over $P'_c$ in \eqref{eq:continuous-f-divergence-reduction} admits following \cite[Theorem 5.1]{ahmadi2012entropic} the strong dual representation stated in (\ref{eq:ahmadi-javid}) for all $x\in X$ and $\hat P\in \mc P(\Sigma)$. We hence obtain with the change of variables $(1-p_s) \eta = \lambda(\mu - r+p_s f^\infty)$ that
\begin{align*}
  & \hat c_{f, r}(x, \hat P)\\
  = & \sup_{p_s\in [0, 1]} \!\inf_{\lambda\in \Re_+, \mu\in \Re} \!\lambda\!\left[\mu\!+\! \int\!\! f^\star\left(\frac{(1-p_s)\ell(x, \xi)\!+\!\lambda(r- \mu -p_s f^\infty)}{(1-p_s)\lambda}\right) \d \hat P(\xi) \right] + p_s\ell^\infty(x)\\
  = & \sup_{p_s\in [0, 1]} \inf_{\lambda\in \Re_+} \lambda \left[\inf_{\mu\in \Re}\mu + \int f^\star\left(\frac{\ell(x, \xi)}{\lambda}-\frac{\mu - r+p_s f^\infty}{1-p_s}\right) \d \hat P(\xi) \right]+ p_s\ell^\infty(x)\\
  = & \sup_{p_s\in [0, 1]} \inf_{\lambda\in \Re_+, \eta\in \Re} (1-p_s)\eta +\lambda(r-p_sf^\infty) + \lambda \int f^\star\left(\frac{\ell(x, \xi)-\eta}{\lambda }\right) \d \hat P(\xi)+ p_s\ell^\infty(x) \\
  = & \sup_{p_s\in [0,1]} \inf_{\lambda\in\Re_+, \eta\in \Re} \eta +\lambda r + \lambda \int f^\star\left(\frac{\ell(x, \xi)-\eta}{\lambda }\right) \d \hat P(\xi)+ p_s(\ell^\infty(x)-\eta-f^\infty\lambda).
\end{align*}
As there always exists a Slater point $\bar \eta\in \Re$, $\bar \lambda\in \Re_+$ so that $\ell^\infty(x) - \bar \eta < f^\infty \bar \lambda$ we have finally following for instance  \cite[Proposition 5.3.1]{bertsekas2009convex} that 
\begin{align*}
  & \hat c_{f, r}(x, \hat P)\\
  = &
                               \inf_{\lambda\in\Re_+, \eta\in \Re} \sup_{p_s\in [0,1]} \eta +\lambda r + \lambda \int f^\star\left(\frac{\ell(x, \xi)-\eta}{\lambda }\right) \d \hat P(\xi)+ p_s(\ell^\infty(x)-\eta-f^\infty\lambda) \\
  = & \inf_{\lambda\in \Re_+, \eta\in\Re} \{L(\eta, \lambda)\!\defn\! \eta \!+\! \lambda r \!+\! \lambda\! \int\!\! f^\star\left(\frac{\ell(x, \xi)\!-\!\eta}{\lambda}\right) \d \hat P(\xi) \!+ \!\max(\ell^\infty(x)\! -\! \eta\! -\! f^\infty \lambda, 0)\}
\end{align*}
For any $\eta$ and $\lambda$ with $\ell^\infty(x) - \eta - f^\infty \lambda\geq 0$, we have that
\[
  L(\eta, \lambda) = \lambda r + \lambda \int f^\star\left(\frac{\ell(x, \xi)-\eta}{\lambda}\right) \d \hat P(\xi) + \ell^\infty(x) - f^\infty \lambda.
\]
As $f^\star$ is increasing we have that the previous function is decreasing in $\eta$. Hence, for any $(\eta, \lambda)$ with $\ell^\infty(x) - \eta - f^\infty \lambda\geq 0$ we have that $L(\eta, \lambda)\geq L(\eta', \lambda)$ with $\ell^\infty(x) - f^\infty \lambda= \eta'\geq \eta$.
Hence it follows that
\begin{align*}
  \hat c_{f, r}(x, \hat P) = & \inf_{\lambda\in \Re_+, \eta\in\Re} \set{\!\eta + \lambda r + \lambda \int\!\! f^\star\left(\frac{\ell(x, \xi)-\eta}{\lambda}\right) \d \hat P(\xi)\!\!}{\!\!\ell^\infty(x) - \eta - f^\infty \lambda\leq 0\!}.
\end{align*}

It remains to show that as $r>0$ we have that the minimum in Equation \eqref{eq:love-dual}
is indeed attained. It is trivial to see that its objective function $g(\lambda, \eta)\defn \int \lambda f^\star\left(\frac{\ell(x,\xi)-\eta}{\lambda}\right) \d \hat P(\xi) + r\lambda +\eta$ is convex and lower semicontinuous.
Hence, the statement would follow if we can show that $g(\lambda, \eta)$ is coercive by \citep[Proposition 3.2.1]{bertsekas2009convex}.
Let us define $g_0(\lambda, \eta) \defn \lambda f^\star\left(\tfrac{-\eta}{\lambda}\right) + r\lambda +\eta$ and observe that $g_0(\lambda, \eta) \leq g(\lambda, \eta)$ for all $x\in X$ and $\xi\in \Sigma$ as $0\leq \ell(x, \xi)$ for all $x\in X$ and $\xi\in \Sigma$ and as $f^\star$ is nondecreasing. It hence suffices to show that $g_0(\lambda, \eta)$ is coercive. Observe that
\begin{align*}
g_{0}(\lambda,\eta)  & =\sup_{t\geq0}(-t\eta-\lambda f(t))+r\lambda+\eta\\
& \geq\max\left\{  -\eta-\lambda f(1),-\tfrac{\eta}{2}-\lambda f\left(\tfrac{1}{2}\right)  ,-2\eta-\lambda f(2)\right\}  +r\lambda+\eta\\
& =\max\left\{  r\lambda,\frac{1}{2}\eta+\left(  r-f\left(  \tfrac{1}{2}\right)  \right)  \lambda,-\eta+(r-f(2))\lambda\right\}  ,
\end{align*}
which implies $g_{0}(\lambda,\eta)$ is coercive (if $\lambda\rightarrow\infty$, then $g_{0}(\lambda,\eta)\geq r\lambda\rightarrow\infty$; if $\lambda$ is bounded and $|\eta|\rightarrow\infty$, we also have $g_{0}(\lambda,\eta)\geq\max\left\{  \frac{1}{2}\eta+\left(  r-f\left(  \frac{1}{2}\right)\right)  \lambda,-\eta+(r-f(2))\lambda\right\}  \rightarrow\infty$).

\subsubsection{Proof of Lemma \ref{lemma:moreau-yosida}}
\label{proof:lemma:moreau-yosida}

\begin{proof}
  We have
  \begin{align*}
    \ell_\gamma(x, \xi) \defn & \sup_{\xi'\in \Sigma} \ell(x, \xi')-\gamma d(\xi, \xi')\\
    = & \sup_{\delta \in \Delta}\sup_{\xi'\in \Sigma, d(\xi, \xi')=\delta} \ell(x, \xi')-\gamma d(\xi, \xi')\\
    = & \sup_{\delta \in \Delta}\sup_{\xi'\in \Sigma, d(\xi, \xi')=\delta} \ell(x, \xi')-\gamma \delta\\
    = & \sup_{\delta \in \Delta}\sup_{\xi'\in \Sigma, d(\xi, \xi')\leq \delta} \ell(x, \xi')-\gamma \delta\\
    = & \sup_{\delta\in \Delta} \ell^\delta(x, \xi) -\gamma\delta.
  \end{align*}
  Here the first equality follows from $\Sigma = \cup_{\delta\in \Delta}\set{\xi'\in\Sigma}{\xi'\in \Sigma, \, d(\xi, \xi')=\delta}$ as $\Sigma$ is compact with diameter $\diam(\Sigma)$. The third equality follows from the fact that we have
  \(
    \ell_\gamma(x, \xi) = \sup_{\delta \in \Delta}\sup_{\xi'\in \Sigma, d(\xi, \xi')=\delta} \ell(x, \xi')-\gamma \delta
    \leq \sup_{\delta \in \Delta}\sup_{\xi'\in \Sigma, d(\xi, \xi')\leq\delta} \ell(x, \xi')-\gamma \delta =: \tilde \ell_\gamma(x, \xi)
  \)
  but also
  \(
    \ell_\gamma(x, \xi) \geq \tilde \ell_\gamma(x, \xi)
  \)
  for all $x\in X$ and $\xi\in\Sigma$ as we will show now. Indeed, let $\{(\delta_k\in \Delta, \xi_k'\in \Sigma)\}$ be such $d(\xi, \xi'_k)\leq\delta_k$ for all $k\geq 1$ and $\lim_{k\to\infty} \ell(x, \xi_k')-\gamma \delta_k = \tilde \ell_\gamma(x, \xi)$ for some arbitrary $x\in X$ and $\xi\in\Sigma$. Consider now the auxiliary sequence $\{(\delta'_k=d(\xi, \xi'_k)\leq \delta_k, \xi_k'\in \Sigma)\}$ then clearly we have $d(\xi, \xi'_k)=\delta'_k$ for all $k\geq 1$ and hence $\ell_\gamma(x, \xi) \geq \lim_{k\to\infty} \ell(x, \xi_k')-\gamma \delta'_k \geq \lim_{k\to\infty} \ell(x, \xi_k')-\gamma \delta_k = \tilde \ell_\gamma(x, \xi)$ as $\gamma\geq 0$.
\end{proof}

\subsubsection{Proof of Proposition \ref{lemma:wasserstein:approximation}}
\label{proof:lemma:wasserstein:approximation}

Introduce the auxiliary distance function
\(
d^K(\xi, \xi') \defn \min \set{\delta\in\Delta_K}{d(\xi, \xi')\leq \delta}
\)
for all $\xi\in \Sigma$ and $\xi'\in\Sigma$.
From the construction of the auxiliary distance function it is clear that if $d(\xi, \xi')\leq \epsilon$ we have sandwich inequalities $d(\xi, \xi') \leq d^K(\xi, \xi') \leq d(\xi, \xi')+\tfrac{\epsilon}{K}$ for any $\xi'\in \Sigma$ and $\xi\in \Sigma$. Likewise, if $d(\xi, \xi')\leq \epsilon$ we have that $d(\xi, \xi') \leq d^K(\xi, \xi') \leq \left(\tfrac{\diam(\Sigma)}{\epsilon}\right)^{1/K}  d(\xi, \xi')$ for any $\xi'\in \Sigma$ and $\xi\in \Sigma$. Hence, we have
\[
  d(\xi, \xi') \leq d^K(\xi, \xi') \leq \left(\tfrac{\diam(\Sigma)}{\epsilon}\right)^{1/K}  d(\xi, \xi') +\tfrac{\epsilon}{K}
\]
for any $\xi'\in \Sigma$ and $\xi\in \Sigma$.
Let
\begin{align}
  W^K_d(P', P) = & \inf \set{\int  d^K(\xi', \xi) \, \d \gamma(\xi', \xi)}{\gamma \in \Gamma(P', P)} \nonumber\\
  \leq & \inf \set{\int \left(\tfrac{\diam(\Sigma)}{\epsilon}\right)^{1/K} d(\xi', \xi) + \frac{\epsilon}{K} \, \d \gamma(\xi', \xi)}{\gamma \in \Gamma(P', P)}\nonumber\\
  \label{eq:Wasserstein-approximation-inequality}
  \leq &  \left(\tfrac{\diam(\Sigma)}{\epsilon}\right)^{1/K} W_d(P', P) + \frac{\epsilon}{K}
\end{align}
denote the associated auxiliary optimal transport distance.

Let
\(
\hat c_{d^K, \epsilon}(x, \hat P)
\)
denote the predictor in (\ref{eq:W-dro}) associated with the auxiliary optimal transport distance $W^K_{d}$. Remark that by construction the auxiliary distance function $d^K$ is lower semicontinous as $d$ is lower semicontinous and we have also $d(\xi,\xi)=0$ for all $\xi\in \Sigma$.
From Theorem \ref{thm:kuhn-dual} it follows hence that
\begin{align*}
  \hat c_{d^K, \epsilon}(x, \hat P) = & \min \set{\int \left[\sup_{\delta\in \Delta} \sup_{\xi'\in \Sigma, d^K(\xi, \xi')\leq \delta}\ell(x,\xi') -\gamma \delta\right] \,\d \hat P(\xi)+\epsilon\gamma}{\gamma\geq 0}\\
  = & \min \set{\int \left[\sup_{\delta\in \Delta_K} \sup_{\xi'\in \Sigma, d^K(\xi, \xi')\leq \delta}\ell(x,\xi') -\gamma \delta\right] \,\d \hat P(\xi)+\epsilon\gamma}{\gamma\geq 0} \\
  = & \hat c^K_{d, \epsilon}(x, \hat P)
\end{align*}
for all $x\in X$ and $\hat P\in \mc P(\Sigma)$.

We finally show that
\[
  \hat c_{d^K, \epsilon}(x, \hat P) \geq \hat c_{d, \epsilon}(x, \hat P)\left(\left(\tfrac{\diam(\Sigma)}{\epsilon}\right)^{1/K}+\tfrac{1}{K}\right)^{-1}
\]
for all $x\in X$ and $\hat P\in \mc P(\Sigma)$ establishing the result. Consider indeed any arbitrary $P\in \mc P(\Sigma)$ so $W_d(\hat P, P)\leq \epsilon$. From the convexity of $W_d^K$ we have $W_d^K(\hat P, (1-\lambda) \hat P + \lambda P) \leq  (1-\lambda) W_d^K(\hat P, \hat P) + \lambda W_d^K(\hat P, P) = \lambda W_d^K(\hat P, P)\leq \lambda(\left(\tfrac{\diam(\Sigma)}{\epsilon}\right)^{1/K} W_d(\hat P, P) +\tfrac{\epsilon}{K}) \leq \lambda(\left(\tfrac{\diam(\Sigma)}{\epsilon}\right)^{1/K} \epsilon + \tfrac{\epsilon}{K})\leq \epsilon$ for $$0\leq \lambda \leq \left(\left(\tfrac{\diam(\Sigma)}{\epsilon}\right)^{1/K}+\tfrac{1}{K}\right)^{-1}\leq 1.$$
We have hence
\begin{align*}
  \hat c_{d^K, \epsilon}(x, \hat P) & \geq \int \ell(x,\xi) \, \d \left[(1-\lambda)\hat P + \lambda P\right](\xi) \\
  \geq &(1-\lambda) \int \ell(x,\xi) \, \d \hat P + \lambda \int \ell(x,\xi) \,  \d P (\xi)\geq \lambda \int \ell(x,\xi) \,  \d P (\xi).
\end{align*}
Considering now a sequence of points $P_k$ feasible in Equation (\ref{eq:W-dro}) with $$\lim_{k\to\infty} \int \ell(x,\xi) \,  \d P_k (\xi) = \hat c_{d, \epsilon}(x, \hat P)$$ establishes the desired result.

\subsubsection{Proof of Theorem \ref{thm:main-type-I}}
\label{proof:main-type-I}

\begin{proof}
  We clearly have the decomposition
  \[
    \hat c^{\epsilon}_{f, d,  r}(x, \hat P) = \sup \set{\hat c_{f, r}(x, Q)}{Q\in \mc P(\Sigma),~W_d(\hat P, Q)\leq \epsilon}
  \]
  where $\hat c_{f, r}(x, Q)$ is defined in (\ref{eq:f-dro}). Following Theorem \ref{thm:love-dual} we have thus
  \begin{align*}
    &\hat c^{\epsilon}_{f, d, r}(x, \hat P)\\
    = & \sup\Big\{\!\inf \set{\textstyle\int \lambda f^\star\left(\frac{\ell(x, \xi)-\eta}{\lambda}\right) \d Q(\xi) + r\lambda +\eta\!\!}{\!\!\eta\in\Re, \lambda\in\Re_+, \ell^\infty(x) - \eta \leq f^\infty \lambda} : \\
    & \hspace{24em}  Q\in \mc P(\Sigma),~W_d(\hat P, Q)\leq \epsilon \Big\}.
  \end{align*}
  Consider the saddle point function
  \(
    L(\lambda, \eta; Q)=\int \lambda f^\star\left(\tfrac{(\ell(x, \xi)-\eta)}{\lambda}\right) + r\lambda +\eta\, \d Q(\xi).
  \)
  The function $L(\lambda, \eta; Q)$ is clearly convex and lower semicontinuous in $(\lambda, \eta)$ and linear (concave) in $Q$.
  As showed in the proof of Theorem \ref{thm:love-dual} we can restrict the dual variables to a compact set and hence  a standard Minimax theorem \citep[Theorem 4.2]{sion1958general} implies that
  \begin{align}
    \label{eq:minimax-equality}
    &\hat c^{\epsilon}_{f, d, r}(x, \hat P)\\
    = & \inf \Big\{\sup \set{\textstyle\int \lambda f^\star\left(\frac{\ell(x, \xi)-\eta}{\lambda}\right) + r\lambda +\eta\, \d Q(\xi) }{ Q\in \mc P(\Sigma),~W_d(\hat P, Q)\leq \epsilon } : \nonumber\\
                                          & \hspace{20.5em} \eta\in\Re, ~\lambda\in\Re_+, ~\ell^\infty(x) - \eta \leq f^\infty \lambda \Big\}.\nonumber
  \end{align}
  Let $\tilde \ell(x, \xi) \defn \lambda f^\star\left(\tfrac{(\ell(x, \xi)-\eta)}{\lambda}\right) + r\lambda +\eta$.
  We clearly have that
  \begin{align*}
    \tilde \ell_\gamma(x, \xi) = & r\lambda +\eta +\sup_{\delta\in \Delta}\sup_{\xi'\in \Sigma,\, d(\xi, \xi')\leq \delta }\lambda f^\star\left(\tfrac{(\ell(x, \xi')-\eta)}{\lambda}\right) -\gamma \delta\\
    = & r\lambda +\eta +\sup_{\delta\in \Delta}\lambda f^\star\left(\tfrac{(\textstyle\sup_{\xi'\in \Sigma,\, d(\xi, \xi')\leq \delta } \ell(x, \xi')-\eta)}{\lambda}\right) -\gamma \delta\\
    = & r\lambda +\eta +\sup_{\delta\in \Delta}\lambda f^\star\left(\tfrac{(\ell^\delta(x, \xi)-\eta)}{\lambda}\right) -\gamma \delta.
  \end{align*}
    We now show the function $\xi \mapsto \tilde \ell(x, \xi)$ is upper semicontinuous on $\Sigma$ for every $x\in X$ for any $(\eta, \lambda)\in \Re\times\Re_+$ so that $\ell^\infty(x) - \eta \leq f^\infty \lambda$. Observe that we have here indeed $\tfrac{(\ell(x, \xi)-\eta)}{\lambda}\leq \tfrac{(\ell^\infty(x)-\eta)}{\lambda}\leq f^\infty$ for all $\xi\in \Xi$.
    Hence, consider any sequence $\xi_k\in \Xi$ converging to some $\xi_\infty\in \Xi$. We have
    \begin{align*}
      \lim_{k\to\infty} \tilde \ell(x, \xi_k) = \lambda f^\star\left(\lim_{k\to\infty} \tfrac{(\ell(x, \xi_k)-\eta)}{\lambda}\right) + r\lambda +\eta 
    \end{align*}
    as $f^\star$ is continuous on $[-\infty, f^\infty]$; see Lemma \ref{sec:supporting_lemmas}. As we have that the function $\xi \mapsto\ell(x, \xi)$ is itself upper semicontinous and because the function $f^\star$ is nondecreasing, see again Lemma \ref{sec:supporting_lemmas}, we can write
    \begin{align*}
      \lim_{k\to\infty} \tilde \ell(x, \xi_k) \leq \lambda f^\star\left(\tfrac{(\ell(x, \xi_\infty)-\eta)}{\lambda}\right) + r\lambda +\eta = \tilde \ell(x, \xi_\infty).
    \end{align*}
  As we have shown that $\tilde \ell$ is an upper semicontinuous function, we can apply Theorem \ref{thm:kuhn-dual} and  get
  \begin{align*}
    &\hat c^{\epsilon}_{f, d, r}(x, \hat P)\\
    = & \inf \Big\{\inf \set{\textstyle r \lambda +\eta +\sup_{\delta\in \Delta}\lambda f^\star\left(\tfrac{(\ell^\delta(x, \xi)-\eta)}{\lambda}\right) -\gamma \delta\, \d \hat P(\xi) + \gamma\epsilon}{\gamma\geq 0 } : \\
                                          & \hspace{20.5em} \eta\in\Re, ~\lambda\in\Re_+, ~\ell^\infty(x) - \eta \leq f^\infty \lambda \Big\}
  \end{align*}
  from which the claimed result follows immediately.
\end{proof}

\subsubsection{Proof of Proposition \ref{lemma:wasserstein-divergence:approximation}}
\label{proof:wasserstein-divergence:approximation}

Let
\(
\hat c'^{\epsilon, K}_{f, d, r}(x, P)
\)
denote the predictor in (\ref{eq:formulation-type-1}) associated with the auxiliary distance function $d^K$ and auxiliary Wasserstein distance $W^K_{d}$ defined in the proof of Proposition \ref{lemma:wasserstein:approximation}. Remark that by construction the auxiliary distance function $d^K$ is lower semicontinous as $d$ is lower semicontinous and we have also $d(\xi,\xi)=0$ for all $\xi\in \Sigma$.
From Theorem \ref{thm:main-type-I} it follows hence that
\begin{align*}
  &\hat c^{\epsilon}_{f, d^K, r}(x, \hat P)\\
  = &\! \left\{\!\!\!
                                            \begin{array}{r@{~}l}
                                              \inf & \int\! \left[\sup_{\delta\in \Delta}\! \lambda f^\star \left(\tfrac{(\sup_{\xi'\in \Sigma, d^K(\xi, \xi')\leq \delta}\ell(x,\xi')-\eta)}{\lambda}\right)-\delta \gamma\right] \d \hat P(\xi) \!+\! r \lambda \!+\!\eta \!+\!\epsilon\gamma\\
                                              \st & \eta\in \Re, ~ \lambda\in \Re_+, ~\gamma\in \Re_+, \\
                                                   & \ell^\infty(x)-\eta \leq f^\infty \lambda.
                                            \end{array}\right.
  \\
  = &\! \left\{\!\!\!
      \begin{array}{r@{~}l}
        \inf & \int \!\left[\sup_{\delta\in \Delta_K} \!\lambda f^\star \left(\tfrac{(\sup_{\xi'\in \Sigma, d^K(\xi, \xi')\leq \delta}\ell(x,\xi')-\eta)}{\lambda}\right)-\delta \gamma\right]  \d \hat P(\xi) \!+\! r \lambda \!+\!\eta \!+\!\epsilon\gamma\\
        \st & \eta\in \Re, ~ \lambda\in \Re_+, ~\gamma\in \Re_+, \\
             & \ell^\infty(x)-\eta \leq f^\infty \lambda.
      \end{array}\right. \\
  = & \hat c^{\epsilon, K}_{f, d, r}(x, \hat P)
\end{align*}
for all $x\in X$ and $\hat P\in \mc P(\Sigma)$.

We finally show that
\[
  \hat c^{\epsilon}_{f, d^K, r}(x, \hat P) \geq \hat c^\epsilon_{f, d, r}(x, \hat P)\left(\left(\tfrac{\diam(\Sigma)}{\epsilon}\right)^{1/K}+\tfrac{1}{K}\right)^{-1}
\]
for all $x\in X$ and $\hat P\in \mc P(\Sigma)$ establishing the result. Consider indeed any arbitrary $P\in \mc P(\Sigma)$ so $W_d(\hat P, P)\leq \epsilon$. From the convexity of $W_d^K$ we have $W_d^K(\hat P, (1-\lambda) \hat P + \lambda P) \leq  (1-\lambda) W_d^K(\hat P, \hat P) + \lambda W_d^K(\hat P, P) = \lambda W_d^K(\hat P, P)\leq \lambda(\left(\tfrac{\diam(\Sigma)}{\epsilon}\right)^{1/K} W_d(\hat P, P) +\tfrac{\epsilon}{K}) \leq \lambda(\left(\tfrac{\diam(\Sigma)}{\epsilon}\right)^{1/K} \epsilon + \tfrac{\epsilon}{K})\leq \epsilon$ for $0\leq \lambda \leq \left(\left(\tfrac{\diam(\Sigma)}{\epsilon}\right)^{1/K}+\tfrac{1}{K}\right)^{-1}\leq 1$ and $\hat P\in\mc P(\Sigma)$.
Remark that we have
\[
  \hat c^{\epsilon}_{f, d,  r}(x, \hat P) = \sup \set{\hat c_{f, r}(x, Q)}{Q\in \mc P(\Sigma),~W_d(\hat P, Q)\leq \epsilon}
\]
where $\hat c_{f, r}(x, Q)$ is defined in (\ref{eq:f-dro}).
From concavity of $\hat c_{f, r}(x, Q)$ in $Q\in \mc P(\Sigma)$ we have hence
\[
  \hat c^{\epsilon}_{f, d^K,  r}(x, \hat P) \geq \hat c_{f, r}(x, (1-\lambda)\hat P + \lambda P) \geq (1-\lambda) \hat c_{f, r}(x, \hat P) + \lambda \hat c_{f, r}(x, P) \geq \lambda \hat c_{f, r}(x, P).
\]
Considering now a sequence of points $P_k$ feasible in (\ref{eq:W-dro}) so that $\lim_{k\to\infty} \hat c_{f, r}(x, P_k) = \hat c^{\epsilon}_{f, d,  r}(x, \hat P)$ establishes the desired result.

\subsubsection{Proof of Theorem \ref{thm:main-type-II}}

\begin{proof}
  We have the decomposition
  \[
    c^{\epsilon}_{f, d, r}(x, \hat P) = \sup \set{\hat c_{d, \epsilon}(x, Q)}{Q\in \mc P(\Sigma),~D_f(\hat P, Q)\leq r}
  \]
  where $\hat c_{d, \epsilon}(x, Q)$ is defined in (\ref{eq:W-dro}). Following Theorem \ref{thm:kuhn-dual} and $\hat c^{ \epsilon}_{f, d, r}(x, \hat P)\leq \ell^\infty(x)$ for all $x\in X$ we have thus
  \begin{align}
    \label{eq:primal-dual-characterization}
    \hat c^{ \epsilon}_{f, d, r}(x, \hat P) = & \sup \big\{\inf \set{\textstyle\int \left[\sup_{\delta\in \Delta}\ell^{\delta}(x,\xi)-\gamma\delta\right] \d Q(\xi) + \epsilon\gamma}{\gamma\in [0, \frac{\ell^{\infty}(x)}{\epsilon}]}\! :\\
    & \hspace{19em}\! Q\in \mc P(\Sigma),~D_f(\hat P, Q)\leq r \big\}.\nonumber
  \end{align}
  Consider the saddle point function
  \(
    L(\gamma; Q)=\int \left[\sup_{\delta\in \Delta}\ell^{\delta}(x,\xi)-\gamma\delta\right] + \epsilon\gamma \, \d Q(\xi) .
  \)
  The saddle point function $L(\gamma; Q)$ is clearly convex and lower semicontinuous in $\gamma$ and concave (linear) in $Q$.
  As $[0, \ell^\infty(x)/\epsilon]$ is a compact set a standard Minimax theorem \citep[Theorem 4.2]{sion1958general} implies that
  \begin{align*}
    & \hat c^\epsilon_{f, d, r}(\ell, \hat P)\\
    = & \inf \set{\sup \set{\textstyle\int \left[\sup_{\delta\in \Delta}\ell^{\delta}(x,\xi)\!-\!\gamma\delta\right]\!+\!  \epsilon\gamma\, \d Q(\xi) \!\!}{\!\!Q\in \mc P(\Sigma),~D_f(\hat P, Q)\leq \epsilon} \!\!}{\!\!\gamma\geq 0}\!.
  \end{align*}
  Applying Theorem \ref{thm:love-dual} we hence get that
  \begin{align*}
    \hat c^\epsilon_{f, d, r}(x, \!\hat P) \!=\! & \left\{\!\!\!
        \begin{array}{r@{\,}l}
          \inf & \int \lambda f^\star\left(\frac{\left[\sup_{\delta\in \Delta}\ell^{\delta}(x,\xi)-\gamma\delta\right]+ \epsilon\gamma - \eta}{\lambda} \right) \, \d \hat P(\xi) + r\lambda+\eta \\[0.75em]
          \st & \gamma\!\in\!\Re_+, \eta\!\in \!\Re, \lambda\!\in\! \Re_+ , \sup_{\xi\in\Sigma}\!\left[\sup_{\delta\geq 0}\ell^{\delta}(x,\xi)\!-\!\gamma\delta\right]\!+\!\epsilon\gamma \!-\!\eta \!\leq\! f^\infty \lambda.
        \end{array}
        \right.
  \end{align*}
  Finally, we observe that
  \begin{align*}
    \sup_{\xi\in\Sigma}\left[\sup_{\delta\in \Delta}\ell^{\delta}(x,\xi)-\gamma\delta\right] \defn &  \sup_{\xi\in\Sigma, \, \xi'\in \Sigma}\ell(x, \xi') - \gamma d(\xi, \xi') \\
    = &  \sup_{\delta\in \Delta} \sup_{\xi\in\Sigma, \, \xi'\in\Sigma, \, d(\xi, \xi')\leq \delta} \ell(x, \xi') - \gamma d(\xi, \xi')  \\
    = &  \sup_{\delta\in \Delta} \sup_{\xi\in\Sigma, \, \xi'\in\Sigma, \, d(\xi, \xi')\leq \delta} \ell(x, \xi') - \gamma \delta \\
    = & \sup_{\delta\in \Delta} \, \ell^\infty(x) - \gamma\delta = \ell^\infty(x)
  \end{align*}
  from which the claim now follows. The first equality follows from $\Sigma\times \Sigma =\set{\xi\in\Sigma, \, \xi'\in\Sigma}{\exists \delta\geq 0, \, d(\xi, \xi')\leq \delta}$ and $\gamma\geq 0$.
  Clearly we have the upper bound $\sup_{\xi\in\Sigma, \, \xi'\in\Sigma, \, d(\xi, \xi')\leq \delta} \ell(x, \xi') \leq \ell^\infty(x)$ for any $\delta\in \Delta$. On the other hand, we have that $$\sup_{\xi\in\Sigma, \, \xi'\in\Sigma, \, d(\xi, \xi')\leq \delta} \ell(x, \xi') \geq \sup_{\xi\in\Sigma} \ell(x, \xi)=\ell^\infty(x)$$ as $0=d(\xi, \xi)\leq \delta$ for all $\xi\in \Sigma$. Hence, the penultimate equality follows from $\sup_{\xi\in\Sigma, \, \xi'\in\Sigma, \, d(\xi, \xi')\leq \delta} \ell(x, \xi') = \ell^\infty(x)$. The final equality follows from $\gamma\geq 0$.
\end{proof}

\subsubsection{Proof of Proposition \ref{lemma:divergence-wasserstein:approximation}}

\begin{proof}
  Clearly, the inequality $\hat c^{\epsilon, K}_{d, f, r}(x, \hat P)\leq  \hat c^{\epsilon}_{d, f, r}(x, \hat P)$ for all $x\in X$ and $P\in \hat P(\Sigma)$ immediately follows from $\Delta_N\subseteq \Delta$ and the fact that $f^\star$ is a nondecreasing function. On the other hand, the inequality $$\hat c^{\epsilon}_{d, f, r}(x, \hat P)\leq \left(\left(\tfrac{\diam(\Sigma)}{\epsilon}\right)^{1/K}+\tfrac{1}{K}\right) \hat c^{\epsilon, K}_{d, f, r}(x, \hat P)$$ can be shown following the proof of Theorem \ref{thm:main-type-II} verbatim where using Proposition \ref{lemma:wasserstein:approximation} we replace equality (\ref{eq:primal-dual-characterization}) with
  \begin{align*}
    & \hat c^{ \epsilon}_{d,f, r}(x, \hat P)\\
    \leq & \left(\left(\tfrac{\diam(\Sigma)}{\epsilon}\right)^{\frac{1}{K}}+\frac{1}{K}\right)  \sup \big\{\inf \set{\textstyle\int \left[\sup_{\delta\in \Delta_K}\ell^{\delta}(x,\xi)-\gamma\delta\right] \d Q(\xi) + \epsilon\gamma\!\!}{\!\!\gamma\geq 0}:\\
    & \hspace{24.5em} Q\in \mc P(\Sigma),~D_f(\hat P, Q)\leq r \big\}
  \end{align*}
  instead.
\end{proof}

\clearpage

\section{Tables}

\begin{table}[ht]
  \centering
  \caption{Several common $f$-divergences.  }
  \label{table:divergence_functions}
  \begin{tblr}{
      hline{1,9} = {-}{0.08em},
      hline{2} = {-}{},
    }
    Name                        & Symbol       & Function $f(t)$                & Function $f^\star(s)$ \\
    Kullback-Leibler divergence & $D_{\kl}$    & $-\log(t)+t-1$                  & $-\log(1-s)$         \\
    Burg entropy & $D_{\lk}$    & $t \log(t)-t+1$                 & $\exp(s)-1$            \\
    Pearson divergence          & $D_{\chi^2}$ & $t(1-1/t)^2$                    & $2-2\sqrt{1-s}$                       \\
    Total variation             & $TV$         & $|t-1|/2$                       & $\left\{\!\!\!\!\begin{array}{r@{\,}l}
                                                                                     -\frac 12 & {\rm{if}}\, s< -\frac 12\\[0.3em]
                                                                                     s & {\rm{if}} \, s \, \!\in\! [-\frac 12,\! \frac 12]
                                                                                   \end{array}\right.$\\
    Squared Hellinger distance  & $D_{\mathrm{H}}$        & $(1-\sqrt{t})^2$                & $s/(1-s)$                       \\
    Le Cam distance             & $D_{\mathrm{LC}}$     & $(1-t)/(t+1)$                   & $1-s-2\sqrt{-2s}$                       \\
    Jensen-Shannon divergence   & $D_{\mathrm{JS}}$     & $t\log(\frac{2t}{t+1})\!+\!\log(\frac{2}{t+1})$ & $-\log(2-\exp(s))$                       
  \end{tblr}
\end{table}

\end{document}